\documentclass[reqno,10pt]{amsart}
\UseRawInputEncoding
\pdfoutput=1
\usepackage{amsfonts}
\usepackage{graphicx}
\usepackage{amsfonts,amsmath, amssymb}
\usepackage{bbm,mathrsfs,mathscinet}
\usepackage{hyperref}
\usepackage{marginnote}
\usepackage{color}
\usepackage{setspace}

\numberwithin{equation}{section}

\newcommand{\fu}{\mathfrak{u}}
\newcommand{\fb}{\mathfrak{b}}
\newcommand{\fa}{\mathfrak{a}}
\newcommand{\fw}{\mathfrak{w}}
\newcommand{\R}{\mathbb{R}}

\newcommand{\tw}{\tilde{w}}

\newcommand{\Rmnum}[1]{\expandafter\@slowromancap\romannumeral #1@}

\newtheorem{theorem}{Theorem}[section]

\newtheorem{lemma}[theorem]{Lemma}
\newtheorem{proposition}[theorem]{Proposition}
\newtheorem{remark}[theorem]{Remark}
\newtheorem{definition}[theorem]{Definition}
\allowdisplaybreaks

\def\v{\varepsilon}

\def\k{\kappa}

\def\g{\gamma}
\def\d{\delta}

\def\r{\rho}
\def\s{\sigma}

\def\f{\frac}
\def\pa{\partial}

\def\hep{h^{\varepsilon}}

\def\FRep{F^{\varepsilon}_{R}}
\def\pt{\partial_{t}}
\def\nax{\nabla_{x}}
\def\dis{\displaystyle}

\def\intr{\int_{\mathbb{R}^{3}}}
\def\ints{\int_{S^{2}}}
\def\smum{\sqrt{\mu_{M}}}
\def\sp{\shortparallel}

\newcommand{\FL}{\mathbf{L}}
\newcommand{\FP}{\mathbf{P}}
\newcommand{\FI}{\mathbf{I}}
\newcommand{\Fi}{\mathbf{1}}

\begin{document}
	
	\title[Hilbert expansion of soft Boltzmann equation]{Hilbert expansion  of Boltzmann equation with soft potentials and specular boundary condition in  half-space}

	\author[Jing Ouyang]{Jing Ouyang}
	\address[Jing Ouyang]{Academy of Mathematics and Systems Science, Chinese Academy of Sciences, Beijing 100190, China; School of Mathematical Sciences, University of Chinese Academy of Sciences, Beijing 100049, China
	}
	\email{ouyangjing@amss.ac.cn}
	
	\author[Yong Wang]{Yong Wang}
	\address[Yong Wang]{Academy of Mathematics and Systems Science, Chinese Academy of Sciences, Beijing 100190, China; School of Mathematical Sciences, University of Chinese Academy of Sciences, Beijing 100049, China
	}
	\email{yongwang@amss.ac.cn}

\begin{abstract}
Boundary effects play an important role in the study of hydrodynamic limits in the Boltzmann theory.
We justify rigorously the validity of the hydrodynamic limit from the Boltzmann equation of soft potentials to the compressible Euler equations by the  Hilbert expansion with multi-scales. Specifically,  the Boltzmann solutions are  expanded  into three parts: interior part, viscous boundary layer and Knudsen boundary layer. Due to the weak effect of collision
frequency of soft potentials, new difficulty arises when tackling the existence of Knudsen layer solutions with space decay rate,  which has been overcome under some constraint conditions and losing velocity weight arguments.
\end{abstract}

	\keywords{Boltzmann equation, compressible Euler equations, hydrodynamic limit, Hilbert expansion, viscous boundary layer, Knudsen boundary layer}
	\date{\today}
	\maketitle
	
	\setcounter{tocdepth}{2}
	\tableofcontents
	
	\thispagestyle{empty}
	
\section{Introduction and Main Results}
	
\subsection{Introduction}
It is well-known that the Boltzmann equation is closely related to the fluid dynamical systems for both compressible and incompressible flows
since the founding work of Maxwell \cite{Maxwell} and Boltzmann \cite{Boltzmann}. In 1912, Hilbert proposed a systematic formal asymptotic expansion for Boltzmann equation with respect to Knudsen number $\mathscr{K}_n\ll1$. 
In 1916 and 1917, Enskog and Chapman independently proposed a different formal expansion, respectively. Based on Hilbert or Chapman-Enskog expansions, the standard fluid theory can be derived formally, for instance: the compressible Euler and Navier-Stokes equations, the incompressible Euler and Navier-Stokes (Fourier) equations, {\it et. al}.

In the past decades,  great effort has been devoted to the study of the hydrodynamic limit from the Boltzmann equation to the fluid systems. When the solutions of compressible Euler equations are smooth, Caflisch \cite{Caflish} rigorously justified the hydrodynamic limit of Boltzmann equation to the compressible Euler equations by a the truncated Hilbert expansion, see also \cite{Grad,Lachowicz,Nishida,Ukai-Asano}, and \cite{Guo Jang Jiang-1,Guo Jang Jiang} via a recent $L^2$-$L^\infty$ framework. When the solutions of compressible Euler equations are consisted by the basic wave patterns (singularities), the convergence has been established in \cite{Huang-Wang-Yang,Huang-Wang-Yang-1,Huang-Wang-Wang-Yang,XZ,Yu}  in one-dimension case, and \cite{WWZ} for multi-dimensional planar rarefaction wave.
There are also lots of literatures on the hydrodynamic limit of Boltzmann equation to  the incompressible fluid equations, see \cite{Bardos,Bardos-2,Golse-Saint-Raymond,Bardos-Ukai,Guo2006,Guo-Liu,Jiang-Masmoudi,Masmoudi-Raymond,Wu} for incompressible Navier-Stokes equations,  \cite{Jang-Kim,HWWX} for incompressible Euler equations, and the references cited therein.

	\smallskip
	
All of the above-mentioned works on the compressible Euler limit were carried out in either spatially periodic domain or the whole space. However, in many important physical models, the physical boundaries occur naturally, and the boundary effects play an important role in the study of hydrodynamic limits in the Boltzmann theory. For initial boundary value problem,  by a formal analysis,  Sone \cite{Sone-2007} showed that the solutions contains three kinds of parts, i.e., interior part, viscous boundary layer and Knudsen boundary layer. Recently, Based on a systematic study of the viscous and Knudsen layers
and the $L^2-L^\infty$ framework, Guo-Huang-Wang \cite{GHW} first justified rigorously the validity of the Hilbert expansion for the hard sphere Boltzmann equation with specular reflection boundary condition in half-space, which leads to derivations of both compressible Euler equations and acoustic equations, see \cite{Jiang-Luo-Tang-1} for Maxwell reflection boundary condition of hard potentials and \cite{Jiang-Luo-Tang-2} for diffuse reflection boundary condition of hard sphere. 
	
In the present paper, we aim to justify the hydrodynamic limit to the compressible Euler equations for the Boltzmann equation of soft potentials. The new difficulty for the soft potentials is that it is hard to establish the existence of solution for Knudsen boundary layer with enough space decay rate, which is crucial to close the Hilbert expansion.
	
To our knowledge, for the specular boundary condition,  the known results \cite{GPS-1988,Jiang-Wang} on the existence of Knudsen boundary layer are for hard sphere, and the exponential space decay was also obtained due to the strong effect of collision frequency $\nu\cong 1+|v|$. For the other boundary conditions, we refer the  readers to  \cite{BCN-1986,CLY,UYY,WYY} for hard potentials and \cite{WYY2} for soft potentials with in-flow boundary condition, \cite{HW} for hard sphere with diffuse reflection boundary condition,  \cite{BG} for hard sphere with phase transition, and the references therein.
	
\smallskip

	
	

	\smallskip
	
We consider  the scaled Boltzmann equation
\begin{equation}\label{c1.1}
F_t+v\cdot\nabla_x F=\frac1{\mathscr{K}_n}Q(F, F),
\end{equation}
where $F(t,x,v)\geq 0$ is the density distribution function for the gas particles with position $x\in\mathbb{R}^3_+=\{x\in\mathbb{R}^3: x_3>0\}$  and velocity $v\in\mathbb{R}^3$ at  time $t>0$, and $\mathscr{K}_n>0$ is  Knudsen number which is proportional to the mean free path. The Boltzmann collision term $Q(F_1,F_2)$ on the right is defined in terms of the following bilinear form
	\begin{align}\label{c1.2}
		Q(F_1,F_2)&\equiv\int_{\mathbb{R}^3}\int_{\mathbb{S}^2} B(v-u,\omega)F_1(u')F_2(v')\,{d\omega du}\nonumber\\
		&\qquad-\int_{\mathbb{R}^3}\int_{\mathbb{S}^2} B(v-u,\omega)F_1(u)F_2(v)\,{d\omega du}\nonumber\\
		&:=Q_+(F_1,F_2)-Q_-(F_1,F_2),
	\end{align}
where the relationship between the post-collision velocity $(v',u')$ of two particles with the pre-collision velocity $(v,u)$ is given by
	\begin{equation*}
		u'=u+[(v-u)\cdot\omega]\omega,\quad v'=v-[(v-u)\cdot\omega]\omega,
	\end{equation*}
for $\omega\in \mathbb{S}^2$, which can be determined by conservation laws of momentum and energy
	\begin{equation*}
		u'+v'=u+v,\quad |u'|^2+|v'|^2=|u|^2+|v|^2.
	\end{equation*}
The Boltzmann collision kernel $B=B(v-u,\omega)$ in \eqref{c1.2} depends only on $|v-u|$ and $\theta$ with $\cos\theta=(v-u)\cdot \omega/|v-u|$.   Throughout this paper,  we  consider cutoff soft potential model, i.e.,
	\begin{equation*}
		B(v-u,\omega)=|v-u|^{\kappa}\cdot \beta(\theta),\quad \kappa\in(-3,0),
	\end{equation*}
where we assume the Grad cutoff condition holds, i.e.,
	\begin{align*}
		0\leq \beta(\theta)\leq \beta_{0}|\cos \theta|,
	\end{align*}
for some constant $\beta_{0}>0$.
	
Denote $\vec{n}=(0,0,-1)$ to be  the outward normal of  $\mathbb{R}^3_+$ and the phase boundary in the space $\mathbb{R}^3_+\times\mathbb{R}^3$ as $\gamma:=\partial\mathbb{R}^3_+\times\mathbb{R}^3$. We split $\gamma$ into outgoing boundary $\gamma_+$, incoming boundary $\gamma_-$, and grazing boundary $\gamma_0$:
	\begin{align}
		\begin{split}\nonumber
			\gamma_+=\{(x,v) : x\in \partial\mathbb{R}^3_+, v\cdot \vec{n}=-v_3>0\},\\
			\gamma_-=\{(x,v) : x\in \partial\mathbb{R}^3_+, v\cdot \vec{n}=-v_3<0\},\\
			\gamma_0=\{(x,v) : x\in \partial\mathbb{R}^3_+, v\cdot \vec{n}=-v_3=0\}.
		\end{split}
	\end{align}
In the present paper, we consider the Boltzmann equation with specular reflection boundary conditions, i.e.,
	\begin{equation}\label{c1.3}
		F(t,x,v)|_{\gamma_-}=F(t,x,R_xv),
	\end{equation}
	where
	\begin{equation}\label{c1.4}
		R_xv=v-2\{v\cdot \vec{n}\} \vec{n}=(v_1,v_2, -v_3)^{t}.
	\end{equation}


\subsection{Asymptotic expansion} From the formal analysis in \cite{Sone-2007}, we know that the thickness of viscous boundary layer is $\sqrt{\mathscr{K}_n}$. For simplicity, we  use the new parameter $\v=\sqrt{\mathscr{K}_n}$ and denote the Boltzmann solution to be $F^\v$, then the Boltzmann equation \eqref{c1.1} is rewritten as
\begin{equation}\label{c1.5}
	\partial_t F^\v+v\cdot \nabla_x F^\v =\frac{1}{\v^2} Q(F^\v, F^\v).
\end{equation}

\subsubsection{Interior expansion} We define the interior  expansion
\begin{equation}\label{c1.6}
	F^\v(t,x,v)\sim \sum_{k=0}^{\infty}\v^{k}F_{k}(t,x,v).
\end{equation}
Substituting \eqref{c1.6} into  \eqref{c1.5} and comparing the order of $\v$, one obtains
\begin{align}\label{c1.7}
	\begin{split}
		\dis \f{1}{\v^2}: \,\,\quad\quad\quad\quad\qquad 0&=Q(F_{0},F_{0}),\\[2mm]
		\dis \f{1}{\v}:  \,\,\quad\quad\quad\quad\qquad 0&=Q(F_{0},F_{1})+Q(F_1,F_0),\\[2mm]
		\dis \v^0:\,\, \{\pt+v\cdot\nax\}F_{0}&=Q(F_{0},F_{2})+Q(F_{2},F_{0})+Q(F_1,F_1),\\[2mm]
		\dis \v: \,\,\{\pt+v\cdot\nax\}F_{1}&=Q(F_{0},F_{3})+Q(F_{3},F_{0})+Q(F_{1},F_{2})+Q(F_2,F_1),\\[2mm]
		\dis \quad\quad\quad\quad\quad\quad\quad\ \vdots\\
		\dis \v^{k}:\,\, \{\pt+v\cdot\nax\}F_{k}&=Q(F_{0},F_{k+2})+Q(F_{k+2},F_{0})+\sum_{\substack{i+j=k+2\\  i,j\geq1}}Q(F_{i},F_{j}).
	\end{split}
\end{align}

It follows from  $\eqref{c1.7}_1$ and the celebrated H-theorem that $F_0$ should be a local Maxwellian, i.e.,
\begin{equation}\label{c1.7-0}
	\mu(t,x,v):=F_0(t,x,v)\equiv \f{\rho(t,x)}{[2\pi T(t,x)]^{3/2}}\exp{\left\{-\f{|v-\fu(t,x)|^2}{2T(t,x)}\right\}},
\end{equation}
where $\rho(t,x)$, $\mathfrak{u}(t,x)=(\fu_1,\fu_2,\fu_3)(t,x),$ and $T(t,x)$ are defined by
\begin{align}\nonumber
	\intr F_0 dv =\rho,\quad \intr vF_0dv =\rho \fu, \quad \intr|v|^2F_0 dv =\rho|\fu|^2+3\rho T,
\end{align}
which represent the macroscopic density, velocity and temperature, respectively. Multiplying $\eqref{c1.7}_3$ by $1,v_i,|v|^2$ and integrating on $\R^3$, one obtains that $(\rho, \fu, T)$ satisfies the compressible Euler system
\begin{equation}\label{c1.8}
	\begin{cases}
		\dis \pt\r+\mbox{\rm div} (\r \fu)=0,\\[2mm]
		\dis \pt(\r\fu)+\mbox{\rm div}(\r \fu\otimes \fu)+\nabla p=0,\\[2mm]
		\dis \pt[\r(\f{3T}{2}+\f{|\fu|^2}{2})]+\mbox{\rm div} [\r \fu(\f{3T}{2}+\f{|\fu|^2}{2})]+\mbox{\rm div}(p\fu)=0,
	\end{cases}
\end{equation}
where $p=\rho T$ is the pressure function. For the compressible Euler equations \eqref{c1.8}, we impose the slip boundary condition
\begin{equation}\label{1.12}
	\fu\cdot \vec{n}|_{x_3=0}=\fu_{3}|_{x_3=0}=0.
\end{equation}
and the initial data
\begin{equation}\label{1.12-2}
	(\rho,\fu,T)(0,x)=(1+\delta\varphi_0,\delta \Phi_0, 1+\delta \vartheta_0)(x),
\end{equation}
with $\|(\varphi_0, \Phi_0,\vartheta_0)\|_{H^{s_0}}\leq 1$ where $\delta>0$ is a  parameter and $s_0\geq 3$ is some given positive number.
Choose $\delta_1>0$ so that for any $\delta\in(0,\delta_1]$, the positivity of $1+\delta\varphi_0$ and $1+\delta \vartheta_0$ is guaranteed. Then for each $\delta\in(0,\delta_1]$, there is a family of classical solutions $(\rho^{\delta},\fu^{\delta}, T^{\delta})\in C([0,\tau^{\delta}]; H^{s_0}(\mathbb{R}^3_+)) \cap C^1([0,\tau^{\delta}]; H^{s_0-1}(\mathbb{R}^3_+))$  of the compressible Euler equations \eqref{c1.8}-\eqref{1.12-2} such that $\rho^{\delta}>0$ and $T^{\delta}>0$.

\smallskip

For later use, we define the linearized collision operator $\mathbf{L}$ by
\begin{equation}\label{c1.9}
	\mathbf{L}\mathfrak{h}=-\frac{1}{\sqrt{\mu}}\Big\{Q(\mu,\sqrt{\mu} \mathfrak{h})+Q(\sqrt{\mu} \mathfrak{h},\mu)\Big\}.
\end{equation}
Denote the null space of $\FL$ as $\mathcal{N}$, it is clear that
\begin{align*}
	\mathcal{N}=\text{span}\{\chi_0,\chi_1,\chi_2,\chi_3,\chi_4\},
\end{align*}
where
\begin{align*}
	\chi_0=\f{1}{\sqrt{\rho}}\sqrt{\mu},\quad \chi_i=\f{1}{\sqrt{\rho T}}(v_i-\fu_i)\sqrt{\mu},\quad \chi_4=\f{1}{\sqrt{6\rho}}(\f{|v-\fu|^2}{T}-3)\sqrt{\mu}.
\end{align*}
For each $k\geq 1$, decompose $\displaystyle f_{k}:=\frac{F_k}{\sqrt{\mu}}$ as
\begin{align}\label{c1.12}
	f_{k}
	&=\mathbf{P} f_{k} +\{\mathbf{I-P}\} f_{k} \nonumber\\
	&\equiv \left\{\frac{\rho_k}{\sqrt{\rho}} \chi_0+\sum_{j=1}^3\sqrt{\frac{\rho}{T}} u_{k,j}\cdot \chi_j+\sqrt{\frac{\rho}{6}}\frac{\theta_k}{T} \chi_4 \right\}+\{\mathbf{I-P}\}f_{k}\nonumber\\
	&\equiv \left\{ \frac{\rho_k}{\rho}+u_{k}\cdot \frac{v-\fu}{T}+\frac{\theta_k}{6T}(\frac{|v-\fu|^2}{T}-3)\right\}\sqrt{\mu}+\{\mathbf{I-P}\}f_{k},
\end{align}
where $\FP$ is the macroscopic projection onto $\mathcal{N}$.
\smallskip

\subsubsection{ Viscous boundary layer expansion}
Generally, the solution of  interior expansion $F_i , i=1,2,\cdots$ do not  satisfy the specular reflection boundary conditions.
To overcome the difficulty coming from the boundary condition, the boundary layer expansion is needed, see \cite{GHW} and \cite{Sone-2002,Sone-2007}.

 We define the scaled normal coordinate:
\begin{equation}\label{c1.13}
	y:=\frac{x_3}{\v}.
\end{equation}
For simplicity of presentation, we denote
\begin{equation}\label{c1.14}
	x_{\shortparallel}=(x_1,x_2),\quad  \nabla_{\sp}=(\partial_{x_1},\partial_{x_2})\quad\mbox{and}\quad  v_{\sp}=(v_1,v_2).
\end{equation}
Motivated by \cite[Section 3.4.1]{Sone-2007}, we define the viscous boundary layer expansion as
\begin{align*}
	\bar{F}^\v(t,x_\sp, y)\sim \sum_{k=1}^\infty \v^k \bar{F}_k(t,x_\sp, y,v).
\end{align*}
Plugging $F^\v+\bar{F}^\v$  into the Boltzmann equation \eqref{c1.5} and comparing the order of $\v$, then using \eqref{c1.7},  in the neighborhood of  physical boundary, we have
\begin{equation}\label{c1.15}
	\begin{split}
		\dis& \f{1}{\v}:  \qquad\quad 0=Q(\mu_0,\bar{F}_{1})+Q(\bar{F}_1,\mu_0),\\[2mm]
		\dis & \v^0:\quad v_3 \frac{\partial\bar{F}_1}{\partial y }=[Q(\mu_0,\bar{F}_{2})+Q(\bar{F}_{2},\mu_0)]+y [Q(\partial_3\mu_0, \bar{F}_1)+Q( \bar{F}_1,\partial_3\mu_0)]\\
		&\qquad\qquad\qquad\qquad+Q(F_1^0,\bar{F}_1)+Q(\bar{F}_1,F_1^0)+Q(\bar{F}_1,\bar{F}_1),\\
		&\dis \quad\quad\quad\quad\quad\quad\quad\ \vdots\\
		\dis & \v^{k}: \quad  \{\partial_t+v_\sp\cdot\nabla_\sp\}\bar{F}_{k}+v_3\frac{\partial\bar{F}_{k+1}}{\partial y }=Q(\mu_0,\bar{F}_{k+2})+Q(\bar{F}_{k+2},\mu_0)\\ &\qquad\qquad\quad+\sum_{\substack{l+j=k+2\\    1\leq l\leq \fb,\, j\geq1}}\frac{y^l}{l!} \big[Q(\partial_3^l\mu_0, \bar{F}_{j})+Q( \bar{F}_{j},\partial_3^l\mu_0) \big] \\
		&\qquad\qquad\quad+\sum_{\substack{i+j=k+2\\   i, j\geq1}}\big[Q(F_i^0,\bar{F}_j)+Q(\bar{F}_j,F_i^0)+Q(\bar{F}_i,\bar{F}_j)\big] \\
		&\qquad\qquad\quad+\sum_{\substack{i+j+l=k+2\\  1\leq l\leq \fb,\, i, j\geq1}} \frac{y^l}{l!} \big[Q(\partial_3^lF_i^0, \bar{F}_{j})+Q( \bar{F}_{j},\partial_3^lF_i^0)\big], \quad \mbox{for}\  k\geq 1,
	\end{split}
\end{equation}
where we have used  the Taylor expansions of $\mu$ and $F_i$ at $x_3=0$, i.e.,
\begin{align}
	\mu(t,x_1,x_2,x_3,v)
	&=\mu_0+\sum_{l=1}^{\fb} \frac{1}{l!} \partial_3^l\mu_0\cdot x_3^l + \frac{x_3^{\fb+1}}{(\fb+1)!} \partial_3^{\fb+1}\tilde{\mu},\label{c1.16}\\
	F_i(t,x_1, x_2, x_3,v)
	&=F_i^0+\sum_{l=1}^{\fb} \frac{1}{l!} \partial_3^l F_i^0\cdot x_3^l + \frac{x_3^{\fb+1}}{(\fb+1)!} \partial_3^{\fb+1}\mathfrak{F}_i,\quad i\geq 1.\label{c1.17}	
\end{align}
Here we have used the simplified notations
\begin{align}\label{c1.18}
	\begin{split}
		\partial_3^l\mu_0:&=(\partial_3^l\mu)(t,x_1, x_2,0,v),
		\quad \partial_3^{\fb+1}\tilde{\mu}:=(\partial_3^{\fb+1}\mu)(t,x_1, x_2, \xi_0,v),\\
		\partial_3^lF_i^0:&=(\partial_3^lF_i)(t,x_1, x_2, 0, v),
		\quad \partial_3^{\fb+1}\mathfrak{F}_i:=(\partial_3^{\fb+1}F_i)(t,x_1, x_2, \xi_i,v),
	\end{split}
\end{align}
for some $\xi_i \in (0,x_3)$ with $i\geq0$.  The number $\fb\in \mathbb{N}_+$ will be chosen later.

\smallskip

For the macro-micro decomposition of viscous and Knudsen boundary layers, we denote the corresponding linearized operator,  macroscopic projection,  and null space as
\begin{align*}
	\mathbf{L}_0=\mathbf{L}(t,x_\sp,0,v),\qquad \mathbf{P}_0=\mathbf{P}(t,x_\sp,0,v), \qquad \mathcal{N}_0=\mathcal{N}(t,x_\sp,0,v).
\end{align*}
It is noted that $\mathbf{L}_0, \mathbf{P}_0$ and $\mathcal{N}_0$ are independent of normal variables. We define
\begin{equation}\label{c1.19}
	\bar{f}_k:=\frac{\bar{F}_k}{\sqrt{\mu_0}},
\end{equation}
then it holds that
\begin{align*}
	\bar{f}_k&=\mathbf{P}_0\bar{f}_k+\{\mathbf{I-P_0}\}\bar{f}_k\nonumber\\
	&= \left\{ \frac{\bar\rho_k}{\rho^0}+\bar{u}_{k}\cdot \frac{v-\fu^0}{T^0}+\frac{\bar{\theta}_k}{6T^0}(\frac{|v-\fu^0|^2}{T^0}-3)\right\}\sqrt{\mu_0}+\{\mathbf{I-P_0}\}\bar{f}_k,
\end{align*}
where and whereafter we always use the notation $(\rho^0,\fu^0,T^0):=(\rho,\fu,T)(t,x_\sp,0)$.

\smallskip

Throughout the present paper, we always assume the far-field condition
\begin{equation}\label{c1.20}
	\bar{f}_k(t,x_\sp,y,v)\rightarrow 0, \quad \mbox{as}\  y\rightarrow+\infty.
\end{equation}

\subsubsection{ Knudsen boundary layer expansion}
To construct the solution satisfying the boundary condition at higher orders, we still need the Knudsen boundary layer. We define the new scaled normal coordinate:
\begin{equation}\nonumber
	\eta:=\frac{x_3}{\v^2}.
\end{equation}
The  Knudsen boundary layer expansion is defined as
\begin{equation}\nonumber
	\hat{F}^\v(t,x_\sp, \eta)\sim \sum_{k=1}^\infty \v^k \hat{F}_k(t,x_\sp, \eta,v).
\end{equation}
Plugging $F^\v+\bar{F}^\v+\hat{F}^\v$  into \eqref{c1.5} and comparing the order of $\v$,  then using \eqref{c1.7} and \eqref{c1.15}, one obtains
\begin{equation}\label{c1.22}
	\begin{split}
		\dis& \f{1}{\v}:  \quad v_3 \frac{\partial\hat{F}_1}{\partial \eta }-\big[Q(\mu_0,\hat{F}_{1})+Q(\hat{F}_1,\mu_0)\big]=0,\\[2mm]
		\dis & \v^0:\quad \  v_3 \frac{\partial\hat{F}_2}{\partial \eta }-\big[ Q(\mu_0,\hat{F}_{2})+Q(\hat{F}_{2},\mu_0) \big]\\
		&\qquad\quad= Q(F_1^0+\bar{F}^0_1,\hat{F}_1)+Q(\hat{F}_1,F_1^0+\bar{F}^0_1)
		+Q(\hat{F}_1,\hat{F}_1),\\[2mm]
		&\dis \quad\quad\quad\quad\quad\quad\quad\ \vdots \\
		\dis & \v^{k}:\quad   v_3\frac{\partial\hat{F}_{k+2}}{\partial \eta }-\big[ Q(\mu_0,\hat{F}_{k+2})+Q(\hat{F}_{k+2},\mu_0) \big] \\
		&\qquad=-\{\partial_t+v_\sp\cdot\nabla_\sp\}\hat{F}_{k}+\sum_{\substack{j+2l=k+2\\  1\leq l\leq \fb, j\geq1}} \frac{\eta^l}{l!} \big[Q(\partial_3^l\mu_0, \hat{F}_{j})+Q( \hat{F}_{j},\partial_3^l\mu_0) \big] \\
		&\qquad\quad+\sum_{\substack{i+j=k+2\\ i,j\geq1}}\big[Q(F_i^0+\bar{F}_i^0,\hat{F}_j)+Q(\hat{F}_j,F_i^0+\bar{F}_i^0)+Q(\hat{F}_i,\hat{F}_j)\big] \\
		&\qquad\quad+\sum_{\substack{i+j+2l=k+2\\ i,j\geq1, 1\leq l\leq \fb}}\frac{\eta^l}{l!} \big[Q(\partial_3^l F_i^0, \hat{F}_{j})+Q( \hat{F}_{j},\partial_3^l F_i^0) \big]   \\
		&\qquad\quad+\sum_{\substack{i+j+l=k+2\\ i,j\geq1, 1\leq l\leq \fb}}\frac{\eta^l}{l!} \big[Q(\partial_y^l \bar{F}_i^0, \hat{F}_{j})+Q( \hat{F}_{j},\partial_y^l \bar{F}_i^0) \big],
		\quad \mbox{for}\  k\geq 1,
	\end{split}
\end{equation}
where we have used \eqref{c1.16}-\eqref{c1.17} and  the Taylor expansion of $\bar{F}_i$
\begin{equation*}
	\bar{F}_i(t, x_1, x_2, y,v)
	=\bar{F}_i^0+\sum_{l=1}^{\fb} \frac{1}{l!} \partial_y^l \bar{F}_i^0\cdot y^l + \frac{y^{\fb+1}}{(\fb+1)!} \partial_y^{\fb+1}\bar{\mathfrak{F}}_i,
\end{equation*}
with
\begin{align}\label{c1.23}
	\begin{split}
		\partial_y^l\bar{F}_i^0:&=(\partial_y^l\bar{F}_i)(t,x_1, x_2, 0, v),\\
		\partial_y^{\fb+1}\bar{\mathfrak{F}}_i:&=(\partial_y^{\fb+1}\bar{F}_i)(t,x_1, x_2, \bar{\xi}_i,v),
	\end{split}
	\qquad \mbox{for}\,\, 0\leq l\leq \fb,
\end{align}
for some $\bar{\xi}_i\in [0,y]$. For later use, we denote $\displaystyle\hat{f}_{k}:=\frac{\hat{F}_{k}}{\sqrt{\mu_0}}$.

\subsection{Hilbert expansion}
Now we consider the Boltzmann solution with the following  Hilbert expansion of multi-scales
\begin{align}\label{c1.24}
	F^\v=\mu(t,x,v)+\sum_{i=1}^{N} \v^i F_i(t,x,v)+\sum_{i=1}^{N} \v^i\bar{F}_i(t,x_\sp,\frac{x_3}{\v},v)
	+\sum_{i=1}^{N} \v^i\hat{F}_i(t,x_\sp,\frac{x_3}{\v^2},v)+\v^5 F^\v_{R}.
\end{align}
Using \eqref{c1.7}, \eqref{c1.15} and \eqref{c1.22}, one can obtain the equation of $F_R^{\v}$
\begin{align}\label{c1.25}
	&\pt F^{\v}_{R}+v\cdot \nabla_x F^{\v}_R-\frac{1}{\v^2}\{Q(\mu, \FRep)+Q(\FRep, \mu)\}\nonumber\\
	&=\v^3Q(\FRep, \FRep)+\sum_{i=1}^N \v^{i-2}\{Q(F_i+\bar{F}_i+\hat{F}_i,F^\v_R)+Q(F^\v_R, F_i+\bar{F}_i+\hat{F}_i)\}\nonumber\\
	&\quad+ R^\v+\bar{R}^\v+\hat{R}^\v,
\end{align}
where $R^\v,\bar{R}^\v$ and $\hat{R}^\v$ are defined in \eqref{H2}-\eqref{H4}.\\

\smallskip

The main purpose of the present paper is to  establish the validity of the Hilbert expansion for the Boltzmann equation around the local Maxwellian $\mu$ determined by compressible Euler equations \eqref{c1.8}. For later use, we define
\begin{equation}\label{c1.28}
	F^{\v}_{R}=\sqrt{\mu} f^\v_R.
\end{equation}
To use the $L^2$-$L^\infty$ framework \cite{Guo Jang Jiang-1,Guo2010}, we also introduce a global Maxwellian
\begin{equation*}
	\mu_{M}:=\frac{1}{(2\pi T_{M})^{3/2}}\exp{\left\{-\frac{|v|^{2}}{2T_{M}}\right\}},
\end{equation*}
where $T_{M}>0$ satisfies the condition
\begin{equation}\label{c1.29}
	T_{M}<\min_{x\in\mathbb{R}^3_+} T(t,x)\leq \max_{x\in\mathbb{R}^3_+} T(t,x)<2T_{M}.
\end{equation}
Using \eqref{c1.29}, one can easily deduce that there exists a positive constant $C>0$ such that for some $\frac{1}{2}<\alpha<1$ , the following holds:
\begin{equation}\label{c1.30}
	\frac1C\mu_{M}\leq \mu(t,x,v)\leq C\mu_{M}^{\alpha}.
\end{equation}
We further define
\begin{equation}\label{c1.31}
	F_{R}^{\v}=\{1+|v|^{2}\}^{-\frac{\mathfrak{k}}{2}}\sqrt{\mu_{M}} h^\v_R\equiv\frac{1}{\varpi_\mathfrak{k}(v)}\sqrt{\mu_{M}}\hep_R,
\end{equation}
with $\mathfrak{k}\geq0$ and $\varpi_\mathfrak{k}:=(1+|v|^2)^{\f{\mathfrak{k}}{2}}$.

\begin{theorem}\label{theorem}
	Let $\tau^\d>0$ be the life-span of smooth solution of compressible Euler equations \eqref{c1.8}. Let  $\mathfrak{k}\geq16$, $N\geq 6$ and $\fb\geq 5$.  We assume the initial data
	\begin{align}
		F^{\v}(0,x,v)&=\mu(0,x,v)+\sum_{i=1}^{N}\v^i \left\{F_i(0,x,v)+ \bar{F}_i(0,x_\sp, \frac{x_3}{\v},v)+ \hat{F}_i(0,x_\sp,\frac{x_3}{\v^2},v)\right\}\nonumber\\
		&\qquad+\v^5F^\v_R(0,x,v)\geq 0,\nonumber
	\end{align}
	and $F_i(0), \bar{F}_i(0),  i=1,\cdots, N$  satisfy the regularity and compatibility conditions described in  Proposition \ref{prop5.1}, and
	\begin{align*}
		\Big\|(\frac{F^\v_R}{\sqrt{\mu}})(0)\Big\|_{L^2_{x,v}}+\v^3\Big\|(\varpi_\mathfrak{k} \frac{F^\v_R}{\sqrt{\mu_M}})(0)\Big\|_{L^\infty_{x,v}}
		<\infty.
	\end{align*}
	Then there exists a small positive constant $\v_0>0$  such that the IBVP problem \eqref{c1.5} and \eqref{c1.3} has a unique solution for  $\v\in(0,\v_0]$ over the time interval $t\in [0,\tau^\d]$ in the following form of expansion
	\begin{align}\label{c1.32}
		F^{\v}(t,x,v)&=\mu(t,x,v)+\sum_{i=1}^{N}\v^i \left\{F_i(t,x,v)+ \bar{F}_i(t,x_\sp, \frac{x_3}{\v},v)+ \hat{F}_i(t,x_\sp,\frac{x_3}{\v^2},v)\right\}\nonumber\\
		&\qquad+\v^5F^\v_R(t,x,v)\geq 0,
	\end{align}
	with
	\begin{align}
		&\sup_{t\in[0,\tau^\d]} \left\{\Big\|\frac{F^\v_{R}(t)}{\sqrt{\mu}}\Big\|_{L^2_{x,v}}+\v^3\Big\|\varpi_\mathfrak{k}(v)\frac{F^\v_{R}(t)}{\sqrt{\mu_M}}\Big\|_{L^\infty_{x,v}}\right\}\leq C(\tau^\d)<\infty.\label{c1.33}
	\end{align}
	Here the functions $F_i(t,x,v), \bar{F}_i(t,x_{\sp},y,v)$ and $\hat{F}_i(t,x_{\sp},\eta,v)$ are the interior expansion, viscous and Knudsen boundary layers   respectively constructed in Proposition \ref{prop5.1}.
\end{theorem}
\begin{remark}
	From  \eqref{c1.32}-\eqref{c1.33} and the uniform estimates in  Proposition \ref{prop5.1}, it is direct to check that
	\begin{align}\nonumber
		\sup_{t\in[0,\tau^\d]}\left\{ \Big\|\Big(\frac{F^\v-\mu}{\sqrt{\mu}}\Big)(t)\Big\|_{L^2(\R_+^3\times\R^3)}+ \Big\|\varpi_\mathfrak{k} \, \Big(\frac{F^\v-\mu}{\sqrt{\mu_M}}\Big)(t)\Big\|_{L^\infty(\R_+^3\times\R^3)}
		\right\}\leq C \v\rightarrow 0.
	\end{align}
	Hence we have established  the hydrodynamic limit from the Boltzmann equation to the compressible Euler system for the half-space  problem.
\end{remark}

\begin{remark}
For simplicity of presentation, we only give details of proof for the Boltzmann equation of soft potentials in the present paper. And we point out that it is also valid for the cases of hard potentials by similar arguments.
\end{remark}

\medskip

Now we briefly comment the key points of present paper. To estimate the microscopic part of interior expansions and viscous boundary layers, we need some decay property on  pseudo-inverse linear operator $\FL^{-1}$ and $\FL_0^{-1}$. For $-\f32<\kappa\leq1$,  the authors \cite{Jiang-Luo-Tang} obtained
	\begin{align*}
		|\mu^{-\f{q}{2}}\FL^{-1}\mathfrak{g}(v)|\lesssim\|\mu^{-\f{q'}{2}}\mathfrak{g}\|_{L^\infty_v},\quad 0<q'<1.
	\end{align*}
Due to the strong singularity, it is hard to establish above estimate for
$-3<\kappa\leq \f32$. In this paper, by observing the feature of Hilbert expansion on interior parts and viscous boundary layers, we can get the following control
\begin{align}\label{c1.44}
	&|\mu^{-\f{q}{2}}\FL^{-1}\mathfrak{g}(v)|^2
	 \lesssim \sum\limits_{0\leq \alpha\leq N}\|\partial_v^\alpha\{\mu^{-\f{q}{2}}\FL^{-1}\mathfrak{g}\}\|_{L^2}^2
	 \lesssim \sum\limits_{0\leq \alpha\leq N}\|\nu^{-1}\mu^{-\f{q'}{2}}\partial_v^\alpha \mathfrak{g}\|_{L^2}^2,\,\, N\geq2,\, 0<q<q'<1,
\end{align}
by losing velocity derivatives, see section 2.2 for details. We point out that the losing velocity derivatives is natural since the interior parts and viscous boundary layers always possess enought regularity with respect to $v\in \R^3$.


\smallskip

The construction of  Knudsen layers  is more delicate for Boltzmann equation of soft potentials. Noting \eqref{c1.22}, to solve the Knudsen layer, it is equivalent to study the following linear boundary value problem
\begin{align}\label{c1.45}
	\begin{cases}
		v_3\partial_\eta f+ {\nu^0}(v) f-K^0f=\mathfrak{g},\\
		f(0,v)|_{v_3>0}=f(0,R_\eta v),\\
		\lim\limits_{\eta\rightarrow\infty}f(\eta,v)=0.
	\end{cases}
\end{align}
Especially, by noting the right hand side (RHS) of \eqref{c1.22}, we have to get at least space polynomial decay for the solution  of \eqref{c1.45} to continue the construction of Hilbert expansion. For hard sphere case, one can obtain even exponential decay with the help of  strong effect of collision frequency $\nu(v)\cong 1+|v|$. However, it is hard for the cases of soft potentials since the effect of collision frequency $\nu(v)\cong (1+|v|)^{\kappa}\to 0$ as $|v|\to \infty$ is very weak.


\smallskip

To solve \eqref{c1.45}, we first establish the {\it a prior} uniform $L^\infty$ estimate  for an approximate problem (see \eqref{D2.42}), i.e.,
\begin{align}\label{1.34}
	\|w_lf^\lambda\|_{L^\infty_{\eta,v}}+|w_lf^\lambda|_{L^\infty(\gamma_+)}\leq C \Big(\|\sqrt{\nu^0}f^\lambda\|_{L^2_{\eta,v}} +   \|(\nu^0)^{-1}w_lg\|_{L^\infty_{\eta,v}}+|w_{l+4}r|_{L^\infty(\gamma_+)}\Big),
\end{align}
see Lemma \ref{lem3.4} for details. Here we point out that the constant in \eqref{1.34} is independent of the length of domain $\Omega=[0,d]$. For soft potentials, since the collision effect is weak, the key point is to take the number of collisions with the boundaries to depend on $v$, that is $k=\tilde{k}_0|v_3| (1+|v|)^{|\kappa|}$ with $\tilde{k}_0\gg 1$.

Under some constraint conditions, we can have the following $L^2_{\eta,v}$ decay estimate  by lossing velocity weight arguments
\begin{equation}\label{D2.94-1}
\int_0^d(1+\eta )^n\|w_l\bar{f}\|_{\nu}^2d\eta \leq C_n\int_0^d(1+\eta )^{2p_n}\|w_{l+2n+2}g\|_{L^2_v}^2d\eta ,\quad p_n>\f{n}{2}+1,
\end{equation}
see Lemmas \ref{lem3.10}-\ref{lem3.11} for details. For the space decay rate in $L^\infty_{\eta,v}$,  we multiply \eqref{c1.45} by $(1+\eta)^n$ to obtain
    \begin{align*}
    	v_3\partial_\eta\{(1+\eta)^n f\}+{\nu^0}(v)\{(1+\eta)^nf\}-K^0\{(1+\eta)^nf\}=(1+\eta)^n\mathfrak{g}+nv_3(1+\eta)^{n-1}f,
    \end{align*}
which yields that
\begin{align*}
	\|w_l\, (1+\eta)^nf\|_{L^\infty_{\eta,v}}
	&\lesssim \|w_{l+4}\,(1+\eta)^{n-1}f\|_{L^\infty_{\eta,v}}	+\|(\nu^0)^{\f12}(1+\eta)^nf\|_{L^2_{\eta,v}}\\
	&\quad +\|(\nu^0)^{-1}w_l\,(1+\eta)^n\mathfrak{g}\|_{L^\infty_{\eta,v}}.
\end{align*}
Then, using \eqref{D2.94-1} and an induction arguments on $n$,  we finally obtain that
\begin{align*}
		\|w_l\,(1+\eta)^nf\|_{L^\infty_{\eta,v}}\lesssim \|w_{l+4n+4}\,(1+\eta)^{q_n}\mathfrak{g}\|_{L^\infty_{\eta,v}},\quad \text{for}\ q_n>n+\f32.
\end{align*}
With above estimates, we obtain the existence of Knudsen boundary layer problem with enough space decay estimate in $L^\infty_{\eta,v}$.

\smallskip

With the help of above estimates on $\FL^{-1}$ and Knudsen boundary layer,
by the same arguments as in \cite{GHW}, we can establish Hilbert expansion of Boltzmann equation of soft potentials with multi-scales in half-space.

\medskip

The paper is organized as follows. In section \ref{sec2}, we give some basic estimates on collision operator and establish the decay estimate of $\FL^{-1}$ for soft potentials. Section \ref{section2} is devoted to existence of Knudsen boundary layer of soft potentials with enough space decay rate. In section \ref{section4}, we construct the Hilbert expansion of soft Boltzmann equation and prove Theorem \ref{theorem}.

\medskip

\noindent{\bf Notations.} Throughout the present paper, $C$ denotes a generic positive constant  and  vary from line to line. And $C(a),C(b),\cdots$ denote the generic positive constants depending on $a,~b,\cdots$, respectively, which also may vary from line to line. We use $\langle\cdot ,\cdot \rangle$ to denote the standard $L^2$ inner product in $\R^3_v$.
$\|\cdot\|_{L^2}$ denotes the standard $L^2(\mathbb{R}^3_+\times\mathbb{R}^3_v)$-norm, $\|\cdot\|_{L^\infty}$ denotes the $L^\infty(\mathbb{R}^3_+\times\mathbb{R}^3_v)$-norm and $\|\cdot\|_{\nu}$ denotes $\langle \nu\cdot,\cdot\rangle^{\f12}$.

\

\section{Some Estimates for Soft Boltzmann Operators}\label{sec2}
\subsection{Preliminaries}
It follows from 
\eqref{c1.9} that $\FL \mathfrak{h}=\nu(v) \mathfrak{h}-K\mathfrak{h}$ with
\begin{equation}\label{nu}
	\nu(v)=\iint_{\R^3\times \mathbb{S}^2}B(v-u,\omega)\mu(u)\,d\omega du\cong (1+|v|)^{\kappa} \quad\text{and}\quad
	K\mathfrak{h}=K_{1}\mathfrak{h}-K_{2}\mathfrak{h},
\end{equation}
where
\begin{equation}\label{K1}
	(K_{1}\mathfrak{h})(v)=\mu^{\frac{1}{2}}(v)\iint_{\R^3\times \mathbb{S}^2}\mathfrak{h}(u)\mu^{\frac{1}{2}}(u)B(v-u,\omega)\,d\omega du,
\end{equation}
and
\begin{align}
	(K_{2}\mathfrak{h})(v)&=\f{1}{\sqrt{\mu}}\iint_{\R^3\times \mathbb{S}^2}B(v-u,\omega)\mu(v')[\sqrt{\mu(u')}\mathfrak{h}(u')+\mu(u')\sqrt{\mu(v')}\mathfrak{h}(v')]dud\omega\nonumber\\
	&=\mu^{\frac{1}{2}}(v)\iint_{\R^3\times \mathbb{S}^2}[\mu^{-\frac{1}{2}}(v')\mathfrak{h}(v')+\mu^{-\frac{1}{2}}(u')\mathfrak{h}(u')]\mu(u)B(v-u,\omega)\,d\omega du,\nonumber\\
	&=2\mu^{\frac{1}{2}}(v)\iint_{\R^3\times \mathbb{S}^2}\mu^{-\frac{1}{2}}(v')\mathfrak{h}(v')\mu(u)B(v-u,\omega)\,d\omega du.\label{K2}
\end{align}
By standard arguments as in \cite{GL}, we can rewrite $K_i\mathfrak{h}=\int_{\R^3} k_i(u,v)\mathfrak{h}(u)du$ with
\begin{align*}
	&\displaystyle k_1(u,v)=C|v-u|^\kappa \mu^{\f12}(v)\mu^{\f12}(u),\\
	&\displaystyle 0\leq k_2(u,v)\leq \f{C_{\kappa}}{|u-v|^{\frac{3-\kappa}{2}}} \exp{\Big\{-\f18|v-u|^2-\f18\f{(|v-\fu|^2-|u-\fu|^2)^2}{|u-v|^2}\Big\}}.
\end{align*}
It is well-known that there is a positive constant $c_0>0$ such that
\begin{align*}
	\langle \FL \mathfrak{h},\mathfrak{h}\rangle \geq c_0\|\{\FI-\FP\}\mathfrak{h}\|_{\nu}^2.
\end{align*}

For soft potentials, motivated by \cite{Guo2}, we define a monotone cutoff function $\chi_z(s)\in C^\infty(0,\infty)$ satisfying
\begin{equation}\label{chi}
	\chi_z(s)\equiv 0\,\,\text{for }0\leq s\leq z,\quad \chi_z(s)\equiv 1\,\,\text{for }s\geq 2z,\quad 0\leq \chi_z(s)\leq 1\,\,\text{for all }s>0,
\end{equation}
where $z$ is a parameter.
  Define
\begin{align*}
	(K^m\mathfrak{h})(v)&=\int_{\R^3}\int_{\mathbb{S}^2}	 B(|v-u|,\theta)\tilde{\chi}_m(|v-u|)\sqrt{\mu(u)\mu(v)}\mathfrak{h}(u)d\omega du\\
	&-\int_{\R^3}\int_{\mathbb{S}^2}B(v-u,\theta)\tilde{\chi}_m(|v-u|)\sqrt{\mu(u)\mu(u')}\mathfrak{h}(v')d\omega du\\
	&-\int_{\R^3}B(|v-u|,\theta)\tilde{\chi}_m(|v-u|)\sqrt{\mu(u)\mu(v')}\mathfrak{h}(u')d\omega du\\
	&=:K_1^m\mathfrak{h}(v)-K_2^m\mathfrak{h}(v),
\end{align*}
and $K^c=K-K^m$, where $\tilde{\chi}_m=1-\chi_m$. We denote
\begin{align}\label{a46}
	(K^m\mathfrak{h})(v)=\int_{\R^3}k^m(v,u)\mathfrak{h}(u)du,\quad (K^c\mathfrak{h})(v)=\int_{\R^3}k^c(v,u)\mathfrak{h}(u)du.
\end{align}

\begin{lemma}[{\cite{Duan}}]\label{lem3.2} For any $0<m\leq1$, it holds that
	\begin{align}\label{D2.16}
		|(K^m\mathfrak{h})(v)|\leq Cm^{3+\kappa}e^{-\f{|v-\fu|^2}{10}}\|\mathfrak{h}\|_{L^\infty_{v}},
	\end{align}
	where $C>0$ is independent of $m$. The kernels $k^m(v,u)$ and $k^c(v,u)$ satisfy
	\begin{align}\label{D2.17}
		|k^m(v,u)|\leq C_\kappa\{|v-u|^{\kappa}+|v-u|^{-\f{3-\kappa}{2}}\}e^{-\f{|v-\fu|^2+|u-\fu|^2}{16}},
	\end{align}
	and
	\begin{align}\label{D2.18}
		|k^c(v,u)|\leq & \f{C_\kappa m^{a(\kappa-1)}}{|v-u|^{1+\f{(1-a)}{2}(1-\kappa)}}\f{1}{(1+|v-\fu|+|u-\fu|)^{a(1-\kappa)}}e^{-\f{|v-u|^2}{10}}e^{-\f{||v-\fu|^2-|u-\fu|^2|^2}{16|v-u|^2}}\nonumber\\
		&+C|v-u|^\kappa e^{-\f{|v-\fu|^2}{4}}e^{-\f{|u-\fu|^2}{4}},
	\end{align}
	where $a\in[0,1]$ is an arbitrary constant and $C_\kappa$ depending only on $\kappa$.
\end{lemma}

\begin{remark}
The original version of Lemma \ref{lem3.2} was proved in \cite{Duan} for the global Maxwellian. And it is direct to check that it is still valid for local Maxwellian. We omit the details for simplicity of presentation.
\end{remark}

Denote
\begin{align*}
	\tilde{w}(v)=(1+|v|^2)^{\f{l}{2}}\mu^{-\fa} \quad \mbox{with}\,\, l\geq 0, 0\leq \fa<\f12,
\end{align*}
and
\begin{align*}
K_{\tilde{w}}^c h\equiv	\tilde{w}K^c(\f{h}{\tilde{w}})=\int_{\R^3}k_{\tilde{w}}^c(v,u)h(u)du.
\end{align*}
Then, from Lemma \ref{lem3.2},  it is clear that
\begin{align}
	&\int_{\R^3}|k_{\tilde{w}}^c(v,u)|e^{\f{|v-u|^2}{32}}du\leq Cm^{\kappa-1}(1+|v|)^{\kappa-2},\label{D2.19}\\
	&\int_{\R^3}|k_{\tilde{w}}^c(v,u)|e^{\f{|v-u|^2}{32}}du\leq C(1+|v|)^{-1},\label{D2.20}
\end{align}
where $C$ is a constant independent of $m$.

\begin{lemma}[\cite{Liu-Yang}]\label{lem2.10}
Let $\Gamma(\mathfrak{h},\mathfrak{g})=\f{1}{\sqrt{\mu}}Q(\sqrt{\mu}\mathfrak{h},\sqrt{\mu}\mathfrak{g})
$. 	For $\kappa\in(-3,0)$, it holds that
\begin{align}\label{D2.22}
		\Big|\int_{\R^3}\Gamma(\mathfrak{g}_1,\mathfrak{g}_2)\mathfrak{g}_3dv\Big|\leq C\{\|\mathfrak{g}_3\|_{\nu}\|\mathfrak{g}_2\|_{\nu}\|\varpi_k \mathfrak{g}_1\|_{L^\infty}+\|\mathfrak{g}_3\|_{\nu}\|\mathfrak{g}_1\|_{\nu}\|\varpi_k \mathfrak{g}_2\|_{L^\infty}\}, \quad k>\f32.
	\end{align}
\end{lemma}

\smallskip

\subsection{Estimate for $\mathbf{L}^{-1}$}\label{sec2.2}
To consider the derivatives for operators $K_1,K_2$, we denote $\xi:=u-v$. Then one can rewrite $K_1\mathfrak{h}$ as
$$
\begin{aligned}
	K_{1}\mathfrak{h}(v)&=\iint_{\R^3\times \mathbb{S}^2}|u-v|^{\kappa}\beta(\theta)\mu^{\frac{1}{2}}(v)\mu^{\frac{1}{2}}(u)\mathfrak{h}(u)\,d\omega du\\
	&=\iint_{\R^3\times \mathbb{S}^2}|\xi|^{\kappa}\beta(\theta)\mu^{\frac{1}{2}}(v)\mu^{\frac{1}{2}}(v+\xi)\mathfrak{h}(v+\xi)\,d\omega d\xi,
\end{aligned}
$$
which yields that
$$
\begin{aligned}
	\partial_{v}^{\alpha}(K_{1}\mathfrak{h})&=\iint_{\R^3\times \mathbb{S}^2}|\xi|^{\kappa}\beta(\theta)\mu^{\frac{1}{2}}(v)\mu^{\frac{1}{2}}(v+\xi)\partial_{v}^{\alpha}\mathfrak{h}(v+\xi)\,d\omega d\xi\\
	&\quad+\sum\limits_{0\leq \alpha'<\alpha}C_{\alpha}^{\alpha'}\iint_{\R^3\times \mathbb{S}^2}|\xi|^{\kappa}\partial_{v}^{\alpha-\alpha'}\big(\beta(\theta)\mu^{\frac{1}{2}}(v)\mu^{\frac{1}{2}}(v+\xi)\big)\partial_{v}^{\alpha'}\mathfrak{h}(v+\xi)\,d\omega d\xi,
\end{aligned}
$$
where $\alpha=(\alpha_{1},\alpha_{2},\alpha_3)$ is the multi-index, and  $\partial_{v}^{\alpha}:=\partial_{v_1}^{\alpha_{1}}\partial_{v_{2}}^{\alpha_{2}}\partial_{v_{3}}^{\alpha_{3}}$. For small positive number $\epsilon$, it is direct to check that
$$
|\partial_{v}^{\alpha-\alpha'}\big(\beta(\theta)\mu^{\frac{1}{2}}(v)\mu^{\frac{1}{2}}(v+\xi)\big)|\leq C_{\epsilon,N}\mu^{\frac{1}{2}-\epsilon}(v)\mu^{\frac{1}{2}-\epsilon}(v+\xi).
$$
Hence one obtains
\begin{equation}\label{pk1}
	\begin{aligned}
		|\partial_{v}^{\alpha}(K_1\mathfrak{h})|&\leq \iint_{\R^3\times \mathbb{S}^2}|\xi|^{\kappa}\beta(\theta)\mu^{\frac{1}{2}}(v)\mu^{\frac{1}{2}}(v+\xi)|\partial_{v}^{\alpha}\mathfrak{h}(v+\xi)|\,d\omega d\xi\\
		&\quad +C_{\epsilon,N}\sum\limits_{0\leq \alpha'<\alpha}\int_{\R^3}|\xi|^{\kappa}\mu^{\frac{1}{2}-\epsilon}(v)\mu^{\frac{1}{2}-\epsilon}(v+\xi)|\partial_{v}^{\alpha'}\mathfrak{h}(v+\xi)|\,d\xi\\
		&\leq \iint_{\R^3\times \mathbb{S}^2}|v-u|^{\kappa}\beta(\theta)\mu^{\frac{1}{2}}(v)\mu^{\frac{1}{2}}(u)|\partial_{u}^{\alpha}\mathfrak{h}(u)|\,d\omega du\\
		&\quad +C_{\epsilon,N}\sum\limits_{0\leq \alpha'<\alpha}\int_{\R^3}|v-u|^{\kappa}\mu^{\frac{1}{2}-\epsilon}(v)\mu^{\frac{1}{2}-\epsilon}(u)|\partial_{u}^{\alpha'}\mathfrak{h}(u)|\,du\\
		&=:I_{1}+I_{2}.
	\end{aligned}
\end{equation}

It follows from $\mu^{\frac{1}{2}}(v)\mu^{\frac{1}{2}}(u)=\mu^{\frac{1}{2}}(v')\mu^{\frac{1}{2}}(u')$ and \eqref{K2} that
$$
\begin{aligned}
	K_2\mathfrak{h}&=2\iint_{\R^3\times \mathbb{S}^2}\mu^{\frac{1}{2}}(u')\mu^{\frac{1}{2}}(u)\mathfrak{h}(v')|v-u|^{\kappa}\beta(\theta)\,d\omega du\\
	&=2\iint_{\R^3\times \mathbb{S}^2}\mu^{\frac{1}{2}}(v+\xi_{\perp})\mu^{\frac{1}{2}}(v+\xi)\mathfrak{h}(v+\xi_{\shortparallel})|\xi|^{\kappa}\beta(\theta)\,d\omega d\xi,
\end{aligned}
$$
which implies that
$$
\begin{aligned}
	\partial_{v}^{\alpha}K_2\mathfrak{h}&=2\iint_{\R^3\times \mathbb{S}^2}\mu^{\frac{1}{2}}(v+\xi_{\perp})\mu^{\frac{1}{2}}(v+\xi)\partial_{v}^{\alpha}\mathfrak{h}(v+\xi_{\shortparallel})|\xi|^{\kappa}\beta(\theta)\,d\omega d\xi\\
	&\qquad +2\sum\limits_{0\leq \alpha'<\alpha}C_{\alpha}^{\alpha'}\iint_{\R^3\times \mathbb{S}^2}\partial_{v}^{\alpha-\alpha'}\big(\mu^{\frac{1}{2}}(v+\xi_{\perp})\mu^{\frac{1}{2}}(v+\xi)\big)\partial_{v}^{\alpha'}\mathfrak{h}(v+\xi_{\shortparallel})|\xi|^{\kappa}\beta(\theta)\,d\omega d\xi,
\end{aligned}
$$
where $\xi_\shortparallel:=[(u-v)\cdot\omega]\omega$ and $\xi_\perp:=\xi-\xi_\shortparallel$.
It is clear that
$$
|\partial_{v}^{\alpha-\alpha'}\big(\mu^{\frac{1}{2}}(v+\xi_{\perp}\big)\mu^{\frac{1}{2}}(v+\xi))|\leq C_{\epsilon,N}\mu^{\frac{1}{2}-\epsilon}(v+\xi_{\perp})\mu^{\frac{1}{2}-\epsilon}(v+\xi).
$$
Then we can obtain
\begin{equation}\label{pk2}
	\begin{aligned}
		|\partial_{v}^{\alpha}(K_2\mathfrak{h})|&\leq 2\iint_{\R^3\times \mathbb{S}^2}|\xi|^{\kappa}\beta(\theta)\mu^{\frac{1}{2}}(v+\xi_{\perp})\mu^{\frac{1}{2}}(v+\xi)|\partial_{v}^{\alpha}\mathfrak{h}(v+\xi_{\shortparallel})|\,d\omega d\xi\\
		&\quad +C_{\epsilon,N}\sum\limits_{0\leq \alpha'<\alpha}\iint_{\R^3\times \mathbb{S}^2}|\xi|^{\kappa}\beta(\theta)\mu^{\frac{1}{2}-\epsilon}(v+\xi_{\perp})\mu^{\frac{1}{2}-\epsilon}(v+\xi)|\partial_{v}^{\alpha'}\mathfrak{h}(v+\xi_{\shortparallel})|\,d\omega d\xi\\
		&\leq2 \iint_{\R^3\times \mathbb{S}^2}|v-u|^{\kappa}\beta(\theta)\mu^{\frac{1}{2}}(u')\mu^{\frac{1}{2}}(u)|\partial_{v'}^{\alpha}\mathfrak{h}(v')|\,d\omega du\\
		&\quad +C_{\epsilon,N}\sum\limits_{0\leq \alpha'<\alpha}\iint_{\R^3\times \mathbb{S}^2}|v-u|^{\kappa}\beta(\theta)\mu^{\frac{1}{2}-\epsilon}(u')\mu^{\frac{1}{2}-\epsilon}(u)|\partial_{v'}^{\alpha'}\mathfrak{h}(v')|\,d\omega du\\
		&=:J_{1}+J_{2}.
	\end{aligned}
\end{equation}

  Recall \eqref{chi}, we split $I_1$ as
$$
I_1=I_1^{1-\chi_r}+I_1^{\chi_r},
$$
where
\begin{equation}\label{K1chi}
	I_{1}^{\Lambda}=\iint_{\R^3\times \mathbb{S}^2}\Lambda(|v-u|)|v-u|^{\kappa}\beta(\theta)\mu^{\frac{1}{2}}(v)\mu^{\frac{1}{2}}(u)|\partial_{u}^{\alpha}\mathfrak{h}(u)|\,d\omega du, \quad \Lambda=1-\chi_r \ \text{or}\ \chi_r.
\end{equation}
Similarly, we shall also use the notations $I_1^\Lambda,J_1^\Lambda,J_2^\Lambda$ and $(\partial_v^\alpha K_i\mathfrak{h})^\Lambda$.
 \begin{definition}[{Pseudo-inverse operator of $\FL$}]
 	The inverse operator
 \begin{align*}
 	(\FL|_{\mathcal{N}^\perp})^{-1}:\mathcal{N}^\perp\rightarrow \mathcal{N}^\perp,
 \end{align*}	
 	is called the pseudo-inverse operator of the linearized Boltzmann collision operator $\FL$, which
 	is briefly denoted by $\FL^{-1}$.
 \end{definition}

\begin{proposition}\label{thmA.1}
	Let $N\geq 2$, $-3<\kappa<0$, $0<q<q'<1$. Assume that $\mathfrak{g}\in \mathcal{N}^{\perp}$ and
	\begin{align*}
	\sum\limits_{|\alpha|\leq N}\|\nu^{-1}\mu^{-\frac{q'}{2}}\partial_{v}^{\alpha}\mathfrak{g}\|_{L_{v}^2}\leq\infty,
\end{align*}
	then it holds that
	\begin{align}\label{b2.7}
	|\mu^{-\frac{q}{2}}\FL^{-1}\mathfrak{g}(v)|\leq C\sum\limits_{|\alpha|\leq N}\|\partial_{v}^{\alpha}(\mu^{-\frac{q}{2}}\FL^{-1}\mathfrak{g})\|_{L_{v}^2}\leq C\sum\limits_{|\alpha|\leq N}\|\nu^{-1}\mu^{-\frac{q'}{2}}\partial_{v}^{\alpha}\mathfrak{g}\|_{L_{v}^2},\quad \text{for }v\in \R^3.
	\end{align}
	Here the constant $C=C(\rho,\fu,T,q',N)>0$.
\end{proposition}
Before proving Proposition \ref{thmA.1}, we first give some useful estimates on $K$.
\begin{lemma}\label{lve}
	Let $N\in\mathbb{N}$, $|\alpha|\leq N$, $0<q<1$, $-3<\kappa<0$. For any $r>0$, it holds that
	\begin{equation}\label{lve1}
		|\langle \nu^{-1}\mu^{-\frac{q}{2}}(\partial_v^\alpha K\mathfrak{h})^{1-\chi_r},\mu^{-\frac{q}{2}}\partial_v^\alpha \mathfrak{h}\rangle|\leq C\sum\limits_{0\leq \alpha'\leq \alpha}\|\partial_{v}^{\alpha'}\mathfrak{h}\|_{\nu}^2,
	\end{equation}
where $C$ depends on $\rho,\fu,T,r,q$ and $N$.
\end{lemma}

\noindent{\bf Proof.}
It suffices to prove
\begin{align*}
	|\langle \nu^{-1}\mu^{-\frac{q}{2}}H^{1-\chi_r},\mu^{-\frac{q}{2}}\partial_v^\alpha\mathfrak{h}\rangle|\leq C\sum\limits_{0\leq \alpha'\leq \alpha}\|\partial_{v}^{\alpha'}\mathfrak{h}\|_{\nu}^2,\quad H=I_1,I_2,J_1,J_2,
\end{align*} 	
where $I_1,I_2,J_1,J_2$ are the ones defined in \eqref{pk1}-\eqref{pk2}.
For $I_{1}$ and $I_2$, by using similar calculations as in \cite[Lemma 3.2]{Jiang-Luo-Tang}, we have
	\begin{equation}\label{a14}
		 |\langle\nu^{-1}\mu^{-\frac{q}{2}}I_{1}^{1-\chi_r},\mu^{-\frac{q}{2}}\partial_{v}^{\alpha}\mathfrak{h}\rangle|\leq C\|\partial_{v}^{\alpha}\mathfrak{h}\|_{\nu}^2,
	\end{equation}
and
	\begin{equation}\label{1.3}
		\begin{aligned}
			 |\langle\nu^{-1}\mu^{-\frac{q}{2}}I_2^{1-\chi_r},\mu^{-\frac{q}{2}}\partial_{v}^{\alpha}\mathfrak{h}\rangle|&\leq C\sum\limits_{0\leq \alpha'<\alpha}\left(\int_{\R^3}|\partial_{u}^{\alpha'}\mathfrak{h}(u)|^2\nu(u)\,du\right)^{\frac{1}{2}}\\
			&\qquad\qquad \times\left(\int_{\R^3}|\partial_{v}^{\alpha}\mathfrak{h}(v)|^2\nu(v)\,dv\right)^{\frac{1}{2}}\\
			&\leq
			\frac{1}{2}\|\partial_{v}^{\alpha}\mathfrak{h}\|_{\nu}^2+C\sum\limits_{0\leq \alpha'<\alpha}\|\partial_{v}^{\alpha'}\mathfrak{h}\|_{\nu}^2,
		\end{aligned}
	\end{equation}
where $C$ depends on ${\rho,\fu,T,r,q}$ and  we have chosen $0<\epsilon<\frac{1-q}{2}$ in the expression of $I_2^{1-\chi_r}$. For $|\langle \nu^{-1}\mu^{-\frac{q}{2}}(J_{i})^{1-\chi_r},\mu^{-\frac{q}{2}}\partial_v^{\alpha}\mathfrak{h}\rangle|$, using similar arguments as in \cite[Lemma 3.2]{Jiang-Luo-Tang}, one can get
	\begin{equation*}
		 |\langle\nu^{-1}\mu^{-\frac{q}{2}}J_i^{1-\chi_r},\mu^{-\frac{q}{2}}\partial_{v}^{\alpha}\mathfrak{h}\rangle|\leq C\sum\limits_{0\leq \alpha'\leq \alpha }\|\partial_{v}^{\alpha}\mathfrak{h}\|_{\nu}^2,\quad i=1,2,
	\end{equation*}
	which, together with \eqref{a14} and \eqref{1.3}, yields \eqref{lve1}. Therefore the proof of Lemma \ref{lve} is completed.
$\hfill\Box$

\begin{lemma}\label{hve}
	Let $N\in \mathbb{N}$, $|\alpha|\leq N$, $0 < q<q' < 1$, $- 3 < \kappa <0$. For any $r > 0$, there exists a constant $C>0$ such that the following estimates hold:
		\begin{equation}\label{K1-chi}
			\begin{aligned}
				\big| \langle \nu^{-1} \mu^{- \frac{q}{2}} (\partial_{v}^{\alpha}K_{1} \mathfrak{h})^{\chi_r}, \mu^{- \frac{q}{2}}\partial_{v}^{\alpha}\mathfrak{h} \rangle \big| &\leq C \exp  \big(- \frac{(1-q) r^2}{32 T} \big)\sum\limits_{0\leq \alpha'\leq  \alpha}\| \mu^{- \frac{q}{2}}\partial_{v}^{\alpha'}\mathfrak{h} \|^2_{L^2_v},
			\end{aligned}
		\end{equation}
	and
		\begin{align}\label{gamma-SP}
				| \langle \nu^{-1} \mu^{- \frac{q}{2}} (\partial_{v}^{\alpha}K_2\mathfrak{h})^{\chi_r} , \mu^{- \frac{q}{2}} \partial_{v}^{\alpha}\mathfrak{h} \rangle |& \leq C\Big\{ \frac{1}{1 + r}
				\sum\limits_{0\leq \alpha'\leq \alpha} \| \mu^{- \frac{q}{2}} \partial_{v}^{\alpha'}\mathfrak{h} \|^2_{L^2_v}  +  \exp (\frac{2 q r^2}{T})\sum\limits_{0\leq \alpha'\leq \alpha} \| \partial_{v}^{\alpha'}\mathfrak{h} \|^2_{\nu} \,\Big\}.
		\end{align}
The constant $C$ depends only on $\rho,\fu,T,q,N$.
\end{lemma}

\noindent{\bf Proof.}
	We divide it into several steps. We point out that the constants $C$ in the proof do not depend on $r$.

	{\it Step 1.} Estimates on $|\langle\nu^{-1}\mu^{-\frac{q}{2}}(\partial_{v}^\alpha K_1\mathfrak{h})^{\chi_r},\mu^{-\frac{q}{2}}\partial_{v}^{\alpha}\mathfrak{h}\rangle|$. Noting \eqref{pk1} and the definition of $\chi_r(s)$, we have
	\begin{align*}
		& \big| \langle \nu^{-1} \mu^{- \frac{q}{2}} I_{1}^{\chi_r}, \mu^{- \frac{q}{2}}\partial_{v}^{\alpha}\mathfrak{h} \rangle \big| \nonumber\\
		&\leq C \Big( \iint_{|u - v| \geq r} |\partial_{v}^{\alpha}\mathfrak{h} (v) \mu^{- \frac{q}{2}} (v)|^2 \mu^\frac{3(1-q)}{8} (v) \mu^\frac{1+q}{2} (u) |u - v|^\kappa d u d v \Big)^\frac{1}{2} \nonumber\\
		& \quad \times \Big( \iint_{|u - v| \geq r} |\partial_{u}^{\alpha}\mathfrak{h} (u) \mu^{- \frac{q}{2}} (u)|^2 \mu^\frac{3(1-q)}{8} (v) \mu^\frac{1+q}{2} (u) |u - v|^\kappa d u d v \Big)^\frac{1}{2} \,,
	\end{align*}
	where we have used $|\nu^{-1} (v) \mu^\frac{1-q}{8} (v)| \leq C(\rho, \fu, T,q) < \infty$. Then it is direct to obtain
	\begin{equation}\label{1.4}
		\big| \langle \nu^{-1} \mu^{- \frac{q}{2}} I_{1}^{\chi_r}, \mu^{- \frac{q}{2}}\partial_{v}^{\alpha}\mathfrak{h} \rangle \big|\leq C \exp{(-\frac{(1-q)r^2}{32T})}\|\mu^{-\frac{q}{2}}\partial_{v}^{\alpha}\mathfrak{h}\|_{L^2_v}^2.
	\end{equation}
Taking $0<\epsilon<\frac{1-q}{8}$, one has
	\begin{equation*}
		\begin{aligned}
			\big| \langle \nu^{-1} \mu^{- \frac{q}{2}} I_{2}^{\chi_r}, \mu^{- \frac{q}{2}}\partial_{v}^{\alpha}\mathfrak{h} \rangle \big|&\leq C \exp{(-\frac{(1-q)r^2}{32T})}\|\mu^{-\frac{q}{2}}\partial_{v}^{\alpha}\mathfrak{h}\|_{L^2_v}\sum\limits_{0\leq \alpha'<\alpha}\|\mu^{-\frac{q}{2}}\partial_{v}^{\alpha'}\mathfrak{h}\|_{L^2_v}\\
			&\leq C\exp{(-\frac{(1-q)r^2}{32T})}\sum\limits_{0\leq \alpha'\leq \alpha}\|\mu^{-\frac{q}{2}}\partial_{v}^{\alpha'}\mathfrak{h}\|_{L^2_v}^2,
		\end{aligned}
	\end{equation*}
	which, together with \eqref{1.4}, yields \eqref{K1-chi}.
	
	{\it Step 2.} Recall \eqref{pk2} and \cite[Lemma 3.3]{Jiang-Luo-Tang}, it is direct to have
	\begin{equation*}
		| \langle \nu^{-1}\mu^{- \frac{q}{2}} J_{1}^{\chi_r}, \mu^{- \frac{q}{2}}\partial_{v}^{\alpha}\mathfrak{h} \rangle|\leq C\Big\{\frac{1}{1+r}\|\mu^{-\frac{q}{2}}\partial_{v}^{\alpha}\mathfrak{h}\|_{L^2_v}^2+\exp(\frac{2qr^2}{T})\|\partial_{v}^{\alpha}\mathfrak{h}\|_{L^2_v}^2\Big\}.
	\end{equation*}
	For $| \langle \nu^{-1} \mu^{- \frac{q}{2}} J_{2}^{\chi_r}, \mu^{- \frac{q}{2}}\partial_{v}^{\alpha}\mathfrak{h} \rangle|$, taking $0<\epsilon<\frac{1-q}{
		2}$, we can obtain
	\begin{align*}
			&|\langle \mu^{-\frac{q}{2}}J_{2}^{\chi_r},\mu^{-\frac{q}{2}}\partial_{v}^{\alpha}\mathfrak{h}\rangle|\nonumber\\
			&\leq \frac{C}{1+r}\|\mu^{-\frac{q}{2}}\partial_{v}^{\alpha}\mathfrak{h}\|_{L^2_v}\sum\limits_{0\leq \alpha'<\alpha}\|\mu^{-\frac{q}{2}}\partial_{v}^{\alpha'}\mathfrak{h}\|_{L^2_v}\nonumber+C\exp(\frac{2qr^2}{T})\|\partial_{v}^{\alpha}f\|_{\nu}\sum\limits_{0\leq \alpha'<\alpha}\|\partial_{v}^{\alpha'}\mathfrak{h}\|_{\nu}\nonumber\\
			&\leq \frac{C}{1+r}\sum\limits_{0\leq \alpha'\leq\alpha}\|\mu^{-\frac{q}{2}}\partial_{v}^{\alpha'}\mathfrak{h}\|_{L^2_v}^2 +C\exp(\frac{2qr^2}{T})\sum\limits_{0\leq \alpha'\leq \alpha}\|\partial_{v}^{\alpha'}\mathfrak{h}\|_{\nu}^2.
	\end{align*}
	Therefore the proof of Lemma \ref{hve} is completed.
$\hfill\Box$

\begin{lemma}[Weighted hypocoercivity of  $\partial_{v}^{\alpha}\FL$]\label{HL}
	Let $N\in \mathbb{N}$, $|\alpha|\leq N$, $0<q<1$ and $-3<\kappa<0$. Then there is a constant $C=C(\rho,\fu,T,q,N)>0$ such that
	\begin{equation}\label{1.13}
		\langle\nu^{-1}\mu^{-\frac{q}{2}}\partial_{v}^{\alpha}(\FL \mathfrak{h}),\mu^{-\frac{q}{2}}\partial_{v}^{\alpha}\mathfrak{h}\rangle\geq \frac{1}{2}\|\mu^{-\frac{q}{2}}\partial_{v}^{\alpha}\mathfrak{h}\|_{L^2_v}^2-C\sum\limits_{0\leq \alpha'<\alpha}\|\mu^{-\frac{q}{2}}\partial_{v}^{\alpha'}\mathfrak{h}\|_{L^2_v}^2-C\sum\limits_{0\leq \alpha'\leq \alpha}\|\partial_{v}^{\alpha}\mathfrak{h}\|_{\nu}^2.
	\end{equation}

\end{lemma}

\noindent{\bf Proof.}	A direct calculation shows that $$
	\nu^{-1}(\partial_{v}^{\alpha}\FL \mathfrak{h})=\nu^{-1}\partial_{v}^{\alpha}(\nu \mathfrak{h})-\nu^{-1}(\partial_{v}^{\alpha}K)^{1-\chi_r}\mathfrak{h}-\nu^{-1}(\partial_{v}^{\alpha}K_1\mathfrak{h})^{\chi_r}+\nu^{-1}(\partial_{v}^{\alpha}K_2\mathfrak{h})^{\chi_r}.$$
	Noting  $\nu^{-1}(v)\partial_{v}^{\alpha}(\nu(v))\leq C_{N}$, we have
	\begin{equation}\label{1.17}
		\begin{aligned}
			\langle \nu^{-1}\mu^{-\frac{q}{2}}(\partial_{v}^{\alpha}(\nu \mathfrak{h})),\mu^{-\frac{q}{2}}\partial_{v}^{\alpha} \mathfrak{h}\rangle
			&\geq \frac{7}{8}\|\mu^{-\frac{q}{2}}\partial_{v}^{\alpha}\mathfrak{h}\|_{L^2_v}^2-C_{N}\sum\limits_{0\leq \alpha'<\alpha}\|\mu^{-\frac{q}{2}}\partial_{v}^{\alpha'}\mathfrak{h}\|_{L^2_v}^2.
		\end{aligned}
	\end{equation}
	Then it follows from \eqref{1.17} and Lemmas \ref{lve}--\ref{hve} that
	\begin{equation*}
		\begin{aligned}
			&\langle \nu^{-1}\mu^{-\frac{q}{2}}(\partial_{v}^{\alpha}\FL \mathfrak{h}),\mu^{-\frac{q}{2}}\partial_{v}^{\alpha}\mathfrak{h}\rangle\\
			&\geq \frac{7}{8}\|\mu^{-\frac{q}{2}}\partial_{v}^{\alpha}\mathfrak{h}\|_{L^2_v}^2-C\big\{\sum\limits_{0\leq \alpha'<\alpha}\|\mu^{-\frac{q}{2}}\partial_{v}^{\alpha'}\mathfrak{h}\|_{L^2_v}^2-\sum\limits_{0\leq \alpha'\leq\alpha}\|\partial_{v}^{\alpha'}\mathfrak{h}\|_{\nu}^2\big\}\\
			&\quad -C(\rho,\fu,T,q,N)[\exp(-\frac{(1-q)r^2}{32T})+\frac{1}{1+r}]\|\mu^{-\frac{q}{2}}\partial_{v}^{\alpha}\mathfrak{h}\|_{L^2_v}^2.
		\end{aligned}
	\end{equation*}
	Taking $r$ large enough, one gets \eqref{1.13}. Therefore the proof of Lemma \ref{HL} is completed.
$\hfill\Box$	
	
\medskip

For later use, we recall a result on the hypercoercivity 
in \cite{Guo2006}.
\begin{lemma}[{\cite{Guo2006}}]\label{HL1}
	Let $-3<\kappa<0$ and $|\alpha|\leq N$. Then there exists a constant $C(\rho,\fu,T,N)>0$ such that
	$$
	\langle\partial_{v}^{\alpha}(\FL \mathfrak{h}),\partial_{v}^{\alpha}\mathfrak{h}\rangle\geq \frac{1}{2}\|\partial_{v}^{\alpha}\mathfrak{h}\|_{\nu}^2-C\|\mathfrak{h}\|_{\nu}^2.
	$$
\end{lemma}

\noindent\textbf{Proof of Proposition \ref{thmA.1}}. Let $\mathfrak{g}\in\mathcal{N}^\perp$, we denote $\mathfrak{h}:=\FL ^{-1}\mathfrak{g}$, that is, $\FL \mathfrak{h}=\nu \mathfrak{h}-K\mathfrak{h}=\mathfrak{g}$. By Sobolev's embedding theorem, we have for $N\geq 2$ that
\begin{equation}\label{1.21}
	\begin{aligned}
		|\mu^{-\frac{q}{2}}\FL ^{-1}\mathfrak{g}|&\leq C\sum\limits_{|\alpha|\leq N}\|\partial_{v}^{\alpha}(\mu^{-\frac{q}{2}}\mathfrak{h})\|_{L^{2}_v}\\
		&=C\sum\limits_{|\alpha|\leq N}\|\mu^{-\frac{q}{2}}\partial_{v}^{\alpha}\mathfrak{h}\|_{L^2_v}+C\sum\limits_{|\alpha|\leq 2}\sum\limits_{0\leq \alpha'<\alpha}C_{\alpha}^{\alpha'}\|(\partial_{v}^{\alpha-\alpha'}\mu^{-\frac{q}{2}})(\partial_{v}^{\alpha'}\mathfrak{h})\|_{L^2_v}\\
		&\leq C\sum\limits_{|\alpha|\leq N}\|\mu^{-\frac{q'}{2}}\partial_{v}^{\alpha}\mathfrak{h}\|_{L^2_v},
	\end{aligned}
\end{equation}
where we have used the fact that
$$
\mu^{\frac{q'}{2}}\sum\limits_{0\leq \alpha'<2}(\partial_{v}^{\alpha-\alpha'}\mu^{-\frac{q}{2}})\leq C\qquad \text{for any }0<q<q'<1,
$$
in the last inequality.

It follows from Lemmas \ref{HL}-\ref{HL1} and $\|\mathfrak{h}\|_{\nu}^2\lesssim \langle\FL \mathfrak{h},\mathfrak{h}\rangle$ that
\begin{equation*}
	\begin{aligned}
		\|\mu^{-\frac{q'}{2}}\partial_{v}^{\alpha}\mathfrak{h}\|_{L^2_v}^2&\leq 2\langle\nu^{-1}\mu^{-\frac{q'}{2}}\partial_{v}^{\alpha}(\FL \mathfrak{h}),\mu^{-\frac{q'}{2}}\partial_{v}^{\alpha}\mathfrak{h}\rangle+C\sum\limits_{0\leq \alpha'<\alpha}\|\mu^{-\frac{q'}{2}}\partial_{v}^{\alpha'}\mathfrak{h}\|_{L^2_v}^2\\
		&\quad+C\|\partial_{v}^{\alpha}\mathfrak{h}\|_{\nu}^2\\
		&\leq 16\|\nu^{-1}\mu^{-\frac{q'}{2}}\partial_{v}^{\alpha}\mathfrak{g}\|_{L^2_v}^2+\frac{1}{4}\|\mu^{-\frac{q'}{2}}\partial_{v}^{\alpha}\mathfrak{h}\|_{L^2_v}^2+C\sum\limits_{0\leq \alpha'<\alpha}\|\mu^{-\frac{q'}{2}}\partial_{v}^{\alpha'}\mathfrak{h}\|_{L^2_v}^2\\
		&\quad +C\|\nu^{-1}\partial_{v}^{\alpha}\mathfrak{g}\|_{L^{2}_v}^2+\frac{1}{4}\|\mathfrak{h}\|_{\nu}^2,
	\end{aligned}
\end{equation*}
which, together with \eqref{1.21}, yields that
$$
|\mu^{-\frac{q}{2}}\FL ^{-1}\mathfrak{g}|\leq C\sum\limits_{|\alpha|\leq N}\|\mu^{-\frac{q'}{2}}\partial_{v}^{\alpha}\mathfrak{h}\|_{L^2_v}\leq C\sum\limits_{|\alpha|\leq N}\|\nu^{-1}\mu^{-\frac{q'}{2}}\partial_{v}^{\alpha}\mathfrak{g}\|_{L^2_v},
$$
where the constant depend only $\rho,\fu,T$ and $q'$. Thus the proof of Proposition \ref{thmA.1} is finished. $\hfill\square$

\begin{remark}\label{rmk2.6}
	Denote $\FL_0 \mathfrak{h}=\nu^0(v) \mathfrak{h}-K^0\mathfrak{h}= \FL \mathfrak{h}|_{x_3=0}$. Similarly, we  have $K^0\mathfrak{h}=K^{0,c}\mathfrak{h}+K^{0,m}\mathfrak{h}$. Define
	\begin{align}\label{D2.18-0}
		w_l(v)=(1+|v|^2)^{\f{l}{2}}\mu_0^{-\fa},
	\end{align}
then we can define $k_w^{0,m},k_w^{0,c}$ similarly as in section 2.1. It is obvious that  we can have similar results as in \eqref{D2.16}-\eqref{D2.20} for $K^{0,m}, K^{0,c},k_w^{0,m},k_w^{0,c}$.

\smallskip

For $\FL_0$, one also has
	\begin{align}\label{D2.21}
		\langle \FL_0\mathfrak{h}, \mathfrak{h}\rangle 
		\geq c_1\|\{\FI-\FP_0\}\mathfrak{h}\|_\nu^2,
	\end{align}
	since $\nu^0\cong \nu\cong  (1+|v|)^\k$.
	$\FL_0^{-1}$ can be defined as
	\begin{align*}
		(\FL_0|_{\mathcal{N}_0^\perp})^{-1}:\mathcal{N}_0^\perp\rightarrow \mathcal{N}_0^{\perp}.
	\end{align*}
	Let $\mathfrak{g}\in\mathcal{N}_0^{\perp}$, from Proposition \ref{thmA.1}, it is direct to know that
	\begin{align}\label{b2.8}
		|\mu_0^{-\frac{q}{2}}\FL_0^{-1}\mathfrak{g}(v)|\leq C\sum\limits_{|\alpha|\leq N}\|\partial_{v}^{\alpha}(\mu_0^{-\frac{q}{2}}\FL_0^{-1}\mathfrak{g})\|_{L_{v}^2}\leq C\sum\limits_{|\alpha|\leq N}\|(\nu^0)^{-1}\mu_0^{-\frac{q'}{2}}\partial_{v}^{\alpha}\mathfrak{g}\|_{L_{v}^2},\quad \text{for }v\in \R^3.
	\end{align}
	where the constant $C=C(\rho^0,\fu^0,T^0,q',N)>0$. These estimates will be used to study the viscous and Knudsen boundary layers.	
\end{remark}
\begin{remark}
We point out that all above results for soft potentials in this section are also valid for hard potentials. The proofs are very similar.
\end{remark}

\medskip

	\section{Existence of a Steady Linear Boltzmann Equation}\label{section2}
 To construct the Knudsen layer solutions, we study equation:
	\begin{align}\label{D2.1}
		\begin{cases}
			v_3\partial_\eta \mathfrak{f}+\FL_0 \mathfrak{f}=S, \quad(\eta,v)\in[0,\infty)\times\R^3,\\		 \mathfrak{f}(0,v)|_{v_3>0}=\mathfrak{f}(0,R_x v)+f_b(v),\\
			\lim\limits_{\eta\rightarrow\infty}\mathfrak{f}(\eta,v)=0.\\
		\end{cases}
	\end{align}
where $S$ is a given function and $f_b(v)$ is defined only for $v_3< 0$, and we always assume that it is extended to be $0$ for $v_3>0$. For soft potentials, there has not been work for Knudsen layer solutions with specular reflection boundary.

\smallskip

	\begin{theorem}\label{thm1.1}
Recall $w_l$ in \eqref{D2.18-0}. Assume $l>2,\ 0\leq\fa<\f12$,
\begin{align}\label{D2.3}
\begin{split}
\int_{\R^3}(1,v_1-\fu^0_1,v_2-\fu^0_2,|v-\fu^0|^2)\sqrt{\mu_0} Sdv=0,\\	 \int_{\R^3}(1,v_1-\fu^0_1,v_2-\fu^0_2,|v-\fu^0|^2)v_3\sqrt{\mu_0} f_bdv=0,
\end{split}
\end{align}
and
\begin{align}\label{D2.4}
\begin{split}
\|(1+\eta)^{q_k}w_{l+4k+4}S\|_{L^\infty_{\eta,v}}<\infty,\\	\|w_{l+4k+5}f_b\|_{L^\infty_{v}}<\infty,
\end{split}
\qquad \text{for} \quad k\in \mathbb{N}_+,\ q_k>\max\{3,k+\f32\},
\end{align}
then there exists a  unique solution $\mathfrak{f}$ of \eqref{D2.1} such that
\begin{align}\label{D2.5}
&\|(1+\eta)^k w_l\mathfrak{f}\|_{L^\infty_{\eta,v}} +|(1+\eta)^kw_l\mathfrak{f}(0,\cdot)|_{L^\infty_{v}}\nonumber\\
&\leq C\Big\{\|(1+\eta)^{q_k}w_{l+4k+4}S\|_{L^\infty_{\eta,v}}+\|w_{l+4k+5}f_b\|_{L^\infty_v}\Big\},
\end{align}
where $C>0$ is a positive constant.
\end{theorem}

\begin{remark}
As indicated in \cite{GKTT}, in general, it is hard to obtain  the normal derivatives estimates for the boundary value problem \eqref{D2.1}. Fortunately, it is easy to obtain the tangential and time derivatives estimates for the solution of \eqref{D2.1}, i.e.,
\begin{align}\label{D2.7}
& \sum_{i+j\leq r} \|(1+\eta)^kw_{l} \partial_t^{i}\nabla_\sp^{j}\mathfrak{f}(t,x_\sp,\cdot,\cdot)\|_{L^\infty_{\eta,v}}
+\|w_{l} \partial_t^{i}\nabla_\sp^{j}\mathfrak{f}(t,x_\sp,0,\cdot)\|_{L^\infty_v}\nonumber\\
&\leq C  \sum_{i+j\leq r}\Big\{ \|w_{l+4k+5} \partial_t^{i}\nabla_\sp^{j}f_b(t,x_\sp,\cdot)\|_{L^\infty_v}+
\| (1+\eta)^{q_{k}}w_{l+4k+4} \partial_t^{i}\nabla_\sp^{j}S\|_{L^\infty_{\eta,v}}\Big\},
\end{align}
provided the right hand side of \eqref{D2.7} is bounded and $q_{k}>\max\{3,k+\f{3}{2}\}$. We point out that such an  estimate \eqref{D2.7} is enough for us to establish the Hilbert expansion. To prove the estimate \eqref{D2.7}, we study  the equation of $\partial_t^{i}\nabla_\sp^{j}(\sqrt{\mu_0}\mathfrak{f})$. It is direct to check that the new source term and boundary perturbation term satisfy the solvability conditions in Theorem \ref{thm1.1}, hence one can obtain the estimate for $\partial_t^{i}\nabla_\sp^{j}(\sqrt{\mu_0}\mathfrak{f})$ by applying Theorem \ref{thm1.1}, therefore  \eqref{D2.7} follows immediately.
		
		Moreover, taking $L^\infty_{x_\sp}\cap L^2_{x_\sp}$ over \eqref{D2.7}, one obtains
		\begin{align}\label{D2.8}
			& \sum_{i+j\leq r} \sup_{t\in[0,\tau^\d]} \Big\{\|(1+\eta)^kw_{l}\partial_t^{i}\nabla_\sp^{j}\mathfrak{f}(t)\|_{L^\infty_{x_\sp,\eta,v}\cap L^2_{x_\sp}L^\infty_{\eta,v}}+\|w_{l}  \partial_t^{i}\nabla_\sp^{j}\mathfrak{f}(t,\cdot,0,\cdot)\|_{L^\infty_{x_\sp,v}\cap L^2_{x_\sp}L^\infty_{v}}\Big\}\nonumber\\
			&\leq C \sup_{t\in[0,\tau^\d]}\Big\{ \sum_{i+j\leq r} \Big\{ \|w_{l+4k+5}\partial_t^{i}\nabla_\sp^{j}f_b(t)\|_{L^\infty_{x_\sp,v}\cap L^2_{x_\sp}L^\infty_{v}} \nonumber\\
			&\qquad\qquad +
			\sum_{i+j\leq r} \| (1+\eta)^{q_k}w_{l+4k+4} \partial_t^{i}\nabla_\sp^{j}S(t)\|_{L^\infty_{x_\sp,\eta,v}\cap L^2_{x_\sp}L^\infty_{\eta,v}}\Big\}\Big\},\quad q_k>\max\{3,k+\f{3}{2}\}.
		\end{align}
	\end{remark}	
	 Let $\Upsilon(\eta)$ be a monotonic smooth cut-off function
	\begin{align*}
		\Upsilon(\eta)\equiv 1,\ \mbox{for} \ \eta\in[0,1],\quad \mbox{and}\quad \Upsilon(\eta)\equiv 0,\ \mbox{for}\ \eta\in[2,+\infty).
	\end{align*}
	Define
	\begin{align*}
		f(x,v):=\mathfrak{f}(x,v)-\Upsilon(\eta)f_b(v), \quad x=(x_{\sp},\eta,v)
	\end{align*}
	then \eqref{D2.1} is rewritten as
	\begin{align}\label{D2.2}
		\begin{cases}
			v_3\partial_\eta f+\FL_0 f=g:=S-v_3\partial_\eta\Upsilon(\eta)f_b(v)-\Upsilon(\eta)\FL_0 f_b, \quad(\eta,v)\in[0,\infty)\times\R^3,\\
			f(0,v)|_{v_3>0}=f(0,R_\eta v),\\
			\lim\limits_{\eta\rightarrow\infty}f(\eta,v)=0,\\
		\end{cases}
	\end{align}
	where $x_{\sp}$ is regarded as parameters. The conditions \eqref{D2.3} deduces that
	\begin{align}\label{D2.2-0}
			\int_{\R^3}(1,v_1-\fu^0_1,v_2-\fu^0_2,|v-\fu^0|^2)\sqrt{\mu_0}g \ dv=0.
	\end{align}

	Define the viscosity and thermal conductivity coefficients by
	\begin{align}\label{D2.24}
		\begin{split}
			\mu(T^0)&:=T^0\langle \mathcal{A}_{31}^0,\  \mathbf{L}_0^{-1} \mathcal{A}_{31}^0 \rangle\equiv T^0\langle\mathcal{A}_{ij}^0,\  \mathbf{L}_0^{-1} \mathcal{A}_{ij}^0 \rangle,\quad \forall i\neq j,\\[1.5mm]
			\kappa(T^0)&:=\frac23T^0\langle \mathcal{B}_{3}^0,\  \mathbf{L}_0^{-1} \mathcal{B}_{3}^0\rangle \equiv \frac23T^0 \langle\mathcal{B}_{i}^0,\  \mathbf{L}_0^{-1} \mathcal{B}_{i}^0 \rangle,
		\end{split}
	\end{align}
where $ i,j=1,2,3$ and $\mathcal{A}_{ij}^0,$ $\mathcal{B}_{i}^0$ are
\begin{align}\label{D2.23}
\begin{split}
\mathcal{A}_{ij}^0&:=\left\{\f{(v_i-\fu_i^0)(v_j-\fu_j^0)}{T^0}-\d_{ij}\f{|v-\fu^0|^2}{3T^0}\right\}\sqrt{\mu_0},\\
\mathcal{B}_{i}^0&:=\f{v_i-\fu_i^0}{2\sqrt{T^0}}\left(\f{|v-\fu^0|^2}{T^0}-5\right)\sqrt{\mu_0}.
\end{split}
\end{align}
Using Lemma 4.4 in \cite{Bardos-2}, one has $\langle T^0\mathcal{A}_{33}^0, \mathbf{L}_0^{-1} \mathcal{A}_{33}^0\rangle=\frac43\mu(T^0) $.
		
\subsection{Approximate solutions and uniform estimate}
This section is devoted to the existence result for the linearized problem \eqref{D2.2}. To prove the existence of solution, we first consider a truncated approximate problem with penalized term:
	\begin{align}\label{D2.25}
		\begin{cases}
			\d f^\d+ v_3 \partial_\eta f^\d+\FL_0 f^\d=g,\\
			f^\d(\eta,v)|_{\gamma_-}=f^\d(\eta,R_\eta v),\\
		\end{cases}
		( \eta,v)\in \Omega_d\times \R^3,
	\end{align}
	where $\Omega_d=(0,d),\, d\geq1$ and $\delta\in(0,1]$. We define
	\begin{equation*}
		h^\d(\eta ,v):=w_l(v) f^\d(\eta ,v), 
	\end{equation*}
	then 
	\eqref{D2.25} can be rewritten as
	\begin{align}\label{D2.26}
		\begin{cases}
			\dis \d h^\d+v_3\partial_\eta  h^\d+\nu^0(v) h^\d=K^0_{w_l} h^\d+w_lg,\\[2mm]
			\dis h^\d(\eta ,v)|_{\gamma_-}=h^\d(\eta ,R_\eta v),
		\end{cases}
	\end{align}
	where $K^0_{w_l}h=w_lK^0(\f{h}{w_l}).$ Then it is clear that
	\begin{align}\label{D2.27}
		K^0_{w_l}h(v)=\int_{\R^3} k^0_{w_l}(v,u) h(u) du\quad\mbox{with}\quad
		k^0_{w_l}(v,u)=w_l(v) k^0(v,u) w_l(u)^{-1}.
	\end{align}
	
	
	For the approximate problem  \eqref{D2.26},  the most difficult part is to obtain the $L^\infty_{\eta ,v}$-bound.
	Motivated by \cite{DHWZ}, multiplying \eqref{D2.26} by $(1+|v|^2)^\f{|\kappa|}{2}$, one gets
	\begin{align}\label{D2.28}
		\begin{cases}
			\dis (\nu^0(v)+\d) (1+|v|^2)^{\f{|\kappa|}{2}}h^\d+v_3(1+|v|^2)^{\f{|\k|}{2}}\partial_\eta  h^\d=(1+|v|^2)^{\f{|\kappa|}{2}}K^0_{w_l} h^\d+(1+|v|^2)^{\f{|\kappa|}{2}}w_lg,\\[2mm]
			\dis h^\d(\eta ,v)|_{\gamma_-}=h^\d(\eta ,R_\eta v).
		\end{cases}
	\end{align}
		Denote $\hat{\nu}_\d=(\nu^0(v)+\d)(1+|v|^2)^{\f{|\kappa|}{2}},\hat{\nu}=\nu^0(v)(1+|v|^2)^{\f{|\kappa|}{2}},\hat{v}_3=v_3(1+|v|^2)^{\f{|\kappa|}{2}}$, then $\eqref{D2.28}_1$ becomes
		\begin{align*}
				\dis \hat{\nu}_\d h^\d+\hat{v}_3\partial_\eta  h^\d=(1+|v|^2)^{\f{|\kappa|}{2}}K^0_{w_l} h^\d+(1+|v|^2)^{\f{|\kappa|}{2}}w_lg.
		\end{align*}
		For given $(t,\eta ,v),$ let $[X(s),V(s)]$ 
		be the speeded backward characteristics for \eqref{D2.28}. Then $[X(s),V(s)]$ is determined by
		\begin{align*}
			\begin{cases}
				\dis \frac{dX(s)}{ds}=\hat{V}_3(s):=V_3(s)(1+|V|^2)^{\f{|\kappa|}{2}},
				\quad \frac{dV(s)}{ds}=0,\\[2mm]
				[X(t),V(t)]=[\eta,v],
			\end{cases}
		\end{align*}
		which yields that
		\begin{equation*}
			[X(s),V(s)]=[X(s;t,\eta,v),V(s;t,\eta,v)]=[\eta-(t-s)\hat{v}_3,v].
		\end{equation*}
	Now for each $(\eta,v)$ with $\eta\in \bar{\Omega}_d$ and $v_3\neq 0,$ we define its {backward exit time} $t_{\mathbf{b}}(\eta,v)\geq 0$ to be the last moment at which the
	back-time straight line $[X({-\tau};0,\eta,v),V({-\tau};0,\eta,v)]$ remains in $\bar{\Omega}$:
	\begin{equation*}
		t_{\mathbf{b}}(\eta,v)={\sup\{s \geq 0:\eta-\tau \hat{v}_3\in\bar{\Omega}_d\text{ for }0\leq \tau\leq s\}.}
	\end{equation*}
	We also define the last position
	\begin{equation*}
		\eta_{\mathbf{b}}(\eta,v)=\eta(t_{\mathbf{b}})=\eta-t_{\mathbf{b}}(\eta,v)\hat{v}_3\in \partial \Omega_d .
	\end{equation*}
	It is obvious that $X(s)$, $t_{\mathbf{b}}(\eta,v)$ and $\eta_{\mathbf{b}}(x,v)$ are independent of  the horizontal velocity $v_{\sp}:=(v_1,v_2)$.\\
	
	Let $\eta\in \bar{\Omega}_d$, $(\eta,v)\notin \gamma _{0}\cup \g_{-}$ and
	$
	(t_{0},\eta_{0},v_{0})=(t,\eta,v)$, we inductively define
	\begin{equation*}
		(t_{k+1},\eta_{k+1},v_{k+1})=(t_{k}-t_{\mathbf{b}}(\eta_{k},v_{k}),\eta_{\mathbf{b}}(\eta_{k},v_{k}), R_{\eta_{k+1}}v_{k}),\quad k\geq1,
	\end{equation*}
	and the back-time cycle as
	\begin{equation}\label{D2.29}
		\left\{\begin{aligned}
			X_{cl}(s;t,\eta ,v)&=\sum_{k}\Fi_{(t_{k+1},t_{k}]}(s)\{\eta _{k}-\hat{v}_{k,3}(t_{k}-s)\},\\[1.5mm]
			V_{cl}(s;t,\eta ,v)&=\sum_{k}\Fi_{(t_{k+1},t_{k}]}(s)v_{k}.
		\end{aligned}\right.
	\end{equation}
	Clearly, for $k\geq 1$ and $(\eta ,v)\notin \g_0\cup \g_-$, it holds that
	\begin{equation}\label{D2.30}
		\begin{split}
			&\eta _k=\frac{1-(-1)^k}{2} \eta _1+\frac{1+(-1)^k}{2} \eta _2,\quad  v_{k,\sp}=v_{0,\sp}, \quad v_{k,3}=(-1)^{k} v_{0,3},\\
			&t_k-t_{k+1}=t_1-t_2=\frac{d}{|\hat{v}_{0,3}|}>0,\quad \nu^0(v)\equiv \nu^0(v_{k}).
		\end{split}
	\end{equation}

Now we are in a position to  construct solutions to \eqref{D2.25} or equivalently \eqref{D2.26}. We first present a useful $L^\infty$ {\it a priori} uniform estimate which will be used frequently.
\begin{lemma}\label{lem3.4}
For any given $\lambda\in[0,1]$, let $f^\lambda$ be the solution of the following system:
\begin{align}\label{D2.42}
\begin{cases}
\delta f^\lambda+v_3\partial_\eta f^\lambda+\nu^0(v) f^\lambda-\lambda K^0f^\lambda=g,\\[2mm]
f^\lambda(\eta,v)|_{\gamma_-}=(1-\frac1n) f^\lambda(\eta,R_\eta v)+r(\eta,R_\eta v),
\end{cases}
\end{align}
	where $n>1$ is an integer and $g,r$ are given. Assume $\| w_l f^\lambda\|_{L^\infty_{\eta,v}}+| w_l f^\lambda|_{L^\infty(\gamma_+)}<\infty$, $l>2$, then it holds that
\begin{align}\label{D2.59}
	\|w_lf^\lambda\|_{L^\infty_{\eta,v}}+|w_lf^\lambda|_{L^\infty(\gamma_+)}\leq C\|\sqrt{\nu^0}f^\lambda\|_{L^2_{\eta,v}}+C \{\|(\nu^0)^{-1}w_lg\|_{L^\infty_{\eta,v}}+|w_{l+4}r|_{L^\infty(\gamma_+)}\}.
\end{align}
We point out the constant  $C>0$ is independent of $\lambda$, $d$ and $n$.
\end{lemma}

\begin{remark}
For hard potentials, similar uniform estimate has been obtained in \cite{Jiang-Wang}. For soft potentials, since the effect of collision frequency is weak, i.e., $\nu^{0}(v)=(1+|v|)^{\kappa} \to 0$ as $|v|\to \infty$, we have to be more careful. In fact, one has to loss some weight to control the boundary perturbation $r$, see \eqref{D2.59}
\end{remark}

\noindent{\bf Proof.}
Denote $h^\lambda:= w_l f^\lambda$, then it holds that
		\begin{align*}
			\begin{cases}
				\dis	\hat{\nu}_\d h^\lambda+\f{d h^\lambda}{ds}=(1+|v|^2)^{\f{|\kappa|}{2}}\lambda K^0_{w_l} h^\lambda+(1+|v|^2)^{\f{|\kappa|}{2}}w_lg,\\
				h^\lambda(\eta,v)|_{\gamma_-}=(1-\f1n)h^\lambda(\eta,R_\eta v)+w_lr(\eta,R_\eta v).
			\end{cases}
		\end{align*}
		Integrating along the characteristic line, one gets
		\begin{align}\label{D2.47}
			 h^\lambda(\eta,v)&=(1-\f1n)^kh^\lambda(\eta_k,v_k)e^{-\hat{\nu}_\d(v)(t-t_k)}+\lambda\sum\limits_{i=0}^{k-1}(1-\f1n)^i\int_{t_{i+1}}^{t_i}e^{-\hat{\nu}_\d(v)(t-s)}(1+|v|^2)^{\f{|\kappa|}{2}}K^0_{w_l}h^\lambda ds\nonumber\\
			 &+\sum\limits_{i=0}^{k-1}(1-\f1n)^i\int_{t_{i+1}}^{t_i}e^{-\hat{\nu}_\d(v)(t-s)}(1+|v|^2)^{\f{|\kappa|}{2}}w_lgds
			+\sum\limits_{i=0}^{k-1}(1-\f1n)^i (w_l r)(\eta_i,v_{i+1})e^{-\hat{\nu}_\d(v)(t-t_i)}\nonumber\\
			&=: I_1+I_2+I_3+I_4.
		\end{align}
	
\smallskip

Taking $k=\tilde{k}_0|v_3|(1+|v|^2)^{\f{|\kappa|}{2}}$ with  $\tilde{k}_0\gg1$ chosen later. Then it holds that
\begin{align}\label{D.21}
I_1\leq e^{-\nu_0(k-1)t_{\mathbf{b}}}\|h^\lambda\|_{L^\infty_{\eta,v}}\leq e^{-\f{1}{2}\nu_0\tilde{k}_0d} \|h^\lambda\|_{L^\infty_{\eta,v}},
\end{align}
where $\nu_0>0$ is a constant depending on $\rho^0,\fu^0,T^0$. It is obvious that
\begin{align}\label{D.22}
I_3\leq \|(\nu^0)^{-1}w_lg\|_{L^\infty_{\eta,v}}.
\end{align}
For $I_4$, noting $|v_{i}|= |v|$, one has
\begin{align}\label{D.23}
I_4\leq k(1+|v|)^{-4}|w_{l+4}r|_{L^\infty(\gamma_+)}\leq C|w_{l+4}r|_{L^\infty(\gamma_+)}.
\end{align}

To estimate $I_2$, we divide it into two parts:
\begin{align}\label{D2.51}
&\sum_{i=0}^{k-1} (1-\f1n)^i\int_{t_{i+1}}^{t_i}e^{-\hat{\nu}_\d(v) (t-s)}  \lambda (1+|v|^2)^{\f{|\kappa|}{2}}K^0_{w_l} h^\lambda(X_{cl}(s),v_{i}) ds\nonumber\\
&\leq\sum_{i=0}^{k-1} \int_{t_{i+1}}^{t_i}e^{-\hat{\nu}_\d(v) (t-s)}(1+|v|^2)^{\f{|\kappa|}{2}}|K^{0,c}_{w_l} h^\lambda(X_{cl}(s),v_{i})| ds\nonumber\\
&+\sum_{i=0}^{k-1} \int_{t_{i+1}}^{t_i}e^{-\hat{\nu}_\d(v) (t-s)}  (1+|v|^2)^{\f{|\kappa|}{2}}|K^{0,m}_{w_l} h^\lambda(X_{cl}(s),v_{i}) |ds.
\end{align}
For the second term on the RHS of \eqref{D2.51}, one has  from  \eqref{D2.16} that
		\begin{align}\label{D.25}
			\sum_{i=0}^{k-1} \int_{t_{i+1}}^{t_i}e^{-\hat{\nu}_\d(v) (t-s)}   (1+|v|^2)^{\f{|\kappa|}{2}}|K^{0,m}_{w_l} h^\lambda(X_{cl}(s),v_{i})| ds\leq C m^{3+\kappa} e^{-\frac{|v|^2}{20}}\|h^\lambda\|_{L^\infty_{\eta ,v}}.	\end{align}
		For the first term on the RHS of \eqref{D2.51}, we use \eqref{D2.47} again to obtain
		\begin{align}\label{D2.53}
			&\sum_{i=0}^{k-1} \int_{t_{i+1}}^{t_i}e^{-\hat{\nu}_\d(v) (t-s)}   (1+|v|^2)^{\f{|\kappa|}{2}}|K^{0,c}_{w_l} h^\lambda(X_{cl}(s),v_{i})| ds\nonumber\\
			&=\sum_{i=0}^{k-1} \int_{t_{i+1}}^{t_i}e^{-\hat{\nu}_\d(v) (t-s)} (1+|v|^2)^{\f{|\kappa|}{2}}\Big| \int_{\R^3} k^{0,c}_{w_l}(v_i,v') h^\lambda(X_{cl}(s),v') dv'\Big|ds\nonumber\\
			&\leq  \sum_{i=0}^{k-1} \int_{t_{i+1}}^{t_i}e^{-\hat{\nu}_\d(v) (t-s)} (1+|v|^2)^{\f{|\kappa|}{2}} \int_{\R^3} |k^{0,c}_{w_l}(v_i,v')|\times (1+|v'|^2)^{\f{|\kappa|}{2}} dv'ds\nonumber\\
			&\qquad\times \sum_{j=0}^{k'-1} \int_{t'_{j+1}}^{t'_j}e^{-\hat{\nu}_\d(v') (s-s_1)}\int_{\R^3} |k^{0,c}_{w_l}(v'_j,v'') h^\lambda(X_{cl}'(s_1),v'')|dv'' ds_1 \nonumber\\
			&\quad +C \big(m^{3+\kappa} + m^{\kappa-1} e^{-\f{1}{2}\nu_0\tilde{k}_0d}\big) \|h^\lambda\|_{L^\infty_{\eta,v}} + Cm^{\kappa-1} \big\{    \|(\nu^0)^{-1}w_lg\|_{L^\infty_{\eta,v}}+ |w_{i+4}r|_{L^\infty(\gamma_+)}\big\},
		\end{align}
		where we have used \eqref{D.21}-\eqref{D.23}, \eqref{D.25} and \eqref{D2.19}-\eqref{D2.20},  and denoted $X_{cl}'(s_1)=X_{cl}(s_1;s,X_{cl}(s),v')$, and $t'_j,v'_j$ are the corresponding times and velocities for specular cycles. Here  $k'=\tilde{k}_0|v_3'|(1+|v'|^2)^{\f{|\kappa|}{2}}$.

		\smallskip
		
		For the first term on RHS of \eqref{D2.53}, we divide the proof into several cases. \\
		\noindent{\it Case 1. $|v|\geq N$.} Using \eqref{D2.19}, the first term on the RHS of \eqref{D2.53} is bounded by
		\begin{align}\label{D.27}
			&Cm^{\kappa-1}\sum_{i=0}^{k-1} \int_{t_{i+1}}^{t_i}e^{-\hat{\nu}_\d(v) (t-s)} (1+|v|^2)^{\f{|\kappa|}{2}} \int_{\R^3} |k^{0,c}_{w_l}(v_i,v')| (1+|v'|)^{-2}dv'ds \cdot \|h^\lambda\|_{L^\infty_{\eta,v}}\nonumber\\
			&\leq Cm^{2(\kappa-1)}(1+|v|)^{-2}\|h^\lambda\|_{L^\infty_{\eta,v}}\leq C\f{m^{2(\kappa-1)}}{N^2} \|h^\lambda\|_{L^\infty_{\eta,v}},
		\end{align}
		where we have used the fact $|v|\equiv |v_i|$ for $i=0,1,\cdots$. It is important that the constant in \eqref{D.27} is independent of $k$.
		
		\smallskip
		
		\noindent{\it Case 2. $|v|\leq N, |v'|\geq 2N$ or $|v'|\leq 2N, |v''|\geq 3N$.}  Noting $|v_i|=|v|$ and $|v'_j|=|v'|$, we get either $|v_i-v'|\geq N$ or $|v'_j-v''|\geq N$, then either one of the following is valid for some small positive constant $0<c_2\leq \frac{1}{32}$:
		\begin{equation}
			\begin{split}
				\dis |k^{0,c}_{w_l}(v_i,v')|&\leq e^{-c_2N^{2}}|k^{0,c}_{w_l}(v_i,,v') \exp{\left(  c_2|v_i-v'|^{2}\right)}|,\\[1mm]
				\dis |k^{0,c}_{w_l}(v'_j,v'')|&\leq e^{-c_2N^{2}}|k^{0,c}_{w_l}(v_j',,v'')\exp{\left(  c_2|v'_j-v'|^{2}\right)}|,
			\end{split}
		\end{equation}
		which, together with \eqref{D2.19}, yields that
		\begin{align}\label{D2.54}
			&\sum_{i=0}^{k-1} \int_{t_{i+1}}^{t_i}e^{-\hat{\nu}_\d(v) (t-s)}  \left\{\iint_{|v|\leq N, |v'|\geq 2N} + \iint_{|v'|\leq 2N, |v''|\geq 3N}  \right\} (\cdots)dv'' ds_1 dv'ds\nonumber\\
			&\leq Cm^{\kappa-1}e^{-c_2 N^2}\|h^\lambda\|_{L^\infty_{\eta,v}} \leq  \frac{Cm^{\kappa-1}}{N}\|h^\lambda\|_{L^\infty_{\eta,v}}.
		\end{align}
	We also point out that the constant in \eqref{D2.54} is independent of $k$.
		
		\smallskip
		
		\noindent{\it Case 3. $|v|\leq N, |v'|\leq 2N$, $|v''|\leq 3N$}. We denote $\mathcal{D}=\{|v|\leq N,\,|v'|\leq 2N,\,|v''|\leq 3N\}$. Noting $\hat{\nu}_{\d}(v)\geq \nu_0$,  the corresponding part is bounded by
		\begin{align}
			&\sum_{i=0}^{k-1} \int_{t_{i+1}}^{t_i}e^{-\nu_0 (t-s)}  \iint_{\mathcal{D}} |k^{0,c}_{w_l}(v_i,v')k^{0,c}_{w_l}(v'_j,v'')|(1+|v|^2)^{\f{|\kappa|}{2}}(1+|v'|^2)^{\f{|\kappa|}{2}} dv'' dv'ds\nonumber\\
			&\qquad\times \sum_{j=0}^{k'-1} \left(\int_{t'_{j}-\frac{1}{N^6}}^{t'_j}+\int_{t'_{j+1}}^{t'_j-\frac{1}{N^6}}\right)e^{-\nu_0 (s-s_1)} | h^\lambda(X_{cl}'(s_1),v'')| ds_1=:P_1+P_2. \nonumber
		\end{align}
	
			For $P_1$, noting $|v'|\leq 2N$, one has
	\begin{align*}
		P_1&\leq C\f{k' m^{2\kappa-2}}{N^6}\|h^\lambda\|_{L^\infty_{\eta,v}}\leq C\f{\tilde{k}_0m^{2(\kappa-1)}N^4}{N^6}\|h^\lambda\|_{L^\infty_{\eta,v}}\leq  \f{Cm^{2(\kappa-1)}}{N^2}\|h^\lambda\|_{L^\infty_{\eta,v}}.
	\end{align*}

		For $P_2$,  a direct calculation shows
		\begin{align}\label{D2.56}
			P_2&\leq \sum\limits_{i=0}^{k-1}\sum\limits_{j=0}^{k'-1}\int_{t_{i+1}}^{t_i}\int_{t_{j+1}'}^{t_j'-\f{1}{N^6}}\iint_{\mathcal{D}}|k^{0,c}_{w_l}(v_i,v')k^{0,c}_{w_l}(v_j',v'')|(1+|v|^2)^{\f{|\kappa|}{2}}(1+|v'|^2)^{\f{|\kappa|}{2}}
				\nonumber\\
				&\qquad\qquad\times|e^{-\nu_0(t-s_1)}h^\lambda(X_{cl}'(s_1),v'')|dv''dv'ds_1ds\nonumber\\
				&\leq C_N\sum\limits_{i=0}^{k-1}\sum\limits_{j=0}^{k'-1}\int_{t_{i+1}}^{t_i}\int_{t_{j+1}'}^{t_j'-\f{1}{N^6}}e^{-\nu_0(t-s_1)}\Big[\iint_{\mathcal{D}}\nu^0(v'')|f^{\lambda }(X_{cl}'(s_1),v'')|^2dv''dv'\Big]^{\f12},
		\end{align}
		where we used the fact
		\begin{align}\nonumber
			\iint_{\mathcal{D}}  |k^{0,c}_{w_l}(v_i,v')k^{0,c}_{w_l}(v'_j,v'')|^2(1+|v|^2)^{|\kappa|}(1+|v'|^2)^{|\kappa|}{w_l}^2(v'') (\nu^0)^{-1}(v'') dv'dv''\leq C_N.
		\end{align}

		Define $y:=\eta'_j-\hat{v}'_{j,3} (t'_j-s_1)=X_{cl}'$. We have  $\eta'_j=0\,  \mbox{or}\, d$ and $\hat{v}'_{j,3}=(-1)^{j} \hat{v}'_{0,3}$. For $t'_j=t'_j(s_1;s, X_{cl}(s),v')$, it holds that
		\begin{align*}
			s-t'_j=
			\begin{cases}
				\displaystyle \frac{X_{cl}(s)}{|\hat{v}'_{0,3}|}+(j-1) \frac{d}{|\hat{v}'_{0,3}|},\quad \mbox{for}\  v'_{0,3}>0,\\[3mm]
				\displaystyle \frac{d-X_{cl}(s)}{|\hat{v}'_{0,3}|}+(j-1) \frac{d}{|\hat{v}'_{0,3}|},\quad \mbox{for}\  v'_{0,3}<0,
			\end{cases}
		\end{align*}
		which yields that
		\begin{align}\nonumber
			y=
			\begin{cases}
				\eta '_j-(-1)^{j}\Big\{\hat{v}'_{0,3}(s-s_1)-[X_{cl}(s)+(j-1)d] \Big\},\quad  \mbox{for}\  v'_{0,3}>0,\\[2mm]
				\eta '_j-(-1)^{j} \Big\{\hat{v}'_{0,3}(s-s_1)+[j d-X_{cl}(s)]\Big\},\quad  \mbox{for}\  v'_{0,3}<0.
			\end{cases}
		\end{align}
		Since  $\eta '_j=0\,  \mbox{or}\, d$, which is independent of $v_{0,3}'$, thus we have
		\begin{equation}\nonumber
			\left|\frac{dy}{d v'_{0,3}}\right| = (s-s_1)\Big\{(1+|v'|^2)^{\f{|\kappa|}{2}}+|\kappa|(1+|v'|^2)^{\f{|\kappa|}{2}-1}(v_{0,3}')^2\Big\}\geq \frac{1}{N^6}, \quad \mbox{for} \  s_1\in [t'_{j+1},t'_j-\frac{1}{N^6}],
		\end{equation}
		which yields that
		\begin{align*}
			& \left( \iint_{\mathcal{D}}\nu^0(v'') | f^\lambda(\eta '_j-v'_{j,3} (t'_j-s_1),v'')|^2 dv'dv''\right)^{\frac12}\leq C_{m,N} \|\sqrt{\nu^0}f^\lambda\|_{L^2_{\eta ,v}}.
		\end{align*}

		Combining above estimates, then the RHS of \eqref{D2.56} is bounded by
		\begin{equation*}
			\frac{Cm^{2(\kappa-1)}}{N^2} \|h^\lambda\|_{L^\infty_{\eta ,v}}+C_{m,N} \|\sqrt{\nu^0}f^\lambda\|_{L^2_{\eta ,v}}.
		\end{equation*}

		Combining the above estimates, we obtain
		\begin{align*}
			\|h^\lambda\|_{L^\infty_{\eta ,v}}+|h^\lambda|_{L^\infty(\gamma_+)}
			&\leq C \big(m^{\kappa+3} +  m^{\kappa-1} e^{-\f{1}{2}\nu_0\tilde{k}_0d}+\frac{m^{2(\kappa-1)}}{N}\big) \, \big\{\|h^\lambda\|_{L^\infty_{\eta ,v}}+|h^\lambda|_{L^\infty(\gamma_+)}\big\}\nonumber\\
			&\, +C_{m,N} \|\sqrt{\nu^0}f^\lambda\|_{L^2_{\eta ,v}}+C_{m}\{ \|(\nu^0)^{-1}w_lg\|_{L^\infty_{\eta ,v}}+|w_{l+4}r|_{L^\infty(\gamma_+)}\}.
		\end{align*}
		First $m$ sufficiently small, then taking Taking $\tilde{k}_0$ and $N$ suitably large  so that
		\begin{equation}\nonumber
		C \big(m^{\kappa+3} +  m^{\kappa-1} e^{-\f{1}{2}\nu_0\tilde{k}_0d}+\frac{m^{2(\kappa-1)}}{N}\big) \leq \f12,
		\end{equation}
		then one has
		\begin{align*}
			\|h^\lambda\|_{L^\infty_{\eta ,v}}+|h^\lambda|_{L^\infty(\gamma_+)}\leq C \|\sqrt{\nu^0}f^\lambda\|_{L^2_{\eta ,v}}+C  \big\{ \|(\nu^0)^{-1}w_lg\|_{L^\infty_{\eta ,v}}+|w_{l+4}r|_{L^\infty(\gamma_+)}\big\}.
		\end{align*}
Therefore the proof of Lemma \ref{lem3.4} is finished.
	$\hfill\Box$

	\begin{lemma}\label{lem3.3}
		Let $\delta>0, d\geq 1$, $n\geq n_0$, and $l> 2$. Assume $\|(\nu^0)^{-1} w_l g\|_{L^\infty_{\eta,v}}<\infty$.
		Then there exists a unique solution $f^n$ to the following boundary value problem
		\begin{align}\label{D2.33}
			\begin{cases}
				\d f^n+v_3\partial_\eta f^n+\nu^0(v)f^n-K^0f^n=g,\\
				f^n(\eta,v)|_{\gamma_-}=(1-\f1n)f^n(\eta,R_\eta v),
			\end{cases}
		(\eta,v)\in\Omega_d\times\R^3,
		\end{align} satisfying
		\begin{align*}
			\| w_l f^{n}\|_{L^\infty_{\eta,v}}+| w_l f^n|_{L^\infty(\gamma_+)}\leq C_{\delta,d} \|(\nu^0)^{-1} w_l g\|_{L^\infty_{\eta,v}},
		\end{align*}
		where the positive constant $C_{\delta,d}>0$ depends only on $\delta,d$. Moreover, if
		g is continuous in $\Omega_d\times\mathbb{R}^2$, then $f^n$ is continuous away from grazing set $\gamma_0$.
	\end{lemma}
	\noindent{\bf Proof.}  We consider the solvability of the  following boundary value problem
	\begin{equation}\label{D2.34}
		\begin{cases}
			\mathcal{L}_\lambda f:=\delta f+v_3 \partial_\eta f+\nu^0(v) f-\lambda K^0f=g,\\
			f(\eta,v)|_{\gamma_-}=(1-\frac1n) f(\eta,R_\eta v),
		\end{cases}
	\end{equation}
	for $\lambda\in[0,1]$. For brevity we denote $\mathcal{L}_\lambda^{-1}$ to be the solution operator associated with the problem,  $f:=\mathcal{L}_\lambda^{-1} g$ is a solution to the BVP \eqref{D2.34}. Our idea is to prove the existence of $\mathcal{L}_0^{-1}$, and then extend to obtain the existence of  $\mathcal{L}_1^{-1}$ by a continuous argument on $\lambda$. We split the proof  into several steps.   
	
	\medskip
	\noindent{\it Step 1.} In this step, we prove the existence of $\mathcal{L}_0^{-1}$. We consider the following approximate sequence
	\begin{align}\label{D2.35}
		\begin{cases}
			\mathcal{L}_0 f^{i+1}=\delta f^{i+1}+v_3\partial_\eta f^{i+1} + \nu^0(v) f^{i+1}=g,\\
			f^{i+1}(\eta,v)|_{\gamma_-}=(1-\frac1n) f^{i}(\eta,R_\eta v),
		\end{cases}
	\end{align}
	for $i=0,1,2,\cdots$, where we have set $f^0\equiv0$. We will construct $L^\infty$  solutions to \eqref{D2.35} for $i=0,1,2,\cdots$, and 
	establish uniform $L^\infty$-estimates.
	\smallskip
	
	Firstly, multiplying \eqref{D2.35} by $f^{i+1}$ and integrating the resultant equality over $\Omega_d\times\mathbb{R}^3$,  one obtains that
	\begin{align}\label{D2.36}
		& \delta\|f^{i+1}\|^2_{L^2_{\eta,v}}+\frac12|f^{i+1}|^2_{L^2(\gamma_+)}+\|\sqrt{\nu^0}f^{i+1}\|_{L^2_{\eta,v}}\nonumber\\
		&\leq \frac12(1-\frac1n)^2|f^{i}|^2_{L^2(\gamma_+)}+C_\d\|g\|^2_{L^2_{\eta,v}}+ \frac{\d}{2}\|f^{i+1}\|^2_{L^2_{\eta,v}},
	\end{align}
	which yields that
	\begin{align*}
		\delta\|f^{i+1}\|^2_{L^2_{\eta,v}}+|f^{i+1}|^2_{L^2(\gamma_+)}
		\leq (1-\frac1n)^2|f^{i}|^2_{L^2(\gamma_+)}+C_\d \|g\|^2_{L^2_{\eta,v}}.
	\end{align*}
Considering the equation of $f^{i+1}-f^i$, by similar energy estimate as above, one obtains
	\begin{align}\label{D2.38}
		&\delta\|f^{i+1}-f^i\|^2_{L^2_{\eta,v}}+|f^{i+1}-f^i|^2_{L^2(\gamma_+)}\nonumber\\
		&\leq (1-\frac1n)^2|f^{i}-f^{i-1}|^2_{L^2(\gamma_+)} \leq \cdots \leq (1-\frac1n)^{2i} |f^1|^2_{L^2(\gamma_+)}\nonumber\\
		&\leq  C_\d(1-\frac1n)^{2i} \|g\|_{L^2_{\eta,v}}^2<\infty.
	\end{align}
	Noting $1-\frac1n<1$, thus $\{f^{i}\}_{i=0}^\infty$ is a Cauchy sequence in $L^2$, i.e.,
	\begin{align*}
		|f^{i}-f^j|^2_{L^2(\gamma_+)}+\|
		{f^i-f^j}\|^2_{L^2_{\eta,v}}\rightarrow0,\quad\mbox{as} \  i,j\rightarrow\infty,
	\end{align*}
	and we have, for $i=0,1,2,\cdots$,
	\begin{align}\label{D2.39}
		|f^{i}|^2_{L^2(\gamma_+)}+\|f^{i}\|^2_{L^2_{\eta,v}}\leq C_\d\|g\|^2_{L^2_{\eta,v}}.
	\end{align}

	\vspace{1mm}
	
	Next we consider the uniform  $L^\infty_{\eta,v}$ estimate. Let $h^i=w_lf^i$, one has that
	\begin{align*}
		h^{i+1}(\eta,v)e^{\hat{\nu}_\delta(v)t}=(1-\f1n)h^i(\eta_1,v_1)e^{\hat{\nu}_\d(v)t_1}+\int_{t_1}^t e^{\hat{\nu}_\d s}(1+|v|^2)^{\f{|\kappa|}{2}}w_lgds.
	\end{align*}
	Then it is easy to obtain
	\begin{align*}
		\|h^{1}\|_{L^\infty_{\eta,v}}+|h^1|_{L^\infty(\gamma_+)}\leq \|(\nu^0)^{-1}w_lg\|_{L^\infty_{\eta,v}}.
	\end{align*}
Also, by iteration, it holds that
\begin{align*}
	\|h^{i}\|_{L^\infty_{\eta,v}}+|h^i|_{L^\infty(\gamma_+)}\leq C_i\|(\nu^0)^{-1}w_lg\|_{L^\infty_{\eta,v}}, \quad i=0,1,2...,
\end{align*}
where the constants $C_i$ depend on $i$. Taking $h^{i+1}-h^i$, similarly, one obtains
	\begin{align*}
		(h^{i+1}-h^i)(\eta,v)&=(1-\f1n)
		e^{-\hat{\nu}_\d(t-t_1)}(h^i-h^{i-1})(\eta_1,v_1),
	\end{align*}
	which yields that
	\begin{align*}
		&\|h^{i+1}-h^i\|_{L^\infty_{\eta,v}}+|h^{i+1}-h^i|_{L^\infty(\gamma_+)}\leq (1-\f1n)|(h^i-h^{i-1})|_{L^\infty(\gamma_+)}\nonumber\\
		&\leq ...\leq (1-\f1n)^i|h^1|_{L^\infty(\gamma_+)}\leq (1-\f1n)^i\|(\nu^0)^{-1}w_lg\|_{L^\infty_{\eta,v}}.
	\end{align*}
	Since $1-\f1n<1$, $\{h^i\}_{i=0}^{\infty}$ is a Cauchy sequence in $L^\infty$ space. Then there exists a solution to $\mathcal{L}_0^{-1}$ with $L^2$ and $L^\infty$-weighted bound.
	
	\medskip
\noindent{\it Step 2.} Assume $f$ is a solution to \eqref{D2.34} and $\|w_l f\|_{L^\infty_{\eta,v}}+|w_lf|_{L^\infty(\gamma_+)}<\infty$. Multiplying $f$ to \eqref{D2.34} and integrating on $\R^3$, one obtains
\begin{align}\label{D2.43}
	\d\|f\|_{L^2_v}^2+\f{d}{d\eta}\int_{\R^3} v_3f^2\ dv+\lambda c_1\|(\FI-\FP_0)f\|_{\nu}^2\leq \f{\d}{4}\|f\|_{L^2_v}^2+\f{C}{\d}\|g\|_{L^2_v}^2,
\end{align}
where we used
\begin{align*}
	\langle f,\nu^0(v)f\rangle -\lambda\langle f,K^0f\rangle \geq \lambda c_1\|(\FI-\FP_0)f\|_{\nu}^2+C(1-\lambda)\|f\|_{\nu}^2.
\end{align*}
A direct calculation shows that
\begin{align*}
	\int_0^d\int_{\R^3}\f{d}{d\eta}(v_3f^2)\ dvd\eta &=\int_{\R^3} v_3|f|^2(d)\ dv -\int_{\R^3} v_3|f|^2(0) \ dv\nonumber\\
	&=\int_{v_3>0}v_3f^2(d,v)\ dv +\int_{v_3<0} (1-\f1n)^2 v_3|f|^2(d,Rv)\ dv\\
	&\quad -\int_{v_3>0}(1-\f1n)^2 v_3|f|^2(0,Rv)\ dv-\int_{v_3<0}v_3|f|^2(0,v)\ dv\\
	&=\big\{1-(1-\f1n)^2\big\}\Big\{\int_{v_3>0}v_3|f|^2(d,v)\ dv-\int_{v_3<0}v_3|f|^2(0,v)\ dv \Big\}\geq0.
\end{align*}
Hence, integrating \eqref{D2.43} on $[0,d]$, one has that
\begin{align}\label{D2.44}
	\|f\|_{L^2_{\eta,v}}\leq C_\d \|g\|_{L^2_{\eta,v}}.
\end{align}
From Lemma \ref{lem3.4}, we have the {\it a prior } estimate, i.e.
	\begin{align}\label{D2.60}
		\|w_l\mathcal{L}_\lambda^{-1} g\|_{L^\infty_{\eta,v}}+|w_l\mathcal{L}_\lambda^{-1}g|_{L^\infty(\gamma_+)}\leq C_{\d,d}\|(\nu^0)^{-1}w_lg\|_{L^\infty_{\eta,v}}.
	\end{align}
	On the other hand, let $(\nu^0)^{-1} w_l g_1 \in L^\infty_{\eta,v}$ and $(\nu^0)^{-1} w_l g_2 \in L^\infty_{\eta,v}$.  Let  $f^n_1=\mathcal{L}^{-1}_\lambda g_1$ and $f^n_2=\mathcal{L}^{-1}_\lambda g_2$ be the solutions to \eqref{D2.42} with $g$ replaced by $g_1$ and $g_2$, respectively. Then we have that
	\begin{align*}
		\begin{cases}
			\delta (f^{n}_2-f^n_1)+v_3\partial_\eta (f^{n}_2-f^n_1)
			+\nu^0(v) (f^{n}_2-f^n_1)-\lambda K^0(f^{n}_2-f^n_1)=g_2-g_1,\\[2mm]
			(f^{n}_2-f^n_1)(\eta,v)|_{\gamma_-}=(1-\frac1n) (f^{n}_2-f^n_1)(\eta,R_\eta v).
		\end{cases}
	\end{align*}
	
	Using \eqref{D2.44}, we obtain
	\begin{equation}\label{D2.61}
		\|\mathcal{L}^{-1}_\lambda g_2-\mathcal{L}^{-1}_\lambda g_1\|^2_{L^2_{\eta,v}}\leq  C_{\d} \|g_2-g_1\|^2_{L^2_{\eta,v}},
	\end{equation}
	and
	\begin{equation}\label{D2.62}
		\| w_l(\mathcal{L}^{-1}_\lambda g_2-\mathcal{L}^{-1}_\lambda g_1)\|_{L^\infty_{\eta,v}}
		+| w_l(\mathcal{L}^{-1}_\lambda g_2-\mathcal{L}^{-1}_\lambda g_1)|_{L^\infty(\gamma_+)}\leq C_{\delta,d}  \|(\nu^0)^{-1} w_l(g_2-g_1)\|_{L^\infty_{\eta,v}}.
	\end{equation}
	The uniqueness of solution to \eqref{D2.34} also follows from \eqref{D2.62}.
	We point out that the constant $C_{\delta,d}$ does not depend on $\lambda\in[0,1]$. This property is crucial for us to extend $\mathcal{L}_0^{-1}$ to $\mathcal{L}_1^{-1}$ by a bootstrap argument.

	\medskip

	\noindent{\it Step 3.}   Similar to \cite{Jiang-Wang}, by fixed point theory, we can obtain the existence of $\mathcal{L}_\lambda^{-1}$ for $\lambda\in[0,\lambda_\ast]$. Step by step, we can finally obtain the existence of operator $\mathcal{L}_1^{-1}$,  with \eqref{D2.44},
	\eqref{D2.60}, \eqref{D2.61} and \eqref{D2.62}. The continuity is easy to obtain by noting the convergence of $L^\infty_{\eta,v}$.  Therefore we complete the proof of Lemma \ref{lem3.3}.
$\hfill\Box$

	\smallskip
	
	\begin{lemma}\label{lem3.5}
		Let $\delta>0, d\geq 1$ and $l>2$, assume $\|(\nu^0)^{-1} w_l g\|_{L^\infty_{\eta,v}}<\infty$. Then there exists a unique solution $f^\delta$ to the approximate linearized steady Boltzmann equation \eqref{D2.25}. Moreover, it satisfies
		\begin{align}\label{D2.63}
			\|  f^\delta\|_{L^2_{\eta,v}} \leq C_{\delta,d} \|(\nu^0)^{-1} w_l g\|_{L^\infty_{\eta,v}},
		\end{align}
	and
	\begin{align}\label{D2.64}
		\|w_lf^{\d}\|_{L^\infty_{\eta,v}}+|w_lf^\d|_{L^\infty(\gamma_+)}\leq C_{\d,d}\|(\nu^0)^{-1}w_lg\|_{L^\infty_{\eta,v}}.
	\end{align}
		where the positive constant $C_{\delta,d}>0$ depends only on $\delta$ and $d$.
	\end{lemma}
	\noindent{\bf Proof.} Let $f^n$ be the solution of \eqref{D2.33} constructed in Lemma \ref{lem3.3} for $n\geq n_0$ with $n_0$ large enough.   Multiplying \eqref{D2.33} by $f^n$, one obtains that
	\begin{align*}
		\delta\|f^n\|^2_{L^2_{\eta,v}}+c_1 \|{(\FI-\FP_0)}f^n\|^2_{\nu}
		\leq C_{\delta} \|g\|^2_{L^2_{\eta,v}}.
	\end{align*}
	We apply Lemma \ref{lem3.4} to obtain
	\begin{align}\label{D2.65}
		\| w_l f^n\|_{L^\infty_{\eta,v}}{+| w_l f^n|_{L^\infty(\g_+)}}&\leq C \Big\{\|(\nu^0)^{-1} w_l g\|_{L^\infty_{\eta,v}}+\|f^n\|_{L^2_{\eta,v}}\Big\}\leq C_{\delta,d}\|(\nu^0)^{-1} w_l g\|_{L^\infty_{\eta,v}} .
	\end{align}
	Taking the difference $f^{n_1}-f^{n_2}$ with $n_1, n_2\geq n_0$, we know that
	\begin{align}\label{D2.66}
		\begin{cases}
			\delta (f^{n_1}-f^{n_2})+ v_3\partial_\eta (f^{n_1}-f^{n_2})+\FL_0 (f^{n_1}-f^{n_2})=0,\\[2mm]
			(f^{n_1}-f^{n_2})(\eta,v)|_{\gamma_{-}}=(1-\frac{1}{n_1})  (f^{n_1}-f^{n_2})(\eta,R_\eta v)+(\frac{1}{n_2}-\frac{1}{n_1}) f^{n_2}(\eta,R_\eta v).
		\end{cases}
	\end{align}
	Multiplying \eqref{D2.66} by $ f^{n_1}-f^{n_2}$ and integrating it over $\Omega_d\times\mathbb{R}^3$, we can obtain
	\begin{align}\label{D2.67}
		&\delta\|(f^{n_1}-f^{n_2})\|^2_{L^2_{\eta,v}} +c_1\int_0^d \|{(\FI-\FP_0)}(f^{n_1}-f^{n_2})\|^2_{\nu}\ d\eta \nonumber\\
		&\leq C (\frac1{n_1}+\frac1{n_2})\int_{\gamma_-}  |v_3(|f^{n_1}|+|f^{n_2}|) f^{n_2}|dv+C(\f{1}{n_1}-\f{1}{n_2})^2|w_lf^{n_2}|_{L^\infty(\gamma_+)}^2\nonumber\\
		&\leq C (\frac1{n_1}+\frac1{n_2}) (|w_lf^{n_1}|^2_{L^\infty(\gamma_+)}+|w_lf^{n_2}|^2_{L^\infty(\gamma_+)})\nonumber\\
		&\leq C_{\delta,d}(\frac1{n_1}+\frac1{n_2}) \|(\nu^0)^{-1}w_lg\|^2_{L^\infty_{\eta,v}}\rightarrow0,
	\end{align}
	as $n_1,n_2 \rightarrow \infty$,
	where we have used the uniform estimate \eqref{D2.65} in the last inequality.
	Also, by Lemma \ref{lem3.4},
	\begin{align}
		&\|w_{l-4}(f^{n_1}-f^{n_2})\|_{L^\infty_{\eta,v}}+|w_{l-4}(f^{n_1}-f^{n_2})|_{L^\infty(\gamma_+)}\leq  C_{\d,d}\sqrt{\f{1}{n_1}+\f{1}{n_2}}\|(\nu^0)^{-1}w_lg\|_{L^\infty_{\eta,v}}\rightarrow0.
	\end{align}
We denote $f^\delta$ to be the Cauchy limit of $\{f^n\}$ in $L^2$ space and $L^\infty$ space, noting \eqref{D2.65}, then it is direct to check that $f^\delta$ is a weak solution to \eqref{D2.25} and \eqref{D2.63}-\eqref{D2.64} hold. Therefore we have completed the proof of Lemma \ref{lem3.5}.
$\hfill\Box$	
	
	\vspace{1mm}
	
	\
	
 Multiplying \eqref{D2.25} by $(1,v_1-\fu_1^0,v_2-\fu_2^0,|v-\fu^0|^2-3T^0)\sqrt{\mu_0}$ and using \eqref{D2.2-0}, one obtains that
	\begin{align}\label{D2.71}
		\int_0^d\int_{\R^3} (1,v_1-\fu^0_1,v_2-\fu^0_2,|v-\fu^0|^2-3T^0) \sqrt{\mu_0} f^\d(\eta ,v) dvd\eta =0.
	\end{align}
We denote
	\begin{equation*}
		\mathbf{P}_0f^\d(\eta ,v)=\{a^\d(\eta )+b^\d\cdot (v-\fu^0)+c^\d(\eta ) (|v-\fu^0|^2-3T^0)\}\sqrt{\mu_0},
	\end{equation*}
	then it follows from \eqref{D2.71} that
	\begin{equation}\label{D2.72}
		\int_0^d a^\d(\eta ) d\eta =\int_0^d b^\d_1(\eta ) d\eta =\int_0^d b^\d_2(\eta ) d\eta =\int_0^d c^\d(\eta ) d\eta =0.
	\end{equation}

	\begin{lemma}\label{lem3.8}
		Let $d\geq 1$, $l>2$.  Assume \eqref{D2.3} and  $\|(\nu^0)^{-1}w_lg\|_{L^\infty_{\eta ,v}}<\infty$.   Then there exists a  unique solution $f=f(\eta ,v)$  to the linearized steady Boltzmann equation
		\begin{align}\label{D2.74}
			\begin{cases}
				v_3 \partial_\eta f+\FL_0 f=g,\quad (\eta ,v)\in \Omega_d\times \R^3, \\
				f(\eta ,v)|_{\gamma_{-}}=f(\eta ,R_\eta v),
			\end{cases}
		\end{align}
		with
		\begin{equation}\label{D2.75}
			\|w_lf\|_{L^\infty_{\eta ,v}} +|w_lf|_{L^\infty(\gamma_+)} \leq C_d \|(\nu^0)^{-1}w_lg\|_{L^\infty_{\eta ,v}}.
		\end{equation}
	 Moreover, if $g$ is continuous in
	$\Omega_d\times\R^3$, then $f$ is continuous away from grazing set $\gamma_0$.
	\end{lemma}

\noindent{\bf Proof.} Let $f^{\delta}$ be the solution of \eqref{D2.25}  constructed in Lemma \ref{lem3.5}. We shall consider the limit $\delta\to0$ to obtain solution of \eqref{D2.74}. 

By similar arguments as in \cite[Lemma 3.7]{Jiang-Wang}, we can obtain
\begin{align}\label{D.53}
\|\FP_0 f^{\delta}\|^2 \leq C d^{6} \Big( \|(\mathbf{I}-\FP_0) f^{\delta}\|_{\nu}^2 + \|g\|^2_{L_{\eta,v}} \Big).
\end{align}
It is noted that \cite[Lemma 3.7]{Jiang-Wang} was proved for hard sphere case, but the proof can be generalized to both hard and soft potentials without any difficulty.

Multiplying \eqref{D2.25} by $f^\delta$ and integrating   over $\Omega_d\times\mathbb{R}^3$, we have
\begin{align}\label{S3.109}
\delta\|f^\delta\|^2_{L^2_{\eta,v}}+c_0\|(\mathbf{I-P})f^\delta\|^2_{\nu}\leq\vartheta \|f^\delta\|^2_{L^2_{\eta,v}}+C_{\vartheta} \|g\|^2_{L^2_{\eta,v}}.
\end{align}
which, together with  \eqref{D.53} and taking $\vartheta$ small enough (depending on $d$), yields that
\begin{equation}\label{3.94}
	\|\sqrt{\nu^0}f^\delta\|^2_{L^{2}_{\eta,v}} \leq C_{d} \|g\|^2_{L^2_{\eta,v}}.
\end{equation}

Applying \eqref{D2.59} to $f^\delta$ and using \eqref{3.94}, one obtain
\begin{equation}\label{S3.119}
	\|w_lf^\delta\|_{L^\infty_{\eta,v}}{+|w_l f^\delta|_{L^\infty(\g_+)}}\leq C_d\|(\nu^0)^{-1}w_l g\|_{L^\infty_{\eta,v}}.
\end{equation}
Next we consider the convergence of $f^\delta$ as $\delta\rightarrow0+$. For any $\delta_1,\delta_2>0$, we consider the difference $f^{\delta_2}-f^{\delta_1}$ satisfying
\begin{equation}\label{S3.120}
	\begin{cases}
		v_3\partial_\eta (f^{\delta_2}-f^{\delta_1})+ \FL_0 (f^{\delta_2}-f^{\delta_1})=-\delta_2 f^{\delta_2}+\delta_1 f^{\delta_1},\\[2mm]
		(f^{\delta_2}-f^{\delta_1})|_{\gamma_-}=(f^{\delta_2}-f^{\delta_1})(\eta,R_\eta v).
	\end{cases}
\end{equation}
Multiplying \eqref{S3.120} by $f^{\delta_2}-f^{\delta_1}$, integrating the resultant equation and by similar arguments as in \eqref{S3.109}-\eqref{S3.119}, one gets
\begin{align}\label{S3.121}
	\|\sqrt{\nu^0}(f^{\delta_2}-f^{\delta_1})\|^2_{L_{\eta,v}}
	&\leq C_d \|\delta_2 f^{\delta_2}-\delta_1 f^{\delta_1}\|^2_{L^2_{\eta,v}}
	 \leq C_d (\delta_1^2+\delta_2^2)\cdot  \|(\nu^0)^{-1}w_l g\|_{L^\infty_{\eta,v}}^2\rightarrow0,
\end{align}
as $\delta_1$, $\delta_2\rightarrow 0+$. Finally, applying \eqref{D2.59} to $f^{\delta_2}-f^{\delta_1}$ and using \eqref{S3.121}, then we obtain
\begin{align}\label{S3.122}
	&\|w_l (f^{\delta_2}-f^{\delta_1})\|_{L^\infty_{\eta,v}} +| w_l (f^{\delta_2}-f^{\delta_1})|_{L^\infty(\g_+)}\nonumber\\
	&\leq C\Big\{ \|(\nu^0)^{-1}w_l (\delta_2 f^{\delta_2}-\delta_1 f^{\delta_1})\|_{L^\infty_{\eta,v}} +\|\sqrt{\nu^0}(f^{\delta_2}-f^{\delta_1})\|_{L^2_{\eta,v}} \Big\}\nonumber\\
	&\leq C_d(\delta_1+\delta_2)\|(\nu^0)^{-1}w_l g\|_{L^\infty_{\eta,v}}\rightarrow 0,
\end{align}
as $\delta_1$, $\delta_2\rightarrow 0+$,
With \eqref{S3.122}, we know that there exists a function $f$ so that $\|w_l\,(f^{\delta}-f)\|_{L^\infty_{\eta,v}}\rightarrow0$ as $\delta\rightarrow 0+$. And it is  direct to see that $f$ solves \eqref{D2.74}. Also,  \eqref{D2.75} follows immediately from \eqref{S3.119}. The continuity of $f$
directly follows from the $L^\infty_{\eta,v}$-convergence and the continuity of $f^\delta$. Therefore the proof of Lemma \ref{lem3.8} is complete.

$\hfill\Box$

	To obtain the solution for half-space problem, we need some uniform estimates independent of $d$, then we can take the limit $d\rightarrow\infty$. Let $f$ be the solution of \eqref{D2.74}, we denote
	\begin{align*}
		\FP_0 f(\eta ,v)=\big[a(\eta )+b(\eta )\cdot (v-\fu^0)+c(\eta ) (|v-\fu^0|^2-3T^0)\big] \sqrt{\mu_0}.
	\end{align*}
	Multiplying \eqref{D2.74} by $\sqrt{\mu_0}$ and using \eqref{D2.2-0}, we have
	\begin{align}\label{D2.76}
		0=\frac{d}{d\eta } \int_{\R^3} v_3 \sqrt{\mu_0} f(\eta ,v)dv = \frac{d}{d\eta } b_3(\eta )\equiv0.
	\end{align}
	Since $f$ satisfies the specular boundary, it holds that $b_3(\eta )|_{\eta =0}=b_3(\eta )|_{\eta =d}=0$, which, together with \eqref{D2.76}, yields
	\begin{equation}\label{D2.77}
		b_3(\eta )=0,\quad \mbox{for}\  \eta \in[0,d].
	\end{equation}
	
	Let $(\phi_0,\phi_1,\phi_2,\phi_3)$ be some constants chosen later, we define
	\begin{align*}
		\bar{f}(\eta ,v)&:=f(\eta ,v)+[\phi_0+\phi_1 (v_1-\fu^0_1)+\phi_2 (v_2-\fu^0_2)+\phi_3 (|v-\fu^0|^2-3T^0)]\sqrt{\mu_0}\nonumber\\
		&=[\bar{a}(\eta )+\bar{b}_1(\eta )(v_1-\fu^0_1)+\bar{b}_2(\eta )(v_2-\fu^0_2)+\bar{c}(\eta ) (|v-\fu^0|^2-3T^0)]\sqrt{\mu_0}\nonumber\\
		&\qquad +(\mathbf{I}-\FP_0) \bar{f},
	\end{align*}
	where
	\begin{align*}
		\begin{cases}
			\bar{a}(\eta )=a(\eta )+\phi_0,\\
			\bar{b}_i(\eta )=b_i(\eta )+\phi_i, \quad i=1,2,\\
			\bar{c}(\eta )=c(\eta )+\phi_3.
		\end{cases}
	\end{align*}
It follows from \eqref{D2.77} that
	\begin{equation}\label{D2.78}
		\bar{b}_3(\eta )\equiv0\quad \mbox{and}\quad(\mathbf{I}-\FP_0) \bar{f}(\eta ,v)\equiv(\mathbf{I}-\FP_0) f(\eta ,v)\quad \forall \eta \in[0,d].
	\end{equation}
The equation for $\bar{f}$ is
	\begin{align}\label{D2.79}
		\begin{cases}
			v_3 \partial_\eta \bar{f}+\FL_0 \bar{f}=g,\quad  (\eta ,v)\in \Omega_{d}\times\R^3,\\
			\bar{f}(\eta ,v)|_{\gamma_{-}}=\bar{f}(\eta ,R_\eta v).
		\end{cases}
	\end{align}
	Hence it follows from \eqref{D2.75} that
	\begin{equation*}
		\|w_l\,\bar{f}\|_{L^\infty_{\eta ,v}} +|w_{l}\,\bar{f}|_{L^\infty(\gamma_+)} \leq C_d \|(\nu^0)^{-1}w_{l}g\|_{L^\infty_{\eta ,v}}+C_d |(\phi_0,\phi_1,\phi_2,\phi_3)|.
	\end{equation*}
	Multiplying $\eqref{D2.79}_1$ by $(v_1-\fu^0_1,v_2-\fu^0_2,|v-\fu^0|^2-5T^0) \sqrt{\mu_0}$ and using \eqref{D2.2-0}, we get
	\begin{align}\label{D2.80}
		\begin{split}
			\int_{\R^3} v_3 (v_i-\fu^0_i) \sqrt{\mu_0} \bar{f}(\eta ,v)dv&=0, \quad \forall\, \eta \in[0,d],\quad i=1,2,\\
			\int_{\R^3} v_3 (|v-\fu^0|^2-5T^0) \sqrt{\mu_0} \bar{f}(\eta ,v)dv&=0,\quad \forall \, \eta \in[0,d].
		\end{split}
	\end{align}

It follows from \eqref{D2.78} and \eqref{D2.80} that
	\begin{align}\label{D2.82}
	\int_{\R^3} v_3|\FP_0\bar{f}(\eta ,v)|^2 dv\equiv\int_{\R^3} v_3\FP_0\bar{f}(\eta ,v)\cdot (\mathbf{I}-\FP_0)\bar{f}(\eta ,v) dv\equiv0,
	\end{align}
	which yields that
	\begin{align}\label{D2.83}
		\int_{\R^3} v_3|\bar{f}(\eta ,v)|^2 dv=\int_{\R^3} v_3|(\mathbf{I}-\FP_0)\bar{f}(\eta ,v)|^2 dv,\quad \forall \eta \in[0,d].
	\end{align}
	Multiplying \eqref{D2.79} by $\bar{f}$ and using \eqref{D2.83}, \eqref{D2.2-0}, we obtain
	\begin{align*}
		\frac{d}{d\eta } \int_{\R^3} v_3 |(\mathbf{I}-\FP_0)\bar{f}|^2dv+\f12c_1\|(\mathbf{I}-\FP_0)\bar{f}\|^2_{\nu}\leq C \|(\nu^0)^{-\f12}g\|_{L^2_v}^2,
	\end{align*}
	which yields that
	\begin{align}\label{D2.84}
		\int_0^d\|(\FI-\FP_0)\bar{f}\|_{\nu}^2\ d\eta\leq C\int_0^d\|(\nu^0)^{-\f12}g\|_{L^2_{v}}^2d\eta ,
	\end{align}
	where we have used \eqref{D2.2-0} to derive
	\begin{align}\nonumber
		\int_{\R^3} g \bar{f}dv=\int_{\R^3} g (\mathbf{I}-\FP_0)\bar{f} dv\leq \frac12c_1\|(\mathbf{I}-\FP_0)\bar{f}\|_{\nu}^2+C\|(\nu^0)^{-\f12}g\|_{L^2_v}^2.
	\end{align}
	\begin{lemma}\label{lem3.9}
		There exist constants $(\phi_0,\phi_1,\phi_2,\phi_3)$  such that
		\begin{align}\label{D2.85}
			\begin{split}
				&\int_{\R^3} v_3\bar{f}(d,v)\cdot v_3 \sqrt{\mu_0} dv=0,\\
				&\int_{\R^3} v_3\bar{f}(d,v)\cdot \FL_0^{-1} (\mathcal{A}_{3i}^0) dv=0,\  i=1,2,\\
				&\int_{\R^3} v_3\bar{f}(d,v)\cdot \FL_0^{-1} (\mathcal{B}_{3}^0) dv=0.
			\end{split}
		\end{align}
	\end{lemma}
	
\noindent{\bf Proof.} A direct calculation shows that
	\begin{align}
		\int_{\R^3} v_3\bar{f}(\eta ,v)\cdot v_3 \sqrt{\mu_0} dv&=\rho^0T^0\bar{a}(\eta ) +2\rho^0(T^0)^2\bar{c}(\eta ) +T^0\int_{\R^3} \mathcal{A}_{33}^0(v)\cdot (\mathbf{I}-\FP_0)\bar{f}(\eta ,v) dv\nonumber\\
		&=\rho^0T^0\phi_0+2\rho^0(T^0)^2\phi_3 +\rho^0T^0a(\eta )+2\rho^0(T^0)^2c(\eta )\nonumber\\
		&\quad+T^0\int_{\R^3} \mathcal{A}_{33}^0(v)\cdot (\mathbf{I}-\FP_0)f(\eta ,v) dv,\label{D2.86}\\
		\int_{\R^3} v_3\bar{f}(\eta ,v)\cdot \FL_0^{-1} (\mathcal{A}_{31}^0) dv
		&= \mu(T^0) \bar{b}_1(\eta )+\int_{\R^3} v_3 (\mathbf{I}-\FP_0)\bar{f}(\eta ,v)\cdot \FL_0^{-1} (\mathcal{A}_{31}^0)  dv\nonumber\\
		&=\mu(T^0) \phi_1+\mu(T^0) b_1(\eta ) +\int_{\R^3} v_3 (\mathbf{I}-\FP_0)f(\eta ,v)\cdot \FL_0^{-1} (\mathcal{A}_{31}^0)  dv,\label{D2.87} \\
		\int_{\R^3} v_3\bar{f}(\eta ,v)\cdot \FL_0^{-1} (\mathcal{A}_{32}^0) dv
		&=\mu(T^0) \bar{b}_2(\eta )+\int_{\R^3} v_3 (\mathbf{I}-\FP_0)\bar{f}(\eta ,v)\cdot \FL_0^{-1} (\mathcal{A}_{32}^0)  dv\nonumber\\
		&=\mu(T^0) \phi_2+\mu(T^0) b_2(\eta ) +\int_{\R^3} v_3 (\mathbf{I}-\FP_0)f(\eta ,v)\cdot \FL_0^{-1} (\mathcal{A}_{32}^0)  dv,\label{D2.88}\\
		\int_{\R^3} v_3\bar{f}(\eta ,v)\cdot \FL_0^{-1} (\mathcal{B}_3^0) dv
		&=\kappa(T^0) \bar{c}(\eta )+\int_{\R^3} v_3 (\mathbf{I}-\FP_0)\bar{f}(\eta ,v)\cdot \FL_0^{-1} (\mathcal{B}_3^0)  dv\nonumber\\
		&=\kappa(T^0) \phi_3+\kappa(T^0) c(\eta )+\int_{\R^3} v_3 (\mathbf{I}-\FP_0)f(\eta ,v)\cdot \FL_0^{-1} (\mathcal{B}_3^0)  dv,\label{D2.89}
	\end{align}
	where we have used the notations in \eqref{D2.24}.
	
	Using  \eqref{D2.86}-\eqref{D2.89}, then \eqref{D2.85} is equivalent  as
	\begin{equation}\nonumber
		\left(
		\begin{array}{cccc}
			1 & 0 & 0 & 2T^0  \\
			0 & \mu(T^0) & 0 & 0  \\
			0 & 0 & \mu(T^0) & 0  \\
			0 & 0 & 0 & \kappa(T^0)  \\
		\end{array}
		\right)
		\left(
		\begin{array}{c}
			\phi_0 \\
			\phi_1\\
			\phi_2\\
			\phi_3\\
		\end{array}
		\right)
		=-\left(
		\begin{array}{c}
			a(d)+2T^0c(d)+\f{1}{\rho^0}\int_{\R^3}(\mathbf{I}-\FP_0)f(d,v) \cdot \mathcal{A}_{33}^0(v)dv \\[1mm]
			\mu(T^0) b_1(d) +\int_{\R^3} v_3 (\mathbf{I}-\FP_0)f(d,v)\cdot \FL_0^{-1} (\mathcal{A}_{31}^0)  dv\\[1mm]
			\mu(T^0) b_2(d) +\int_{\R^3} v_3 (\mathbf{I}-\FP_0)f(d,v)\cdot \FL_0^{-1} (\mathcal{A}_{32}^0)  dv\\[1mm]
			\kappa(T^0) c(d)+\int_{\R^3} v_3 (\mathbf{I}-\FP_0)f(d,v)\cdot \FL_0^{-1} (\mathcal{B}_3^0)  dv
		\end{array}
		\right).
	\end{equation}
	Noting the matrix is non-singular, hence  $(\phi_0,\phi_1,\phi_2,\phi_3)$ are found. Therefore the proof of Lemma \ref{lem3.9} is completed.
	$\hfill\Box$	
	
\smallskip

	From now on, the proof is quite different with hard sphere case since we do not have $\nu^0\geq \s|v_3|$ for soft cases. Hence it is hard to obtain the space exponential decay as hard sphere case. Our strategy is to get the space decay by losing the particle velocity weight.
	\begin{lemma}\label{lem3.10}
		Let $(\phi_0,\phi_1,\phi_2,\phi_3)$ be the ones determined  in Lemma \ref{lem3.9}, then it holds that
		\begin{equation}\label{D2.90}
			\int_0^d\|\bar{f}\|_{\nu}^2d\eta \leq C\int_0^d\int_{\R^3}(1+\eta )^{2p_0}(\nu^0)^{-1}g^2dvd\eta ,\quad p_0>1.
		\end{equation}
		where the constant $C>0$ is independent of $d$.
	\end{lemma}
	
	\noindent{\bf Proof.}
		Multiplying \eqref{D2.79} by  $\FL_0^{-1} (\mathcal{A}_{31}^0), \FL_0^{-1} (\mathcal{A}_{32}^0) $ and $\FL_0^{-1} (\mathcal{B}_3^0) $, respectively, one has from \eqref{D2.80} that
	\begin{align}
		\frac{d}{d\eta}\left(
		\begin{array}{c}
			\int_{\R^3} v_3\bar{f}(\eta ,v)\cdot \FL_0^{-1} (\mathcal{A}_{31}^0) dv\\[1mm]
			\int_{\R^3} v_3\bar{f}(\eta ,v)\cdot \FL_0^{-1} (\mathcal{A}_{32}^0) dv\\[1mm]
			\int_{\R^3} v_3\bar{f}(\eta ,v)\cdot \FL_0^{-1} (\mathcal{B}_3^0) dv\\
		\end{array}
		\right)
		=
		\left(
		\begin{array}{c}
			\int_{\R^3}  g\cdot \FL_0^{-1} (\mathcal{A}_{31}^0) dv\\[1mm]
			\int_{\R^3} g\cdot \FL_0^{-1} (\mathcal{A}_{32}^0) dv\\[1mm]
			\int_{\R^3} g\cdot \FL_0^{-1} (\mathcal{B}_3^0) dv\\
		\end{array}
		\right).\nonumber
	\end{align}
	Integrating the above system over $[\eta ,d]$ and using $\eqref{D2.85}$, one obtains
	\begin{align}
		\left(
		\begin{array}{c}
			\int_{\R^3} v_3\bar{f}(\eta ,v)\cdot \FL_0^{-1} (\mathcal{A}_{31}^0) dv\\[1mm]
			\int_{\R^3} v_3\bar{f}(\eta ,v)\cdot \FL_0^{-1} (\mathcal{A}_{32}^0) dv\\[1mm]
			\int_{\R^3} v_3\bar{f}(\eta ,v)\cdot \FL_0^{-1} (\mathcal{B}_3^0) dv\\
		\end{array}
		\right)
		=
		-\int_\eta ^d\left(
		\begin{array}{c}
			\int_{\R^3} g\cdot \FL_0^{-1} (\mathcal{A}_{31}^0) dv\\[1mm]
			\int_{\R^3} g\cdot \FL_0^{-1} (\mathcal{A}_{32}^0) dv\\[1mm]
			\int_{\R^3}g\cdot \FL_0^{-1} (\mathcal{B}_3^0) dv\\
		\end{array}
		\right)(z)dz,\nonumber
	\end{align}
	which, together with
	\eqref{D2.87}-\eqref{D2.89} and Proposition \ref{thmA.1}, yields that
	\begin{align}\label{D2.91}
		|(\mu(T^0) \bar{b}_{1},\mu(T^0)  \bar{b}_{2},\kappa(T^0) \bar{c})(\eta )|\leq C\|(\mathbf{I}-\FP_0)f(\eta )\|_{\nu}+C\int_\eta ^d\|g(z)\|_{L^2_v}dz,
	\end{align}
where we used Proposition \ref{thmA.1} to derive the decay estimates for $v_3\FL_0^{-1}(\mathcal{A}_{31}^0),v_3\FL_0^{-1}(\mathcal{A}_{32}^0),v_3\FL_0^{-1}(\mathcal{B}_{3}^0)$. \\

It follows from \eqref{D2.86} that
	\begin{align*}
		\bar{a}(\eta)=-2T^0\bar{c}(\eta)-\f{1}{\rho^0T^0}\int_{\R^3}(\FI-\FP_0)\bar{f}(\eta,v)\cdot v_3^2\sqrt{\mu_0}dv-\f{1}{\rho^0T^0}\int_\eta^d\int_{\R^3}g\cdot v_3\sqrt{\mu_0}dvdz,
	\end{align*}
	which yields that
	\begin{align}\label{D2.92}
		|\bar{a}(\eta)|\leq C\|(\mathbf{I}-\FP_0)f(\eta )\|_{\nu}+C\int_\eta ^d\|g(z)\|_{L^2_v}dz.
	\end{align}
	Using \eqref{D2.84}, \eqref{D2.91}-\eqref{D2.92}, one gets
	\begin{align}\label{D2.93}
		\int_0^d\|\FP_0\bar{f}\|_{\nu}^2d\eta \leq &C\int_0^d\|(\FI-\FP_0)f(\eta )\|_{\nu}^2d\eta +C\int_0^d\Big\{\int_\eta ^d\|g(z)\|_{L^2_v}dz\Big\}^2d\eta \nonumber\\
		\leq &C\int_0^d\int_{\R^3}(\nu^0)^{-1}g^2dvd\eta +C\int_0^d\int_\eta ^d(1+z)^{-2p_0}dzd\eta \cdot\int_0^d\int_{\R^3}(1+\eta )^{2p_0}g^2dvd\eta \nonumber\\
		\leq & C\int_0^d\int_{\R^3}(1+\eta )^{2p_0}(\nu^0)^{-1}g^2dvd\eta , \quad p_0>1.
	\end{align}
We conclude \eqref{D2.90} from \eqref{D2.93} and \eqref{D2.84}. The proof of Lemma \ref{lem3.10} is completed.
$\hfill\Box$	
	
\smallskip

 Since we will encounter some space weight $\eta^l$ in the formulation of Hilbert expansion (see \eqref{c1.22} for details), then we have to derive at least polynomial decay for Knudsen boundary layer so that the analysis can be closed.
	\begin{lemma}\label{lem3.11}
		Let $(\phi_0,\phi_1,\phi_2,\phi_3)$ be the ones determined  in Lemma \ref{lem3.9}, then it holds that
		\begin{equation}\label{D2.94}
			\int_0^d(1+\eta )^k\|w_l\bar{f}\|_{\nu}^2d\eta \leq C_k\int_0^d(1+\eta )^{2p_k}\|w_{l+2k+2}g\|_{L^2_v}^2d\eta ,\quad p_k>\f{k}{2}+1,
		\end{equation}
		where  $k$ is a non-negative integer and the constant $C_k$ depends only on $k$.
	\end{lemma}
	\noindent{\bf Proof.}
		We divide the proof into three steps.\\
		\noindent{\it Step 1.} Let $l$ be any positive constant. From \cite[Corollary 1]{Guo2}, it holds that
		\begin{align}\label{D2.96}
			\langle w_l^2\FL_0 \mathfrak{h},\mathfrak{h}\rangle \geq \f12\|w_l\mathfrak{h}\|_{\nu}^2-C\|\mathfrak{h}\|_{\nu}^2.
		\end{align}
		Multiplying \eqref{D2.79} by $w_l^2\bar{f}$ and integrating on $\R^3$, one has
		\begin{align}\label{D2.97}
			\f12\f{d}{d\eta} \int_{\R^3} v_3w_l^2\bar{f}^2dv+\int_{\R^3}w_l^2\bar{f}\cdot\FL_0 \bar{f}dv=\int_{\R^3}w_l^2\bar{f}gdv.
		\end{align}

		Integrating \eqref{D2.97} on $[0,d]$ and using \eqref{D2.96}, one gets
		\begin{align}\label{D2.98}
			\int_0^d \|w_l\bar{f}\|_{\nu}^2\ d\eta &\lesssim \int_0^d \|\bar{f}\|_{\nu}^2d\eta +\int_0^d\int_{\R^3}w_l^2\bar{f}g\ dvd\eta \nonumber\\
	 &\lesssim \int_0^d \|\bar{f}\|_{\nu}^2d\eta +\int_0^d\|(\nu^0)^{-\f12}w_lg\|_{L^2_v}^2d\eta \nonumber\\
			&\lesssim\int_0^d \|\bar{f}\|_{\nu}^2d\eta +\int_0^d\|w_{l+2}g\|_{L^2_v}^2d\eta \nonumber\\
			&\lesssim\int_0^d(1+\eta )^{2p_0}\|w_{l+2}g\|_{L^2_v}^2d\eta ,\quad p_0>1,
		\end{align}
		where we have used Lemma \ref{lem3.10}.
		
	\smallskip
	
		\noindent{\it Step 2.}
		Multiplying \eqref{D2.79} by $\bar{f}$ and integrating over $\R^3$, one has
		\begin{align*}
			\f12\f{d}{d\eta} \int_{\R^3}v_3\bar{f}^2dv+\int_{\R^3}\bar{f}\FL_0 \bar{f}dv =\int_{\R^3}\bar{f}gdv,
		\end{align*}
		which implies that
		\begin{align}\label{D2.100}
			\f{d}{d\eta} \int_{\R^3}v_3\bar{f}^2dv+c_1\|(\FI-\FP_0)\bar{f}\|_{\nu}^2\lesssim\|(\nu^0)^{-\f12}g\|_{L^2_v}^2.
		\end{align}
		Multiplying \eqref{D2.100} by $(1+\eta )^k$ with $k$ being some  positive integer, we get
		\begin{align}\label{D2.101}
			&\partial_\eta \big\{(1+\eta )^k\int_{\R^3} v_3\bar{f}^2dv\big\}+c_1(1+\eta )^k\|(\FI-\FP_0)\bar{f}\|_{\nu}^2\nonumber\\
			&\lesssim k(1+\eta )^{k-1}\int_{\R^3}v_3\bar{f}^2dv+(1+\eta )^k\|(\nu^0)^{-\f12}g\|_{L^2_v}^2.
		\end{align}
		Then, integrating \eqref{D2.101} on $[0,d]$, one obtains
		\begin{align}\label{D2.102}
			\int_0^d (1+\eta )^k\|(\FI-\FP_0)\bar{f}\|_{\nu}^2d\eta  &\lesssim k\int_0^d(1+\eta )^{k-1}\int_{\R^3}v_3\bar{f}^2dvd\eta +\int_0^d(1+\eta )^k\|(\nu^0)^{-\f12}g\|_{L^2_v}^2d\eta \nonumber\\
			&\lesssim k\int_0^d(1+\eta )^{k-1}\|w_2\bar{f}\|_{\nu}^2d\eta +\int_0^d(1+\eta )^k\|w_2g\|_{L^2_v}^2d\eta.
		\end{align}

		On the other hand, from \eqref{D2.91}-\eqref{D2.92}, one has that
		\begin{align}\label{D2.103}
			\int_0^d(1+\eta )^k\|\FP_0 \bar{f}\|_{L^2_v}^2d\eta &\lesssim_k\int_0^d(1+\eta )^k\|(\FI-\FP_0)\bar{f}\|_{\nu}^2d\eta +\int_0^d(1+\eta)^k\Big\{\int_\eta^d\|g(z)\|_{L^2_v}dz\Big\}^2d\eta\nonumber\\
			&\lesssim_k \int_0^d(1+\eta )^{k-1}\|w_2\bar{f}\|_{\nu}^2d\eta +\int_0^d(1+\eta )^{k}\|w_2g\|_{L^2_v}^2d\eta\nonumber\\
			&\quad+\int_0^d(1+\eta )^{2p_k}\|g\|_{L^2_v}^2d\eta \nonumber\\
			&\lesssim_k\int_0^d(1+\eta )^{k-1}\|w_2\bar{f}\|_{\nu}^2d\eta  +\int_0^d(1+\eta )^{2p_k}\|w_2g\|_{L^2_v}^2d\eta ,
		\end{align}
	where $ p_k>\f{k}{2}+1$.
		It follows from \eqref{D2.102}-\eqref{D2.103} that
		\begin{align}\label{D2.104}
			\int_0^d(1+\eta )^k\|\bar{f}\|_{\nu}^2d\eta \lesssim k\int_0^d(1+\eta )^{k-1}\|w_2\bar{f}\|_{\nu}^2d\eta  +\int_0^d(1+\eta )^{2p_k}\|w_2g\|_{L^2_v}^2d\eta ,\quad p_k>\f{k}{2}+1.
		\end{align}
		\noindent{\it Step 3.} Multiplying \eqref{D2.97} by $(1+\eta )^k$, one has
		\begin{align}\label{D2.105}
			&\f12\f{d}{d\eta} \Big\{(1+\eta )^k\int_{\R^3}v_3w_l^2\bar{f}^2dv\Big\}-\f{k}{2}(1+\eta )^{k-1}\int_{\R^3}v_3w_l^2\bar{f}^2dv\nonumber\\
			&+(1+\eta )^k\int_{\R^3}w_l^2\bar{f}\FL_0 \bar{f}dv=(1+\eta )^k\int_{\R^3}w_l^2\bar{f}gdv.
		\end{align}
		We have from \eqref{D2.96} that
		\begin{align*}
			(1+\eta )^k\int_{\R^3}w_l^2\bar{f}\FL_0 \bar{f}dv\geq \f12(1+\eta )^k\|w_l\bar{f}\|_{\nu}^2-C(1+\eta )^k\|\bar{f}\|_{\nu}^2,
		\end{align*}
		which, together with  \eqref{D2.104}-\eqref{D2.105}, yields that
		\begin{align}\label{D2.106}
			\int_0^d(1+\eta )^k\|w_l\bar{f}\|_{\nu}^2d\eta
			&\lesssim  \int_0^d(1+\eta )^k\|\bar{f}\|_{\nu}^2d\eta +\int_0^d\int_{\R^3}(1+\eta )^k(\nu^0)^{-1}w_l^2g^2dvd\eta \nonumber\\
			&\quad + k \int_0^d\int_{\R^3}(1+\eta )^{k-1}|v_3|w_{l}^2\bar{f}^2dvd\eta \nonumber\\
			&\lesssim k \int_0^d(1+\eta )^{k-1}\|w_2\bar{f}\|_{\nu}^2d\eta +\int_0^d(1+\eta )^{2p_k}\|w_2g\|_{L^2_v}^2d\eta \nonumber\\
			& \quad  +\int_0^d(1+\eta )^k\|w_{l+2}g\|_{L^2_v}^2d\eta + k\int_0^d(1+\eta )^{k-1}\|w_{l+2}\bar{f}\|_{\nu}^2d\eta \nonumber\\
			&\lesssim k\int_0^d(1+\eta )^{k-1}\|w_{l+2}\bar{f}\|_{\nu}^2c\eta+\int_0^d(1+\eta )^{2p_k}\|w_{l+2}g\|_{L^2_v}^2d\eta ,
		\end{align}
		where $ p_k>\f{k}{2}+1.$ Using \eqref{D2.98}, \eqref{D2.106},  and induction arguments, one can deduce that
		\begin{align}\label{K3.90}
			\int_0^d(1+\eta )^k\|w_l\bar{f}\|_{\nu}^2d\eta &\lesssim_k\int_0^d(1+\eta )^{k-1}\|w_{l+2}\bar{f}\|_{\nu}^2d\eta+\int_0^d(1+\eta )^{2p_k}\|w_{l+2}g\|_{L^2_v}^2d\eta \nonumber\\
			&\lesssim_k\int_0^d(1+\eta )^{k-2}\|w_{l+4}\bar{f}\|_{\nu}^2d\eta+\int_0^d(1+\eta )^{2p_k}\|w_{l+4}g\|_{L^2_v}^2d\eta \nonumber\\
			&\lesssim_k \cdots  \lesssim_k\int_0^d\|w_{l+2k}\bar{f}\|_{\nu}^2d\eta+\int_0^d(1+\eta )^{2p_k}\|w_{l+2k}g\|_{L^2_v}^2d\eta \nonumber\\
			&\lesssim_k\int_0^d(1+\eta )^{2p_k}\|w_{l+2k+2}g\|_{L^2_v}^2d\eta,\quad p_k>\f{k}{2}+1.
		\end{align}
		Therefore the proof of Lemma \ref{lem3.11} is completed.
	$\hfill\Box$	
	
	\begin{lemma}\label{lem3.13}
		Let $(\phi_0,\phi_1,\phi_2,\phi_3)$ be the ones determined  in Lemma \ref{lem3.9}, then it holds that
		\begin{equation}\label{D2.110}
			\|(1+\eta )^kw_l\bar{f}\|_{L^\infty_{\eta ,v}}+|w_l\bar{f}|_{L^\infty(\gamma_+)}\leq C_k\|(1+\eta )^{q_k}w_{l+4k+4}g\|_{L^\infty_{\eta ,v}},\quad q_k>k+\f32,
		\end{equation}
		where  $k$ is a non-negative integer, and the constant $C_k$ is independent of $d$.
	\end{lemma}
	\noindent{\bf Proof.}
		Let $h_0=w_l\bar{f}$, then it holds that
		\begin{align}
			\begin{cases}
				v_3\partial_\eta  h_0+\nu^0 h_0 =K^0_{w_l}h_0+w_lg,\\
				h_0(\eta,v)|_{\gamma_-}=h_0(\eta,R_\eta v).
			\end{cases}
		\end{align}
		Applying Lemma \ref{lem3.4}, one has that
		\begin{align}\label{D2.112}
			\|w_l\bar{f}\|_{L^\infty_{\eta,v}}+|w_l\bar{f}|_{L^\infty(\gamma_+)}&\leq C\Big\{\int_0^d\int_{\R^3}\nu^0 \bar{f}^2 dvd\eta \Big\}^{\f12}+C\|(\nu^0)^{-1}w_lg\|_{L^\infty_{\eta,v}}\nonumber\\
			&\leq C\Big\{\int_0^d(1+\eta)^{2p_0}\|w_2g\|_{L^2_v}^2d\eta\Big\}^{\f12}+C\|(\nu^0)^{-1}w_lg\|_{L^\infty_{\eta,v}}\nonumber\\
			&\leq C\|(1+\eta )^{q_0}w_4g\|_{L^\infty_{\eta,v}}+C\|w_{l+3}g\|_{L^\infty_{\eta,v}}\nonumber\\
			&\leq C\|(1+\eta )^{q_0}w_{l+4}g\|_{L^\infty_{\eta,v}},\quad \text{for} \ q_0>p_0+\f12
		\end{align}
		where we have used \eqref{D2.94} to derive
		\begin{align}\label{D2.113}
			\int_0^d(1+\eta )^{2k}\|w_l\bar{f}\|_{\nu}^2d\eta &\lesssim_k \int_0^d (1+\eta)^{2p_{2k}}\|w_{l+4k+2}g\|_{L^2_v}^2\nonumber\\
			&\lesssim_k\|(1+\eta )^{q_k}w_{l+4k+4}g\|_{L^\infty_{\eta ,v}}^2\cdot\int_0^d(1+\eta )^{2p_{2k}-2q_{k}}d\eta \int_{\R^3}w_2^{-2}dv\nonumber\\
			&\lesssim_k\|(1+\eta )^{q_{k}}w_{l+4k+4}g\|_{L^\infty_{\eta ,v}}^2,\quad \text{for} \ q_{k}>p_{2k}+\f12.
		\end{align}

		Let $h_k=(1+\eta )^kw_l\bar{f}$, then it holds that
		\begin{align}
			v_3\partial_\eta  h_k+\nu^0 h_k=K^0_{w_l}h_k+k(1+\eta )^{k-1}v_3w_l\bar{f}+(1+\eta )^kw_lg.
		\end{align}
		Applying Lemma \ref{lem3.4} and  \eqref{D2.113}, one gets that
		\begin{align}\label{D2.115}
			\|h_k\|_{L^\infty_{\eta,v}}+|h_k|_{L^\infty(\gamma_+)}&\leq Ck\|(1+\eta )^{k-1}(\nu^0)^{-1}v_3w_l\bar{f}\|_{L^\infty_{\eta,v}}+C\|(1+\eta )^k(\nu^0)^{-1}w_{l}g||_{L^\infty_{\eta,v}}\nonumber\\
			&\quad+C\|(1+\eta )^k(\nu^0)^{\f12}\bar{f}\|_{L^2_{\eta,v}}\nonumber\\
			&\leq C_k\|(1+\eta )^{k-1}w_{l+4}\bar{f}\|_{L^\infty_{\eta,v}} +C_k\|(1+\eta )^{q_k}w_{\max\{4k+4,l+3\}}g\|_{L^\infty_{\eta,v}},
		\end{align}
		where $q_k>p_{2k}+\f12$. Using \eqref{D2.112}, \eqref{D2.115} and induction arguments, one obtains that
		\begin{align*}
			&\quad \|(1+\eta )^kw_l\bar{f}\|_{L^\infty_{\eta,v}}+|(1+\eta )^kw_l\bar{f}|_{L^\infty(\gamma_+)}\nonumber\\
			&\leq C_k\|(1+\eta )^{k-1}w_{l+4}\bar{f}\|_{L^\infty_{\eta,v}} +C_k\|(1+\eta )^{q_k}w_{\max\{4k+4,l+3\}}g\|_{L^\infty_{\eta,v}}\nonumber\\
			&\leq C_k\|(1+\eta )^{k-2}w_{l+8}\bar{f}\|_{L^\infty_{\eta,v}}+C_k\|(1+\eta )^{q_k}w_{\max\{4k+4,l+7\}}g\|_{L^\infty_{\eta,v}}\nonumber\\
			&...\nonumber\\
			&\leq C_k\|w_{l+4k}\bar{f}\|_{L^\infty_{\eta,v}}+C_k\|(1+\eta )^{q_k}w_{\max\{4k+4,l+4k\}}g\|_{L^\infty_{\eta,v}}\nonumber\\
			&\leq C_k\|(1+\eta )^{q_k}w_{l+4k+4}g\|_{L^\infty_{\eta,v}},\quad q_k>p_{2k}+\f12.
		\end{align*}
		Recall the range of $p_k$ in \eqref{D2.94}, then the proof of Lemma \ref{lem3.13} is finished.
	$\hfill\Box$	
	
	\
	
	With the help of decay estimate in Lemma \ref{lem3.13}, we shall prove Theorem \ref{thm1.1} by taking the limit $d\rightarrow\infty$. From now on, we denote the solution $\bar{f}(\eta ,v)$ of \eqref{D2.79} as $\bar{f}_d(\eta ,v)$ to emphasize the dependence on $d$. We denote
	\begin{equation*}
		\tilde{f}(\eta ,v)=\bar{f}_{d_2}(\eta ,v)-\bar{f}_{d_1}(\eta ,v), \quad 1\leq d_1\leq d_2.
	\end{equation*}
	Then $\tilde{f}$ satisfies the following equation
	\begin{equation}\label{D2.117}
		\begin{cases}
			v_3\partial_\eta  \tilde{f}+\FL_0 \tilde{f}=0,\quad \eta \in[0,d_1],\ v\in\R^3, \\
			\tilde{f}(0,v)|_{v_3>0}=\tilde{f}(0,Rv).
		\end{cases}
	\end{equation}

	\subsection{Proof of Theorem \ref{thm1.1}} We  divide the proof into two steps.\vspace{1mm}
	
	\noindent{\it Step 1. Convergence in $L^2$-norm. } Multiplying \eqref{D2.117} by $\tilde{f}$ and integrating on $[0,d_1]\times\R^3$, one obtains that
	\begin{align}\label{D2.118}
		&\int_0^{d_1}\int_{\R^3}\nu^0|(\mathbf{I}-\FP_0)\tilde{f}(\eta ,v)|^2 dvd\eta \nonumber\\
		&\leq C\int_{\R^3} |v_3|\cdot |\tilde{f}(d_1,v)|^2dv \leq C \big\{\|w_l\bar{f}_{d_2}(d_1)\|^2_{L^\infty_{v}}+|w_l\bar{f}_{d_1}(d_1)|^2_{L^\infty(\gamma_+)}\big\}\nonumber\\
		&\leq C\|(1+\eta )^{\mathfrak{q}}w_{l+4}g\|^2_{L^\infty_{\eta ,v}}\cdot d_1^{-2}, \quad \mathfrak{q}\geq3.
	\end{align}
	
	We still need to control the macroscopic part. Denote
	\begin{align*}
		\FP_0\tilde{f}=[\tilde{a}(\eta )+\tilde{b}_{1}(\eta ) (v_1-\fu^0_1)+\tilde{b}_{2}(\eta ) (v_2-\fu^0_2)+ \tilde{c}(\eta ) (|v-\fu^0|^2-3T^0)]\sqrt{\mu_0}.
	\end{align*}
	Similar as in Lemma \ref{lem3.10}, we can obtain
	\begin{align}
		\left(
		\begin{array}{c}
			\int_{\R^3} v_3\tilde{f}(\eta ,v)\cdot \FL_0^{-1} (\mathcal{A}_{31}^0) dv\\[1mm]
			\int_{\R^3} v_3\tilde{f}(\eta ,v)\cdot \FL_0^{-1} (\mathcal{A}_{32}^0) dv\\[1mm]
			\int_{\R^3} v_3\tilde{f}(\eta ,v)\cdot \FL_0^{-1} (\mathcal{B}_3^0) dv\\
		\end{array}
		\right)
		=\left(
		\begin{array}{c}
			\int_{\R^3} v_3\tilde{f}(d_1,v)\cdot \FL_0^{-1} (\mathcal{A}_{31}^0)(d_1) dv\\[1mm]
			\int_{\R^3} v_3\tilde{f}(d_1,v)\cdot \FL_0^{-1} (\mathcal{A}_{32}^0)(d_1) dv\\[1mm]
			\int_{\R^3} v_3\tilde{f}(d_1,v)\cdot \FL_0^{-1} (\mathcal{B}_3^0)(d_1) dv\\
		\end{array}
		\right),\nonumber
	\end{align}
	which, together with \eqref{D2.87}-\eqref{D2.89}, yields that
	\begin{align}\label{D2.119}
		|(\mu(T^0) \tilde{b}_{1},\mu(T^0)  \tilde{b}_{2},\kappa(T^0) \tilde{c})(\eta )|&\leq C\Big\{\|w_l\bar{f}_{d_2}(d_1)\|_{L^\infty_{v}}+|w_l\bar{f}_{d_1}(d_1)|_{L^\infty(\g_+)}\Big\}+C\|(\mathbf{I}-\FP_0)\tilde{f}(\eta )\|_{\nu}.
	\end{align}
	Integrating \eqref{D2.119} over $[0,d_1]$, using \eqref{D2.110}, \eqref{D2.118}, one has
	\begin{align}\label{D2.120}
		\int_{0}^{d_1} |(\tilde{b}_{1},\tilde{b}_{2},\tilde{c})(\eta )|^2 d\eta \leq C\|(1+\eta )^\mathfrak{q}w_{l+4}g\|^2_{L^\infty_{\eta ,v}} \cdot d_1^{-1}, \quad \mathfrak{q}\geq 3.
	\end{align}

	Multiplying \eqref{D2.117} by $v_3\sqrt{\mu_0}$, we have that
	\begin{align}\nonumber
		\frac{d}{d\eta } \int_{\R^3} \tilde{f}(\eta ,v)\cdot v_3^2 \sqrt{\mu_0} dv=0.
	\end{align}
	Integrating the above equation over $[\eta ,d]$ and using \eqref{D2.86}, one obtains
	\begin{align}\label{D2.121}
		\tilde{a}(\eta )=-2T^0\tilde{c}(\eta )-\f{1}{\rho^0T^0}\int_{\R^3}   (\mathbf{I}-\FP_0)\tilde{f}(\eta ,v)\cdot v_3^2 \sqrt{\mu_0} dv+\f{1}{\rho^0T^0}\int_{\R^3}\tilde{f}(d_1,v)\cdot v_3^2\sqrt{\mu_0}  dv.
	\end{align}
	Using  \eqref{D2.110} \eqref{D2.118} and \eqref{D2.120}, one can get that
	\begin{align*}
		\int_{0}^{d_1} |\tilde{a}(\eta )|^2 d\eta  \leq C\|(1+\eta )^\mathfrak{q}w_{l+4}g\|^2_{L^\infty_{\eta ,v}} \cdot d_1^{-1}, \quad \text{for}\  \mathfrak{q}\geq 3,
	\end{align*}
	which, together with \eqref{D2.118} and \eqref{D2.120}, yields that
	\begin{align}\label{D2.122}
		\int_0^{d_1}\int_{\R^3}\nu^0|\tilde{f}(\eta ,v)|^2 dvd\eta \leq C\|(1+\eta )^\mathfrak{q}w_{l+4}g\|^2_{L^\infty_{\eta ,v}} \cdot d_1^{-1}, \quad \mathfrak{q}\geq 3.
	\end{align}
	
	\vspace{1mm}
	
	\noindent{\it Step 2. Convergence in $L^\infty$-norm.}
	We shall use $t_{k}=t_{k}(t,\eta ,v), X_{cl}(s;t,\eta ,v), \eta _k=\eta _k(\eta,v)$ to be the back-time cycles defined for domain $[0,d_1]\times \R^3$. For later use, we denote $\tilde{h}:=w_l\tilde{f}$. Let $(\eta,v)\in [0,d_1]\times\R^3\backslash (\g_0\cup\g_{-})$, it follows from \eqref{D2.117} that
	\begin{align}\label{D2.123}
		\tilde{h}(\eta ,v)&=e^{-\hat{\nu}(v) (t-t_{k})} \tilde{h}(d_1,v_{k-1}) +\sum_{i=0}^{k-1} \int_{t_{i+1}}^{t_i}e^{-\hat{\nu}(v) (t-s)} (1+|v|^2)^{\f{|\kappa|}{2}} K^0_{w_l} \tilde{h}(X_{cl}(s), v_i)ds,
	\end{align}
	with $k=1$ for $v_{0,3}<0$, and $k=2$ for $v_{0,3}>0$. We   will  use this summation convention in the following of this lemma.
	We always have
	\begin{align}\label{D2.124}
		|e^{-\hat{\nu}(v) (t-t_{k})} \tilde{h}(d_1,v_{k-1}) |
		&\leq C \Big(\|{w_l}\bar{f}_{d_2}(d_1)\|_{L^\infty_v}+|{w_l}\bar{f}_{d_1}(d_1)|_{L^\infty(\gamma_+)}\Big)\nonumber\\
		&\leq C\|(1+\eta )^\mathfrak{q} w_{l+4}g\|_{L^\infty_{\eta ,v}}\cdot d_1^{-1},\quad \mathfrak{q}\geq3.
	\end{align}

	For the second term on RHS of \eqref{D2.123}, we use \eqref{D2.123} again to obtain
	\begin{align}\label{D2.125}
		&\quad \left|\sum_{i=0}^{k-1} \int_{t_{i+1}}^{t_i}e^{-\hat{\nu}(v) (t-s)}(1+|v|^2)^{\f{|\kappa|}{2}}  K^0_{w_l} \tilde{h}(X_{cl}(s), v_i)ds\right|\nonumber\\
		&=\left|\sum_{i=0}^{k-1} \int_{t_{i+1}}^{t_i}e^{-\hat{\nu}(v) (t-s)}(1+|v|^2)^{\f{|\kappa|}{2}}  K^{0,c}_{w_l} \tilde{h}(X_{cl}(s), v_i)ds\right|\nonumber\\
		&+\left|\sum_{i=0}^{k-1} \int_{t_{i+1}}^{t_i}e^{-\hat{\nu}(v) (t-s)}(1+|v|^2)^{\f{|\kappa|}{2}}  K^{0,m}_{w_l} \tilde{h}(X_{cl}(s), v_i)ds\right|\nonumber\\
		&\leq \f14\|\tilde{h}\|_{L^\infty_{\eta ,v}}+C\|(1+\eta )^\mathfrak{q}w_{l+4}g\|_{L^\infty_{\eta ,v}}\cdot d_1^{-1}\nonumber\\
		&+\Big|\sum_{i=0}^{k-1} \int_{t_{i+1}}^{t_i}e^{-\hat{\nu}(v) (t-s)}(1+|v|^2)^{\f{|\kappa|}{2}}  \int_{\R^3}k^{0,c}_{w_l}(v_i,v')(1+|v'|^2)^{\f{|\kappa|}{2}}\nonumber\\
		&\quad \times\sum\limits_{j=0}^{k-1}\int_{t_{j+1}'}^{t_j'}e^{-\hat{\nu}(v')(s-s_1)}\int_{\R^3}k^{0,c}_{w_l}(v_j',v'') \tilde{h}(X_{cl}(s_1), v'')dv''ds_1dv'ds\Big|.
	\end{align}
	where we have used \eqref{D2.124} and denote $X'_{cl}(s_1)=X_{cl}(s_1;s,X_{cl}(s), v')$, $t'_j=t'_j(s_1;s, X_{cl}(s), v')$ and $v'_j$ to be the back-time cycle of $(s,X_{cl}(s), v')$. Then, by the same arguments as in Lemma \ref{lem3.4}, we get
	\begin{align}\label{D2.126}
		&\|\tilde{h}\|_{L^\infty([0,d_1]\times \R^3)}+|\tilde{h}(0)|_{L^\infty(\gamma_+)}\nonumber\\
		&\leq \f12 (\|\tilde{h}\|_{L^\infty([0,d_1]\times \R^3)}+|\tilde{h}(0)|_{L^\infty(\g_+)} )\nonumber\\
		&\quad+Cd_1^{-1}\|(1+\eta )^\mathfrak{q}w_{l+4}g\|_{L^\infty_{\eta ,v}}+C\|(\nu^0)^{\f12}\tilde{f}\|_{L^2([0,d_1]\times\R^3)}\nonumber\\
		&\leq Cd_1^{-\f12}\|(1+\eta )^\mathfrak{q}w_{l+4}g\|_{L^\infty_{\eta ,v}},\quad \mathfrak{q}\geq3.
	\end{align}
	With the help of \eqref{D2.126}, 
	there exists a function $f(\eta ,v)$ with $(\eta ,v)\in\R_+\times \R^3$ so that $\|w(\bar{f}_{d}-f)\|_{L^\infty([0,d]\times\R^3)}\rightarrow0$ as $d \rightarrow \infty$. The uniform bound  \eqref{D2.5} follows from \eqref{D2.110} and the strong convergence in $L^\infty_{\eta ,v}$. It is  direct to see that $f(\eta,v)$ solves \eqref{D2.2}. The continuity of $f$  follows directly from the $L^\infty_{x,v}$-convergence and the continuity of $\bar{f}_{d}$.
	
	\smallskip
	
	For the uniqueness, let $\mathbf{f}_1,\mathbf{f}_2$ be two solutions of \eqref{D2.2} with the bound \eqref{D2.5}, then it holds that
	\begin{align}\label{D2.127}
		\begin{cases}
			v_3 \partial_\eta  (\mathbf{f}_1-\mathbf{f}_2)+\FL_0 (\mathbf{f}_1-\mathbf{f}_2)=0,\\
			\mathbf{f}_i(0,v)|_{v_3>0}=\mathbf{f}_i(0,Rv),\ i=1,2,\\
			\lim_{\eta \rightarrow\infty}\mathbf{f}_i(\eta ,v)=0,\ i=1,2.
		\end{cases}
	\end{align}
	Multiplying \eqref{D2.127} by $(\mathbf{f}_1-\mathbf{f}_2)$, it is direct to prove that
	\begin{align}\nonumber
		\int_0^\infty \|(\mathbf{I}-\FP_0)(\mathbf{f}_1-\mathbf{f}_2)\|_{\nu}^2 d\eta =0.
	\end{align}
	That is, $(\mathbf{f}_1-\mathbf{f}_2)=\FP_0 (\mathbf{f}_1-\mathbf{f}_2)$. Then by the same arguments as \eqref{D2.91}-\eqref{D2.92}, one has that
	\begin{align}\nonumber
		\int_0^\infty \|\FP_0(\mathbf{f}_1-\mathbf{f}_2)\|_{L^2_v}^2 d\eta=0.
	\end{align}
	Thus, we prove $\mathbf{f}_1\equiv\mathbf{f}_2$. 
	
	Finally, let $\mathfrak{f}:=f+\Upsilon(\eta)\,f_{b}(v)$, then it direct to know that $\mathfrak{f}$ solves \eqref{D2.1}. The proof of Theorem \ref{thm1.1} is completed.
	
	\
	
\section{Hilbert Expansions for Boltzmann Equation of Soft Potentials}\label{section4}
In this section, we aim to construct the solutions of Boltzmann equation of soft potentials through Hilbert expansion with multi-scales.

\subsection{Linear parts of Hilbert expansion}
	In this section, we shall construct the soft Boltzmann solutions in the form of Hilbert expansion with multi-scales. Recall $\varpi_{\mathfrak{k}}$ in \eqref{c1.31}, we define the velocity  weight functions
	\begin{align}\label{ewf}
		\tw_{\k_i}(v)=\varpi_{\k_i}(v) \mu^{-\fa_i},
		\quad \fw_{\bar{\k}_i}(v)=\varpi_{\bar{\k}_i}(v) \mu_0^{-\fa_i} \,
		\mbox{and}\
		\fw_{\hat{\k}_i}(v)=\varpi_{\hat{\k}_i}(v) \mu_0^{-\fa_i},
	\end{align}
	for constants $\k_i, \bar{\k}_i, \hat{\k}_i\geq 0,\, 1\leq i\leq N$ and $0\leq \fa_i <\frac{1}{2}$. Note that the weight function $\tw_{\k_i}$ depends on $(t,x)$, while $\fw_{\bar{\k}_i}$ and $\fw_{\hat{\k}_i}$ depend on $(t,x_{\sp} )$.  For later use, we define
	\begin{align}\nonumber
		\quad \hat{x}=(x_{\sp} ,\eta)\in \R_+^3, \quad \nabla_{\hat{x}}:=(\nabla_{\sp} ,\pa_{\eta}),
	\end{align}
	and recall  $\bar{x}=(x_{\sp} ,y) \in \R_+^3,\ \nabla_{\bar{x}}=(\nabla_{\sp} , \pa_y)$ , and the weighted $L^2_{l}$-norm with $(1+y)^l$ weight.
	\begin{proposition} \label{prop5.1}
		Let $\tau^\d>0$ be the life-span of compressible Euler equations. Let $0\leq \fa_i<\frac12$ in \eqref{ewf} and $\fa_i>\fa_{i+1}$. Let  $s_0, s_i, \bar{s}_i, \hat{s}_i,\zeta_i\in \mathbb{N}_+$, $\k_i, \bar{\k}_{i}, \hat{\k}_i \in \R_+$ for $1\leq i\leq N$; and define  $l_j^i:=\bar{l}_i+2(\bar{s}_i-j)$ for $1\leq i \leq N, \  0\leq j\leq \bar{s}_i$. For these parameters, we have chosen $s_i,\bar{s}_i,\hat{s}_i$ such that
		\begin{align}\label{H47}
			&s_0\geq s_1+\fb+6,\quad s_1=\bar{s}_1=\hat{s}_1\gg1;\nonumber\\
			&s_1>s_i>\bar{s}_i>\hat{s}_i\geq s_{i+1}>\bar{s}_{i+1}>\hat{s}_{i+1}\geq ...\gg1,\ i=2,...,N-1;\nonumber\\
			&s_{i+1}\leq \min\{\hat{s}_i,\f12\bar{s}_i-3\},\ \bar{s}_{i+1}\leq s_{i+1}-8-\fb,\ \hat{s}_{i+1}\leq \f12\bar{s}_{i+1}-2-\fb,\ i=1,...,N-1,
		\end{align}
		and taken $l_j^i=\bar{l}_j+2(\bar{s}_i-j)$ with $0\leq j\leq \bar{s}_i$ so that
		\begin{align}\label{H48}
			l_j^N\gg2\fb \quad \mbox{and} \quad l_j^i\geq 2l_j^{i+1}+18+2\fb,\ \mbox{for}\ 1\leq i\leq N-1,
		\end{align}
		and
		\begin{align}\label{H49}
			&\kappa_i\gg \bar{\k}_i\gg\hat{\k}_i\gg\k_{i+1}\gg\bar{\k}_{i+1}\gg\hat{\k}_{i+1}\gg1,\nonumber\\
			&\zeta_{i+1}-\zeta_i\geq \fb+3\quad  \text{and }\quad \zeta_1\gg\zeta_2...\gg\zeta_i...\gg \fb.
		\end{align}
	
		Let $(\rho_{i}, u_i, \theta_i)(0)$ be the initial data for interior expansions, and $(\bar{u}_{i,\sp},\bar{\theta}_{i})(0)$ be the initial data of viscous boundary layer. Assume
		\begin{align}\label{7.1-0}
			\sum_{i=0}^{N}\Big\{\sum_{\gamma+\beta\leq s_{i}}\|\pa_t^\gamma\nabla_x^{\beta}(\rho_{i}, u_{i},  \theta_{i})(0)\|_{L^2_{x}}+\sum_{j=0}^{\bar{s}_{i}} \sum_{j=2\gamma+\beta}   \|\pa_t^{\gamma} \nabla_{\bar{x}}^{\beta} (\bar{u}_{i,\sp},\bar{\theta}_{i})(0)\|^2_{L^2_{l_j^{i}}}
			\Big\}<\infty.
		\end{align}
		And we also assume that the compatibility conditions for initial data $(\rho_{i}, u_i, \theta_i)(0)$ and $(\bar{u}_{i,\sp},\bar{\theta}_{i})(0)$ are satisfied. Then there exist solutions $F_i=\sqrt{\mu}f_i, \, \bar{F}_i=\sqrt{\mu_0}\bar{f}_i, \, \hat{F}_i=\sqrt{\mu_0}\hat{f}_i$ to  interior expansions \eqref{c1.7}, viscous boundary layer \eqref{c1.15} and Knudsen layer solutions \eqref{c1.22} over the time interval  $t\in[0,\tau^\d]$ so that the specular boundary condition is satisfied in the following form:
		\begin{align*}
			 (F_i+\bar{F}_i+\hat{F}_i)(t,x_\sp,0,v_\sp,v_3)|_{v_3>0}=(F_i+\bar{F}_i+\hat{F}_i)(t,x_\sp,0,v_\sp,-v_3).
		\end{align*}
	Moreover, it holds that
		\begin{align}\label{7.2-0}
			&\sup_{t\in[0,\tau^\d]} \sum_{i=1}^N\Bigg\{\sum_{\gamma+\beta\leq s_i} \|\tw_{\k_i}\pa_t^\gamma \nabla_x^{\beta} f_{i}(t)\|_{L^2_xL^\infty_v}
			+\sum_{j=0}^{\bar{s}_i}\sum_{j=2\gamma+\beta} \|\fw_{\bar{\k}_i}\pa_t^\gamma \nabla_{\bar{x}}^{\beta} \bar{f}_{i}(t)\|_{L^2_{l^i_j}L^\infty_v}\nonumber\\
			&\qquad\qquad\qquad+\sum_{\gamma+\beta\leq \hat{s}_i} \| (1+\eta)^{\zeta_i}\fw_{\hat{\k}_i}\pa_t^\gamma \nabla_{\sp} ^{\beta} \hat{f}_{i}(t)\|_{L^\infty_{\hat{x},v}\cap L^2_{x_{\sp} }L^\infty_{\eta,v}}\Bigg\}\nonumber\\
			&\leq C\Bigg( \tau^\d, \|(\varphi_0,\Phi_0, \vartheta_0)\|_{H^{s_0}}+\sum_{i=0}^{N}\sum_{\gamma+\beta\leq s_{i}}\|\pa_t^\gamma\nabla_x^{\beta}(\rho_{i}, u_{i},  \theta_{i})(0)\|_{L^2_{x}}\nonumber\\
			&\qquad\qquad\qquad+\sum_{i=0}^{N}\sum_{j=0}^{\bar{s}_{i}} \sum_{j=2\gamma+\beta}   \|\pa_t^{\gamma} \nabla_{\bar{x}}^{\beta} (\bar{u}_{i,\sp},\bar{\theta}_{i})(0)\|^2_{L^2_{l_j^{i}}}
			\Bigg).
		\end{align}
	\end{proposition}
\noindent{\bf Proof.}
  With the help of Proposition \ref{thmA.1} and Theorem \ref{thm1.1}, by similar arguments as in \cite[Proposition 5.1]{GHW}, one can construct $f_i,\bar{f}_i$ and $\hat{f}_i$.

  Here we explain a little bit on how to use Proposition \ref{thmA.1} and Theorem \ref{thm1.1} for soft potential cases. Noting $f_i,\bar{f}_i$ are smooth on $(t,x,v)$ and $(t,\bar{x},v)$, then by using Proposition \ref{thmA.1}, one can always get the exponential decay on $v$, i.e.
  \begin{align*}
  	|\partial_t\nabla_xf_i|\lesssim\mu^{\f{q}{2}},\quad |\partial_t\nabla_{\bar{x}}\bar{f}_i|\lesssim\mu_0^{\f{q}{2}} \quad  \text{for} \,\, q\in(0,1).
  \end{align*}
With the help of Theorem \ref{thm1.1}, we can construct the solutions for Knudsen boundary layers $\hat{f}_i$ with enough polynomial space decay estimate.
$\hfill\Box$	
	
\subsection{Estimates on the remainder}
We first consider the $L^2$-energy estimate. Recall the definition of $f^\v_R$ in \eqref{c1.28}, we rewrite the equation of $f^\v_R$ as
\begin{align}\label{H1}
	&\dis \partial_t f^\v_R+v\cdot\nabla_x f^\v_R+\frac{1}{\v^2}\mathbf{L}f^\v_R\nonumber\\
	 &=-\frac{\{\partial_t+v\cdot\nabla_x\}\sqrt{\mu}}{\sqrt{\mu}}f^\v_R+\v^{3}\frac1{\sqrt{\mu}}Q(\sqrt{\mu}f^\v_R,\sqrt{\mu}f^\v_R)\nonumber\\
	&\,\,+\sum_{i=1}^N\v^{i-2}\frac1{\sqrt{\mu}}\Big\{Q(F_i+\bar{F}_i+\hat{F}_i,\sqrt{\mu}f^\v_R)
	+Q(\sqrt{\mu}f^\v_R, F_i+\bar{F}_i+\hat{F}_i)\Big\}\nonumber\\
	& \,\, + \frac1{\sqrt{\mu}}R^\v+\frac1{\sqrt{\mu}}\bar{R}^\v+\frac1{\sqrt{\mu}}\hat{R}^\v,
\end{align}
where
\begin{align}\label{H2}
	R^\v&=-\v^{N-6}\{\partial_t+v\cdot \nabla_x\} (F_{N-1}+\v F_N) +\v^{N-6}\sum_{\substack{i+j\geq N+1\\1\leq i,j\leq N}}\v^{i+j-N-1}Q(F_{i},F_{j}),
\end{align}
\begin{align}\label{H3}
	&\bar{R}^\v=-\v^{N-6}\{\partial_t+v_\sp\cdot \nabla_\sp\} (\bar{F}_{N-1}+\v \bar{F}_N)-\v^{N-6} v_3\partial_y \bar{F}_N\nonumber\\
	&\qquad+\v^{N-6}\sum_{\substack{j+l\geq N+1\\ 1\leq j\leq N, \,  1\leq l\leq \fb}} \v^{l+j-N-1}\cdot\frac{y^l}{l!} \big[Q(\partial_3^l\mu_0, \bar{F}_{j})+Q( \bar{F}_{j},\partial_3^l\mu_0) \big]\nonumber\\
	&\qquad+\v^{N-6}\sum_{\substack{i+j\geq N+1\\   1\leq i,j\leq N}} \v^{i+j-N-1} \big[ Q(F_i^0,\bar{F}_j)+Q(\bar{F}_j,F_i^0)+Q(\bar{F}_i,\bar{F}_j) \big]\nonumber\\
	&\qquad+\v^{N-6}\sum_{\substack{i+j+l\geq N+1\\ 1\leq i, j\leq N, \,  1\leq l\leq \fb}} \v^{i+j+l-N-1}\cdot \frac{y^l}{l!} \big[Q(\partial_3^l F_i^0, \bar{F}_{j})+Q( \bar{F}_{j},\partial_3^l F_i^0) \big]\nonumber\\
	&\qquad +  \v^{\fb-5}\frac{y^{\fb+1}}{(\fb+1)!} \sum_{j=1}^N \v^{j-1} [Q(\partial_3^{\fb+1}\tilde{\mu}, \bar{F}_{j})+Q( \bar{F}_{j},\partial_3^{\fb+1}\tilde{\mu})]\nonumber\\
	&\qquad +\v^{\fb-4}\frac{y^{\fb+1}}{(\fb+1)!} \sum_{i,j=1}^N \v^{i+j-2} \big[Q(\partial_3^{\fb+1}\mathfrak{F}_i, \bar{F}_{j})+Q( \bar{F}_{j},\partial_3^{\fb+1}\mathfrak{F}_i)\big],
\end{align}
and
\begin{align}\label{H4}
	&\hat{R}^\v=-\v^{N-6}\{\partial_t+v_\sp\cdot \nabla_\sp\} (\hat{F}_{N-1}+\v \hat{F}_N)\nonumber\\
	&\quad+\v^{N-6}\sum_{\substack{j+2l\geq N+1\\  1\leq j\leq N, 1\leq l\leq \fb}} \v^{j+2l-N-1}\cdot \frac{\eta^l}{l!} \big[ Q(\partial_3^l\mu_0, \hat{F}_{j})+Q( \hat{F}_{j},\partial_3^l\mu_0) \big]\nonumber\\
	&\quad+\v^{N-6}\sum_{\substack{i+j\geq N+1\\  1\leq i,j\leq N}} \v^{i+j-N-1} \big[Q(F_i^0+\bar{F}_i^0,\hat{F}_j)+Q(\hat{F}_j,F_i^0+\bar{F}_i^0)+Q(\hat{F}_i,\hat{F}_j)\big]\nonumber\\
	&\quad+\v^{N-6}\sum_{\substack{i+j+2l\geq N+1\\  1\leq i,j\leq N, 1\leq l\leq \fb}} \v^{i+j+2l-N-1}\cdot \frac{\eta^l}{l!} \big[ Q(\partial_3^l F_i^0, \hat{F}_{j})+Q( \hat{F}_{j},\partial_3^l F_i^0) \big]\nonumber\\
	&\quad+\v^{N-6}\sum_{\substack{i+j+l\geq N+1\\  1\leq i,j\leq N, 1\leq l\leq \fb}} \v^{i+j+l-N-1}\cdot \frac{\eta^l}{l!} \big[ Q(\partial_y^l \bar{F}_i^0, \hat{F}_{j})+Q( \hat{F}_{j},\partial_y^l \bar{F}_i^0) \big]\nonumber\\
	&\quad + \v^{2\fb -4} \frac{\eta^{\fb+1}}{(\fb+1)!} \sum_{j=1}^N \v^{j-1} \big[ Q(\partial_3^{\fb+1}\tilde{\mu}, \hat{F}_{j})+Q( \hat{F}_{j},\partial_3^{\fb+1}\tilde{\mu}) \big]\nonumber\\
	&\quad + \v^{2\fb-3}\frac{\eta^{\fb+1}}{(\fb+1)!} \sum_{i,j=1}^N \v^{i+j-2} \big[ Q(\partial_3^{\fb+1}\mathfrak{F}_i, \hat{F}_{j})+Q( \hat{F}_{j},\partial_3^{\fb+1}\mathfrak{F}_i) \big]\nonumber\\
	&\quad + \v^{\fb-4}\frac{\eta^{\fb+1}}{(\fb+1)!} \sum_{i,j=1}^N \v^{i+j-2} \big[ Q(\partial_3^{\fb+1}\bar{\mathfrak{F}}_i, \hat{F}_{j})+Q( \hat{F}_{j},\partial_3^{\fb+1}\bar{\mathfrak{F}}_i) \big],
\end{align}
where $\pa_3^l\mu_0, \pa_3^{\fb+1}\tilde{\mu}$, $\pa_3^lF_i^0, \pa_3^{\fb+1}\mathfrak{F}_i$ and $\pa_y^l\bar{F}_i^0, \pa_y^{\fb+1}\bar{\mathfrak{F}}_i$ are  defined in \eqref{c1.18}, \eqref{c1.23}. From Proposition \ref{prop5.1}, we know that $f^\v_R$ satisfies specular reflection boundary conditions
\begin{equation}\label{H5}
	f^\v_R(t,x_1,x_2,0,v_1,v_2,v_3)|_{v_3>0}=f^\v_R(t,x_1,x_2,0,v_1,v_2,-v_3).
\end{equation}

\smallskip

\begin{lemma}\label{lem5.1}
	Recall $\alpha$ in \eqref{c1.30}. Let $0<\frac{1}{2\alpha}(1-\alpha)<\fa_i<\frac12$,  $\mathfrak{k}\geq 18$,  $N\geq 6$ and $\fb\geq 5$. Let $\tau^\d>0$ be the life span of compressible Euler solution, then there exists a suitably small constant $\v_0>0$  such that for all $\v\in (0,\v_0)$,  it holds that
	\begin{align}\label{H6}
		 &\frac{d}{dt}\|f^\v_R(t)\|_{L^2}^{2}+\frac{c_{0}}{2\v^2}\|\{\mathbf{I-P}\}f^\v_R(t)\|_{\nu}^{2}\nonumber\\
		&\leq C \big\{1+ \v^{8}\|h_R^\v(t)\|_{L^\infty}^2 \big\}
		\cdot(\|f_R^\v(t)\|_{L^2}^{2}+1),\,\, \mbox{for}\,\, t\in[0,\tau^\d].
	\end{align}
\end{lemma}

\noindent{\bf Proof.} Multiplying \eqref{H1} by $f^\v_R$ and integrating   over $\mathbb{R}^3_+\times\mathbb{R}^3$, one obtains that
\begin{align}\label{H7}
	&\frac12 \frac{d}{dt} \|f^{\v}_R\|_{L^2}^2+\frac{c_0}{\v^2} \|\{\mathbf{I-P}\}f^{\v}_R\|_{\nu}^2\nonumber\\
	&= -\int_{\mathbb{R}^3_+}\intr \frac{\{\partial_t+v\cdot\nabla_x\}\sqrt{\mu}}{\sqrt{\mu}} |f^\v_R|^2  + \v^3\int_{\mathbb{R}^3_+}\intr \frac1{\sqrt{\mu}}Q(\sqrt{\mu}f^\v_R,\sqrt{\mu}f^\v_R) f^{\v}_R  \nonumber\\
	&\,\, +\int_{\mathbb{R}^3_+}\intr\sum_{i=1}^N \v^{i-2}\frac1{\sqrt{\mu}}\Big\{Q(F_i+\bar{F}_i+\hat{F}_i,\sqrt{\mu}f^\v_R)
	+Q(\sqrt{\mu}f^\v_R, F_i+\bar{F}_i+\hat{F}_i)\Big\} f^\v_R  \nonumber\\
	&\,\,+ \int_{\mathbb{R}^3_+}\intr \left\{\frac1{\sqrt{\mu}}R^\v+\frac1{\sqrt{\mu}}\bar{R}^\v+\frac1{\sqrt{\mu}}\hat{R}^\v\right\} f^\v_R  ,
\end{align}
where we have used \eqref{H5} so that the boundary term vanishes. \\

Recall the definition $h^{\v}_R$ in \eqref{c1.31}. For any $\lambda>0$, motivated by  \cite{Guo Jang Jiang}, we  take $\mathfrak{k}\geq 18$ to get
\begin{align*}
	&\int_{\mathbb{R}^3_+}\intr \frac{\{\partial_t+(v\cdot\nabla_x)\}\sqrt{\mu}}{\sqrt{\mu}} |f^\v_R|^2 dv dx \nonumber\\
	&\leq C \int_{\mathbb{R}^3_+}\intr |(\nabla_x\rho,\nabla_x\fu,\nabla_x T)| (1+|v|)^3 |f^\v_R|^2 dv dx \nonumber\\
	&\leq C\left\{\int_{\mathbb{R}^3_+}\int_{|v|\geq \frac{\lambda}{\v^{1/3}}}+  \int_{\mathbb{R}^3_+}\int_{|v|\leq \frac{\lambda}{\v^{1/3}}}\right\} (\cdots) dvdx\nonumber\\
	&\leq C\frac{\lambda}{\v^2} \|\{\mathbf{I-P}\}f^{\v}_R\|_{\nu}^2+C_\lambda (1+\v^4 \|h^\v_R\|_{L^\infty})  \|f^\v_R\|_{L^2},
\end{align*}
where we have used
\begin{align*}
	&\int_{\R^3}\int_{|v|\leq\f{\lambda}{\v^{1/3}}}|\nabla_x(\rho,\fu,T)|(1+|v|)^3|f^\v_R|^2dvdx\nonumber\\
	&\leq \int_{\R^3}\int_{|v|\leq\f{\lambda}{\v^{1/3}}}|\nabla_x(\rho,\fu,T)|(1+|v|)^3\Big\{|\FP f^\v_R|^2+|(\FI-\FP)f^\v_R|^2\Big\}dvdx\nonumber\\
	&\leq C_\lambda\|f_R^\v\|_{L^2}^2+C\|(\FI-\FP)f^\v_R\|_{\nu}^2\cdot\max_{|v|\leq \f{\lambda}{\v^{1/3}}}(1+|v|)^{3-\kappa}\nonumber\\
	&\leq C_\lambda\|f^\v_{R}\|_{L^2}^2+C\f{\lambda}{\v^2}\|(\FI-\FP)f_R^\v\|_{\nu}^2,
\end{align*}
and
\begin{align*}
	&\int_{\R^3}\int_{|v|\geq\f{\lambda}{\v^{1/3}}}|\nabla_x(\rho,\fu,T)|(1+|v|)^3|f^\v_R|^2dvdx\nonumber\\
	&\dis\leq C\|f^\v_R\|_{L^2}\|h^\v_R\|_{L^\infty}\cdot\Big\{\int_{|v|\geq\f{\lambda}{\v^{\f13}}} (1+|v|)^{6-2\mathfrak{k}}dv\Big\}^{\f12}\nonumber\\
	&\leq C_\lambda \v^{\f13\mathfrak{k}-2}\|f_R^\v\|_{L^2}\|h_R^\v\|_{L^\infty}\leq C_\lambda \v^4\|f_R^\v\|_{L^2}\|h_R^\v\|_{L^\infty}.
\end{align*}
Using Lemma \ref{lem2.10}, one has
\begin{align*}
	&\v^3 \int_{\mathbb{R}^3_+}\intr \frac1{\sqrt{\mu}}Q(\sqrt{\mu}f^\v_R,\sqrt{\mu}f^\v_R) f^{\v}_R dvdx\nonumber\\
	&=\v^3 \int_{\mathbb{R}^3_+}\intr \frac1{\sqrt{\mu}}Q(\sqrt{\mu}f^\v_R,\sqrt{\mu}f^\v_R) \{\mathbf{I-P}\}f^{\v}_R dvdx\nonumber\\
	&\leq \v^3 \|\{\mathbf{I-P}\}f^{\v}_R \|_{\nu} \|h^\v_R\|_{L^\infty} \|f^\v_R\|_{L^2}\nonumber\\
	&\leq \frac{\lambda}{\v^2} \|\{\mathbf{I-P}\}f^{\v}_R \|_{\nu}^2+C_{\lambda} \v^8 \|h^\v_R\|_{L^\infty}^2 \|f^\v_R\|_{L^2}^2.
\end{align*}

From \eqref{H47}, we have
\begin{equation*}
	s_N>\bar{s}_N\geq 2\fb+4+\hat{s}_N,\quad \hat{s}_N\geq 1,
\end{equation*}
which, together with Proposition \ref{prop5.1} and Sobolev embedding theorem, yields  that, for $1\leq i\leq N$ and $t\in[0,\tau^\d]$,
\begin{align}\label{3.10}
	\begin{split}
		&\sum_{k=0}^{2\fb+2}\Big\{\left\|\tw_{\k_i}(v)   \nabla^k_{t,x}f_i(t)\right\|_{L^2_{x,v}}+\left\|\tw_{\k_i}\nabla^k_{t,x} f_i(t)\right\|_{L^\infty_{x,v}}\Big\}\leq C_R(\tau^\d),\\
		&\sum_{k=0}^{\fb+2}\Bigg\{\left\|\fw_{\bar{\k}_i} (1+y)^{\fb+9} \nabla^k_{t,\bar{x}}\bar{f}_i(t)\right\|_{L^2_{\bar{x},v}}+\left\|\fw_{\bar{\k}_i} (1+y)^{\fb+9} \nabla^k_{t,\bar{x}}\bar{f}_i(t)\right\|_{L^\infty_{\bar{x},v}} \Bigg\}\leq C_R(\tau^\d), \\
		&\sum_{k=0,1}\Bigg\{\left\|\fw_{\bar{\k}_i} (1+\eta)^{\zeta_i} \nabla_{t,x_{\sp} }^k\hat{f}_i(t)\right\|_{L^2_{\hat{x},v}}+ \left\|\fw_{\bar{\k}_i} (1+\eta)^{\zeta_i} \nabla_{t,x_{\sp} }^k\hat{f}_i(t)\right\|_{L^\infty_{\hat{x},v}}\Bigg\}\leq C_R(\tau^\d),
	\end{split}
\end{align}
where we have denoted
\begin{align}
	C_R(\tau^\d):&=C\Bigg( \tau^\d, \|(\varphi_0,\Phi_0, \vartheta_0)\|_{H^{s_0}}+\sum_{i=0}^{N}\sum_{\gamma+\beta\leq s_{i}}\|\pa_t^\gamma\nabla_x^{\beta}(\rho_{i}, u_{i},  \theta_{i})(0)\|_{L^2_{x}}\nonumber\\
	&\qquad\qquad\qquad+\sum_{i=0}^{N}\sum_{j=0}^{\bar{s}_{i}} \sum_{j=2\gamma+\beta}   \|\pa_t^{\gamma} \nabla_{\bar{x}}^{\beta} (\bar{u}_{i,\sp},\bar{\theta}_{i})(0)\|^2_{L^2_{l_j^{i}}}
	\Bigg).\nonumber
\end{align}
\

\

Recall $\varpi_{\mathfrak{k}}$ in \eqref{c1.30}. For $1\leq i\leq N$, it is clear that
\begin{align}\label{3.9}
	\begin{split}
		\left|\varpi_{\mathfrak{k}}(v) \frac{\sqrt{\mu_0}}{\sqrt{\mu}} \bar{f}_i(t,x_{\sp} ,y,v)\right|
		&\leq C |w_{\mathfrak{k}}(v) \mu_0^{-\fa_i} \bar{f}_i(t,x_{\sp} ,y,v)|\cdot \frac{\mu_0^{\frac12+\fa_i}}{\mu^{\frac12}}\\
		&\leq C |w_{\mathfrak{k}}(v) \mu_0^{-\fa_i} \bar{f}_i(t,x_{\sp} ,y,v)|\cdot (\mu_{M})^{(\frac12+\fa_i)\alpha-\frac12},\\[2mm]
		\left|w_{\mathfrak{k}}(v) \frac{\sqrt{\mu_0}}{\sqrt{\mu}} \hat{f}_i(t,x_{\sp} ,\eta,v)\right|
		&\leq C |w_{\mathfrak{k}}(v) \mu_0^{-\fa_i} \hat{f}_i(t,x_{\sp} ,y,v)|\cdot (\mu_{M})^{(\frac12+\fa_i)\alpha-\frac12}.
	\end{split}
\end{align}
Taking $0<\frac{1}{2\alpha}(1-\alpha)<\fa_i<\frac12$, we have
$ (\frac12+\fa_i)\alpha-\frac12>0$. Taking $\k_i, \bar{\k}_{i}, \hat{\kappa}_{i} > \mathfrak{k}$,
using \eqref{3.9},  \eqref{ewf}, \eqref{3.10} and  \eqref{3.9}, implies that  the third term on RHS of \eqref{H7} is bounded by
\begin{align}
	&\|\{\mathbf{I-P}\}f^{\v}_R \|_{\nu} \|f^\v_R\|_{L^2} \cdot \sum_{i=1}^N \v^{i-2}  \bigg\{ \|\tw_{\k_i}(v) f_i\|_{L^\infty_{x,v}}+\|\fw_{\bar{\k}_i}   \bar{f}_i\|_{L^\infty_{\bar{x},v}}+\|\fw_{\hat{\k}_i}  \hat{f}_i\|_{L^\infty_{\hat{x},v}}\bigg\}\nonumber\\
	&\leq \frac{C(\tau^{\delta})}{\v}\|\{\mathbf{I-P}\}f^{\v}_R \|_{\nu} \|f^\v_R\|_{L^2}
	\leq \frac{\lambda}{\v^2} \|\{\mathbf{I-P}\}f^{\v}_R \|_{\nu}^2+C_\lambda(\tau^{\delta}) \|f^{\v}_R \|_{L^2}^2.\nonumber
\end{align}

From \eqref{H2} and $\eqref{3.10}_1$, a direct calculation shows that
\begin{align}\label{3.12}
	\left(\int_{\mathbb{R}^3_+}\intr |\frac1{\sqrt{\mu}}R^\v|^2 dvdx\right)^{\frac12}\leq C_R(\tau^\d) \v^{N-6}.
\end{align}
Noting the polynomial space decay estimates on viscous and Knudsen boundary layers in Proposition \ref{prop5.1}, by taking $l_j^{i}, \zeta_i\gg \fb$,  it follows from \eqref{H3}, \eqref{H4} and \eqref{3.10} that
\begin{align}\label{3.14}
	\begin{split}
		\left(\int_{\mathbb{R}^3_+}\intr |\frac1{\sqrt{\mu}}\bar{R}^\v|^2 dvdx\right)^{\frac12}&\leq C_R(\tau^\d) (\v^{N-5.5}+\v^{\fb-4.5}),\\
		\left(\int_{\mathbb{R}^3_+}\intr |\frac1{\sqrt{\mu}}\hat{R}^\v|^2 dvdx\right)^{\frac12}&\leq C_R(\tau^\d) (\v^{N-5}+\v^{\fb-3}).
	\end{split}
\end{align}
Combining \eqref{3.12}-\eqref{3.14} and Cauchy inequality, we conclude that
\begin{equation*}
	\left| \int_{\mathbb{R}^3_+}\intr \left\{\frac1{\sqrt{\mu}}R^\v+\frac1{\sqrt{\mu}}\bar{R}^\v+\frac1{\sqrt{\mu}}\hat{R}^\v\right\} f^\v_R dvdx\right|\leq C_R(\tau^\d) \, (\v^{N-6} +\v^{\fb-5}) \|f^\v_R\|_{L^2}.
\end{equation*}
Hence we conclude \eqref{H6} from above estimates. Therefore the proof of Lemma \ref{lem5.1} is completed. $\hfill\Box$

\

Next we consider weighted $L^\infty$-estimate.
As in \cite{Guo Jang Jiang-1,Guo Jang Jiang}, we denote
\begin{equation}
	L_{M}g=-\frac{1}{\smum}\Big\{Q(\mu,\sqrt{\mu_M} g)+Q(\sqrt{\mu_M}  g,\mu)\Big\}=\nu(\mu)g-\hat{K}g,\nonumber
\end{equation}
where  $\nu(\mu)$ is the collision frequency of $\FL$, and  $\hat{K}g=\hat{K}_{2}g-\hat{K}_{1}g$ with
\begin{align}
	\hat{K}_{1}g=&\intr\ints |(v-u)\cdot \omega|\sqrt{\mu_{M}(u)}\frac{\mu(v)}{\sqrt{\mu_{M}(v)}}g(u)dud\omega,
	\nonumber\\
	\hat{K}_{2}g=&\intr\ints |(v-u)\cdot \omega|\mu(u')\frac{\sqrt{\mu_{M}(v')}}{\sqrt{\mu_{M}(v)}}g(v')dud\omega
	\nonumber\\
	&+\intr\ints |(v-u)\cdot \omega|\mu(v')\frac{\sqrt{\mu_{M}(u')}}{\sqrt{\mu_{M}(v)}}g(u')dud\omega.\nonumber
\end{align}
We still divide $\hat{K}_1,\hat{K}_2$ into $\hat{K}_1^m,\hat{K}_1^c,\hat{K}_2^m,\hat{K}_2^c$ as \eqref{a46}.

\begin{lemma}[\cite{Guo Jang Jiang}]\label{lem5.2-1}
	It holds that
	\begin{align*}
		|\hat{K}^mg(v)|\leq Cm^{3+\kappa}\nu(\mu)\|g\|_{L^\infty},
	\end{align*}
and $\displaystyle \hat{K}^cg(v)=\intr l(v,v')g(v') dv'$ where the kernel $l(v,v')$ satisfies
	\begin{equation*}
		|l(v,v')|\leq C_m\f{\exp(|v-u|^2)}{|v-u|(1+|v|+|u|)^{1-\kappa}}.
	\end{equation*}
\end{lemma}

Denoting $\displaystyle K_{\varpi}g\equiv \varpi_\mathfrak{k} \hat{K}(\frac{g}{\varpi_\mathfrak{k}})$, we deduce from \eqref{H1} and \eqref{c1.31} that
\begin{align}\label{4.9}
	\dis &\partial_t h^\v_R+v\cdot\nabla_x h^\v_R+\frac{\nu(\mu)}{\v^2} h^\v_R-\frac1{\v^2}K_\varpi h^\v_{R} \nonumber\\
	 &=\sum_{i=1}^{N}\v^{i-2}\frac{\varpi_{\mathfrak{k}}(v)}{\sqrt{\mu_M(v)}}\Big\{Q(F_i+\bar{F}_i+\hat{F}_i,\frac{\sqrt{\mu_M} h^\v_R}{\varpi_\mathfrak{k}})
	+Q(\frac{\sqrt{\mu_M} h^\v_R}{\varpi_\mathfrak{k}}, F_i+\bar{F}_i+\hat{F}_i)\Big\}\nonumber\\
	&\quad+\v^{3}\frac{\varpi_\mathfrak{k}}{\sqrt{\mu_M}}Q\Big(\frac{\sqrt{\mu_M} h^\v_R}{\varpi_\mathfrak{k}},\frac{\sqrt{\mu_M} h^\v_R}{\varpi_\mathfrak{k}}\Big)+ \frac{\varpi_\mathfrak{k}}{\sqrt{\mu_M}}\big[R^\v+\bar{R}^\v+\hat{R}^\v\big].
\end{align}

Using Lemma \ref{lem5.2-1}, by similar arguments as in \cite[Lemma 6.3]{GHW} (see also \cite[Lemma 2.2]{Guo Jang Jiang}), we can obtain the following $L^\infty$ estimates. Here we omit the details for simplicity of presentation.
\begin{lemma}\label{lem5.2}
	For $t\in[0,\tau^\d]$, it holds that
	\begin{equation*}
		\sup_{0\leq s\leq t}\|\v^{3}h^\v_R(s)\|_{L^\infty}
		\leq C(t) \{\|\v^{3}h^\v_R(0)\|_{L^\infty}+\v^{N-1}+\v^{\fb}\}+\sup_{0\leq s\leq t}\|f^\v_R(s)\|_{L^2}.
	\end{equation*}
\end{lemma}

\subsection{Proof of Theorem \ref{theorem}}

With Lemma \ref{lem5.1} and Lemma \ref{lem5.2}, one can close the proof by the same arguments as in \cite{Guo Jang Jiang}. We omit the details for simplicity of presentation. Therefore  the proof of Theorem \ref{theorem} is complete. $\hfill\Box$

	\
	
	\noindent{\bf Acknowledgments.} Yong Wang's research is partially supported by National Key R\&D Program of China No. 2021YFA1000800, National Natural Science Foundation of China No. 12022114, 12288201, CAS Project for Young Scientists in Basic Research，Grant No. YSBR-031, and Youth Innovation Promotion Association of the Chinese Academy of Science No. 2019002. We thank Weiqiang Wang for his valuable discuss.
	\
	
	\noindent{\bf Conflict of interest.} The authors declare that they have no conflict of interest.
	

\begin{thebibliography}{99}
		
		
		\bibitem{BCN-1986} C. Bardos, R.E. Caflisch, B. Nicolaenko, The Milne and Kramers problems for the Boltzmann equation of a hard sphere gas, Comm. Pure Appl. Math. 39 (1986) 323-352.
		
		\bibitem{Bardos} C. Bardos, F. Golse, C.D.  Levermore,  Fluid dynamic limits of kinetic equations. I. Formal derivations. J. Statist. Phys. 63 (1991), no. 1-2, 323-344.
		
		\bibitem{Bardos-2} C. Bardos, F. Golse, C.D.  Levermore, Fluid dynamic limits of kinetic equations. II. Convergence proofs for the Boltzmann equation. Comm. Pure Appl. Math. 46 (1993), no. 5, 667-753.
		
		\bibitem{BG} N. Bernhoff,  F. Golse,  On the boundary layer equations with phase transition in the kinetic theory of gases. Arch. Ration. Mech. Anal. 240 (2021), no. 1, 51-98.
		
		\bibitem{Bardos-Ukai} C. Bardos, S. Ukai,  The classical incompressible Navier-Stokes limit of the Boltzmann equation. Math. Models Methods Appl. Sci. 1 (1991), no. 2, 235-257.
		
		\bibitem{Boltzmann} Boltzmann, L. Weitere Studien \"{u}ber das W\"armegleichgewicht unter Gasmolek\"ulen. Sitzungs. Akad. Wiss. Wien 66 (1872), 275-370; translated as: Further studies on the thermal equilibrium of gas molecules. Kinetic theory, vol. 2, 88-174. Pergamon, London, 1966.
		
		\bibitem{Caflish} R.E.Caflisch, The fluid dynamic limit of the nonlinear Boltzmann equation. Comm.Pure Appl. Math. 33(1980), No.5, 651-666.
		
		\bibitem{CLY} C.C. Chen, T.P. Liu, T. Yang, Existence of boundary layer solutions to the Boltzmann equation, Analysis and Applications, Vol.2, No.4 (2004) 337-363.
		
		
		\bibitem{Duan} Duan, R.J., Huang, F.M., Yang T., Wang Y., Global well-posedness of the Boltzmann equation with large amplitude initial data, Arch. Rational Mech. Anal. 225 (2017), 375-424.
		
		\bibitem{DHWZ} R.J. Duan, F.M. Huang, Y. Wang, Z. Zhang, Effects of soft interaction and non-isothermal boundary upon long-time dynamics of rarefied gas, Arch. Rational Mech. Anal. 234 (2019) 925-1006.
		
		\bibitem{GL}R.T. Glassey, The Cauchy Problem in Kinetic Theory, Society for Industrial and Applied Mathematics (SIAM),Philadelphia, 1996.
		
		\bibitem{Golse-Saint-Raymond} F. Golse, L.  Saint-Raymond, The Navier-Stokes limit of the Boltzmann equation for bounded collision kernels. Invent. Math. 155 (2004), no. 1, 81-161.
		
		\bibitem{GPS-1988} F. Golse, B. Perthame, C.  Sulem, On a boundary layer problem for the nonlinear Boltzmann equation. Arch. Rational Mech. Anal. 103 (1988), no. 1, 81-96.
		
		\bibitem{Grad} H. Grad, Asymptotic equivalence of the Navier-Stokes and nonlinear Boltzmann equations, in Proc. Sympos. Applied Mathematics, Vol. 17 (Amer. Math. Soc., 1965), pp. 154-183.
		
		
		
		\bibitem{Guo2006} Y. Guo, Boltzmann diffusive limit beyond the Navier-Stokes approximation. Comm. Pure Appl. Math. 59 (2006), no. 5, 626-687.
		
		\bibitem{Guo2010} Y. Guo, Decay and continuity of the Boltzmann equation in bounded domains. Arch. Ration. Mech. Anal. 197
		(2010), no. 3, 713-809.
		
		\bibitem{Guo Jang Jiang-1} Y. Guo, J. Jang, N. Jiang, Local Hilbert expansion for the Boltzmann equation. Kinet. Relat. Models. 2 (2009), no. 1, 205-214.
		
		\bibitem{Guo Jang Jiang} Y. Guo, J. Jang, N. Jiang, Acoustic limit for the Boltzmann equation in optimal scaling. Comm. Pure Appl. Math. 63 (2010), no. 3, 337-361.
		
		
		\bibitem{GKTT} Y. Guo,  C. Kim,  D. Tonon,  A. Trescases,  Regularity of the Boltzmann equation in convex domains,  Invent. Math.,  207 (2017), 115-290.
		
		\bibitem{GHW} Y. Guo, F.M. Huang, Y. Wang,  Hilbert expansion of the  Boltzman equation with specular boundary condition in half-space, Arch. Rational. Mech. Anal.241(2021), 231-309.
		
		
		\bibitem{Guo-Liu} Y. Guo, S.Q. Liu, Incompressible hydrodynamic approximation with viscous heating to the Boltzmann equation, Math. Models Methods Appl. Sci., 27  (2017), no. 12, p. 2261-2296.
		
		
		
		\bibitem{Jiang-Wang} F.M. Huang,  Z.H. Jiang, Y. Wang, Boundary layer solution of the Boltzmann equation for specular boundary condition. Acta Mathematicae Applicatae Sinica,
		Vol. 39, No. 1 (2023) 65-94.
		
		\bibitem{HWWX} F.M. Huang, W.Q. Wang, Y. Wang, F. Xiao, Hilbert expansion of the Boltzmann equation in the incompressible Euler level in a channel, Sci. China Math., doi.org/10.1007/s11425-016-5135-4.
		
		 \bibitem{HW}  F.M. Huang, Y. Wang, Boundary layer solution of the Boltzmann equation for diffusive reflection boundary conditions in half-space. SIAM J. Math. Anal. 54 (2022), no. 3, 3480-3534.
		
		\bibitem{Huang-Wang-Yang} F. M. Huang, Y. Wang,  T. Yang, Hydrodynamic limit of the Boltzmann equation with contact discontinuities, Comm. Math. Phys., 295 (2010), pp. 293-326.
		
		\bibitem{Huang-Wang-Yang-1} F. M. Huang, Y. Wang,  T. Yang, Fluid Dynamic Limit to the Riemann Solutions of Euler Equations: I. Superposition of rarefaction waves and contact discontinuity, Kinet. Relat. Models, 3 (2010), pp. 685-728.
		
		\bibitem{Huang-Wang-Wang-Yang} F. M. Huang, Y. Wang, Y. Wang, T. Yang, The limit of the Boltzmann equation to the Euler equations for Riemann problems. SIAM J. Math. Anal. 45 (2013), no. 3, 1741-1811.
		
		\bibitem{Jang-Kim} J. Jang, C. Kim, Incompressible Euler limit from Boltzmann equation with Diffuse Boundary Condition for Analytic data. Ann. PDE 7 (2021), no.2, Paper No.22, 103 pp.
		
		\bibitem{Jiang-Masmoudi} N. Jiang, N. Masmoudi, Boundary layers and incompressible Navier-Stokes-Fourier limit of the Boltzmann equation in bounded domain I. Comm. Pure Appl. Math. 70 (2017), no. 1, 90-171.
		
		\bibitem{Jiang-Luo-Tang} N. Jiang, Y.L. Luo, S.J. Tang, Grad-Caflisch type decay estimates of pseudo-inverse of linearized Boltzmann operator and application to Hilbert expansion of compressible Euler scaling. arXiv:2206.02677.
		
		\bibitem{Jiang-Luo-Tang-1} N. Jiang, Y.L. Luo, S.J. Tang, Compressible Euler limit from Boltzmann equation with Maxwell reflection boundary condition in half-space, arXiv:2101.11199.
		
		\bibitem{Jiang-Luo-Tang-2} N. Jiang, Y.L. Luo, S.J. Tang, Compressible Euler limit from Boltzmann equation with complete diffusive boundary condition in half-space, arXiv:2104.11964 .
		
		\bibitem{Lachowicz} M. Lachowicz, On the initial layer and the existence theorem for the nonlinear Boltzmann equation. Math. Methods Appl. Sci. 9 (1987), no. 3, 342-366.
		
		\bibitem{Liu-Yang} S.Q.Liu, X.F. Yang, The Initial Boundary Value Problem for the
		Boltzmann Equation with Soft Potential, Arch. Rational Mech. Anal. 223 (2017) 463-541.
		
		\bibitem{Masmoudi-Raymond}  N. Masmoudi, L. Saint-Raymond, From the Boltzmann equation to the Stokes-Fourier system in a bounded domain.
		Comm. Pure Appl. Math. 56 (2003), no. 9, 1263-1293.
		
		\bibitem{Maxwell} J.C. Maxwell,  On the dynamical theory of gases. Philos. Trans. Roy. Soc. London Ser. A 157 (1867), 49-88. Reprinted in The scientific letters and papers of James Clerk Maxwell, vol. II, 1862-1873, 26-78. Dover, New York, 1965.
		
		\bibitem{Guo2} Robert M. Strain, Y. Guo, Exponential decay for soft potentials near Maxwellian, Arch. Rational Mech. Anal. 187 (2008) 287-339.
		
		\bibitem{Nishida} T. Nishida, Fluid dynamical limit of the nonlinear Boltzmann equation to the level of the compressible Euler equation. Comm. Math. Phys. 61 (1978), no. 2, 119-148.		
		
		\bibitem{Sone-2002} Y. Sone, Kinetic Theory and Fluid Dynamics, Birkhäuser Boston, Inc., Boston, MA, 2002.
		
		\bibitem{Sone-2007} Y. Sone, Molecular Gas Dynamics. Theory, Techniques, and Applications, Birkhäuser Boston, Inc., Boston, MA, 2007.
		
		\bibitem{Ukai-Asano} S. Ukai, K. Asano, The Euler limit and the initial layer of the nonlinear Boltzmann equation, Hokkaido Math. J., 12 (1983), pp. 303-324.
		
		\bibitem{UYY} S. Ukai, T. Yang,  S.H. Yu,  Nonlinear boundary layers of the Boltzmann equation. I. Existence. Comm. Math. Phys. 236 (2003), no. 3, 373-393.
		
		\bibitem{WWZ} G.F. Wang, Y. Wang, J.W. Zhou,  Hydrodynamic limit of the Boltzmann equation to the planar rarefaction wave in three dimensional space. SIAM J. Math. Anal. 53 (2021), no. 4, 4692-4726.
		
		\bibitem{Wu} L. Wu, Hydrodynamic limit with geometric correction of stationary Boltzmann equation, J. Differential Equations 260 (2016) 7152-7249.
		
		\bibitem{WYY} W.K. Wang, T. Yang, X.F. Yang, Nonlinear stability of boundary layers of the Boltzmann equation for cutoff hard potentials, J. Math. Phys., 47, 083301(2006).
		
		\bibitem{WYY2} W.K. Wang, T. Yang, X.F. Yang, Existence of boundary layers to the Boltzmann equation with cutoff soft potentials, J. Math. Phys., 48, 073304(2007).
		
		\bibitem{XZ} Z.P. Xin, H.H. Zeng,  Convergence to the rarefaction waves for the nonlinear Boltzmann equation and compressible Navier-Stokes equations. J.Differ. Equ. 249 (2010), 827-871.
		
		\bibitem{Yu} S.H. Yu,  Hydrodynamic limits with shockwaves of theBoltzmann equations. Commun. Pure Appl. Math 58 (2005), 409-443.
	\end{thebibliography}
\end{document}